\documentclass[11pt]{article}

\usepackage{xr-hyper}

\usepackage{./src/preamble} 
\usepackage{./src/macros} 
\usepackage{authblk}
\usepackage[normalem]{ulem}
\usepackage{xcolor}

\usepackage{euscript}
\usepackage{algorithm} 
\usepackage{algpseudocodex}
\usepackage[nodisplayskipstretch]{setspace}

\bibstyle{sort&compress}
\PassOptionsToPackage{sort&compress,round}{natbib}
\definecolor{nb}{RGB}{0,75,200}
\hypersetup{colorlinks, linkcolor=nb, citecolor=nb, urlcolor=nb}
\allowdisplaybreaks
\crefname{equation}{}{}
\crefname{figure}{Figure}{Figure}

\usepackage[most]{tcolorbox}
\newcommand{\cbox}[2]{%
\begin{tcolorbox}[
    enhanced,            
    breakable,           
    colback=#1,          
    colframe=black!90,   
    arc=3mm,             
    boxrule=0.5pt,       
    left=5pt,right=5pt,  
    top=5pt,bottom=5pt
    ]
    {#2}%
\end{tcolorbox}%
}

\usepackage{multirow}

\title{Confidence Sets for Multidimensional Scaling}

\author[1]{Siddharth Vishwanath} 
\author[1,2]{Ery Arias-Castro}
\affil[1]{\small Department of Mathematics, University of California, San Diego} 
\affil[2]{\small Halıcıoğlu Data Science Institute, University of California, San Diego}
\date{}

\begin{document}
\maketitle
\vspace*{-3em}

\begin{abstract}
We develop a formal statistical framework for classical multidimensional scaling (CMDS) applied to noisy dissimilarity data. We establish  distributional convergence results for the embeddings produced by CMDS for various noise models, which enable the construction of \emph{bona~fide} uniform confidence sets for the latent configuration, up to rigid transformations. We further propose bootstrap procedures for constructing these confidence sets and provide theoretical guarantees for their validity. We find that the multiplier bootstrap adapts automatically to heteroscedastic noise such as multiplicative noise, while the empirical bootstrap seems to \mbox{require} homoscedasticity. Either form of bootstrap, when valid, is shown to substantially improve finite-sample accuracy. The empirical performance of the proposed methods is demonstrated through numerical experiments.
\end{abstract}


\section{Introduction}
\label{sec:introduction}

Multidimensional scaling (MDS) is an essential tool in multivariate analysis, and underpins a broad class of unsupervised learning and linear/non-linear dimension reduction techniques. The objective of MDS is to embed a set of $n$ items into a low-dimensional Euclidean space given only an $n \times n$ matrix $\Del = (\del_{ij})$ of pairwise dissimilarities between the items. Specifically, {given an embedding dimension $p$ (often $p = 2$ for visualization purposes),} the goal is to find a configuration $\hX \in \Rnp$ of $n$ points embedded in $\Rp$ such that the pairwise squared Euclidean distances between the points in $\hX$ reproduce the original dissimilarities.

Although MDS is sometimes used as a dimension reduction method similar to principal component analysis (PCA), its scope is more general. PCA operates directly on feature vectors from a data matrix $X \in \R^{n \times q}$ to produce a lower-dimensional representation in $\R^p$ for $p < q$. In contrast, MDS only requires access to a pairwise dissimilarity matrix, $\Del$, between the $n$ items. This allows MDS to be applied in settings where the original data are unavailable, or where relational or proximity data is more natural and meaningful. Such situations naturally arise in many applications, e.g., survey data in psychology, spatial capture-recapture data in ecology, morphological and physiological dissimilarities in biology, and sensor network data in wireless communication, to name a few. For a comprehensive overview of MDS and its applications, see~\citep{borg2005modern,young2013multidimensional}. 

Despite its long history and its broad range of applications, the statistical treatment of MDS has remained relatively underdeveloped. As a result, most applications of MDS have been exploratory in nature---serving primarily as graphical tools for data visualization. The embedded points are often interpreted directly, without any adjustment for the uncertainty arising from sampling variation or measurement noise. The lack of a formal statistical framework in CMDS was recognized {by \cite{ramsay1982some}, who noted:} 
``\emph{Implicit in almost all data analyses is some statement about the manner in which the observation varies about its fitted value.}''
The absence of inferential tools can be problematic, especially in applications where geometric and topological properties of the embedded configurations are used to draw inferences.

We focus on the \textit{classical} multidimensional scaling (CMDS) algorithm, which dates back to the foundational work of \cite{young1938discussion} and later formalized by \cite{torgerson1952multidimensional} and \cite{gower1966some}; CMDS is not only the analog of PCA, but is as central to MDS as PCA is to dimensionality reduction. In this work, we aim to place CMDS within a formal statistical framework and develop methods for constructing \emph{uniform confidence sets} for the latent configuration underlying the observed dissimilarities. These confidence sets account for sampling noise and provide valid simultaneous coverage for all points in the configuration, modulo rigid transformations, necessarily.

\subsection{Contributions}

\begin{algorithm}[t]
    \small
    \caption{\small Classical Multidimensional Scaling (CMDS)}
    \label{alg:cmds}
    \begin{algorithmic}[1]
        \Require Dissimilarity matrix $D \in \R^{n \times n}$, embedding dimension $p$.
        \State Compute $\Dc = -\half HDH$ where $H = I - \frac{1}{n}\onev\onev\tr$ \Comment{Double-centering}
        \State Compute the $p$-largest eigenvalues of $\Dc$, $\h\lambda_1, \cdots, \h\lambda_p$, and corresponding eigenvectors $\h{u}_1, \dots, \h{u}_p$
        \State Set $\hL = \text{diag}(\h\lambda_1, \dots, \h\lambda_p) \in \Rpp$ and $\hU = [\h{u}_1\tr, \dots, \h{u}_p\tr] \in \Rnp$
        \\[0.5em]
        \Return Embedding $\hX = \hU \hL^{1/2}$
    \end{algorithmic}
\end{algorithm}

We place ourselves in the \emph{noisy realizable setting} where the observed dissimilarities ${D = (d_{ij}) \in \Rnn}$ are noisy versions of true squared Euclidean distances $\delta_{ij} = \|x_i - x_j\|^2$ between unknown latent points $x_1, \dots, x_n \in \R^p$, i.e., 
\begin{align}
    d_{ij} = \delta_{ij} + \eps_{ij} \qq{for all } i < j \in [n],\label{noisy-setting}
\end{align}
where $\Eps = (\eps_{ij})$ is a symmetric and hollow\footnote{The assumption that $\Eps$ is hollow ensures that $d_{ii} = 0$ for all $i \in [n]$, which is natural for dissimilarity data.} random noise matrix. Letting $X \in \Rnp$ denote the latent configuration with rows $x_1, \dots, x_n$, and $\Delta(X) = (\delta_{ij})$, we write $D = \Delta(X) + \Eps$.  Given $D$, the CMDS algorithm returns an embedding $\hX = \mds(D, p)$ via \cref{alg:cmds}. From a statistical estimation standpoint, the configuration $X \in \Rnp$ constitutes the unknown parameters of interest, and $\hX$ is an estimator of these parameters. Because $\Delta(X) = \Delta(g(X))$ for any rigid transformation $g \in \euc{p}$, the configuration $X$ is only identifiable up to such transformations.

In order to quantify the uncertainty in the embedding $\hX$, we construct \textit{uniform confidence sets} for $\X$. For a level $\alpha \in (0, 1)$, the set $\conf_\alpha(D) := \prod_{i=1}^n \conf_{\alpha, i}(D) \subset (\Rp)^n$ is a uniform $(1-\alpha)$-confidence set for the configuration $X$ (up to rigid transformations) if
\begin{align}
    \pr\qty\big( \exists\; g \in \euc{p}\;: g(x_i) \in \conf_{\alpha, i}(D) \;\;\forall i \in [n] ) \ge 1 - \alpha.
\end{align}
In other words, with probability at least $1-\alpha$, there exists a (data dependent) rigid transformation $g \in \euc{p}$ such that each transformed latent point $g(x_i)$ is contained in the corresponding confidence region $\conf_{\alpha, i}(D)$ for all $i \in [n]$ simultaneously.

With this background, our main contributions are summarized as follows:
\begin{itemize}%
    \item We prove that, under mild conditions, the maximum deviation of the embedding $\hat{X}$ from the true configuration $X$ (up to a rigid transformation and after appropriate normalization) converges to the Gumbel distribution (\cref{thm:main}). This leads to the construction of plug-in confidence sets for the true configuration $X$, up to rigid transformations, with {bona fide} uniform coverage guarantees (\cref{cor:confidence-pivotal}).
    
    \item {While the plug-in approach guarantees valid inference, it often suffers from limited finite-sample accuracy.} To this end, we propose a \emph{multiplier bootstrap} procedure for constructing these confidence sets. We establish the validity of the bootstrap procedure and show that it achieves {much} improved finite-sample accuracy (\cref{thm:multiplier-bootstrap}). The resulting confidence sets are adaptive to non-identically distributed, and in particular, heteroscedastic noise (\cref{cor:confidence-multiplier-bootstrap}).
    
    \item For the special case of additive \iid{} noise, {which is necessarily homoscedastic,} we establish a similar finite-sample convergence result for the empirical bootstrap (a.k.a. \textit{nonparametric} or \textit{Efron's} bootstrap) procedure and show that the resulting confidence sets also enjoy the {much} improved accuracy (\cref{thm:empirical-bootstrap}).
\end{itemize}

\begin{figure}[t!]
    \centering
    \includegraphics[width=0.65\linewidth]{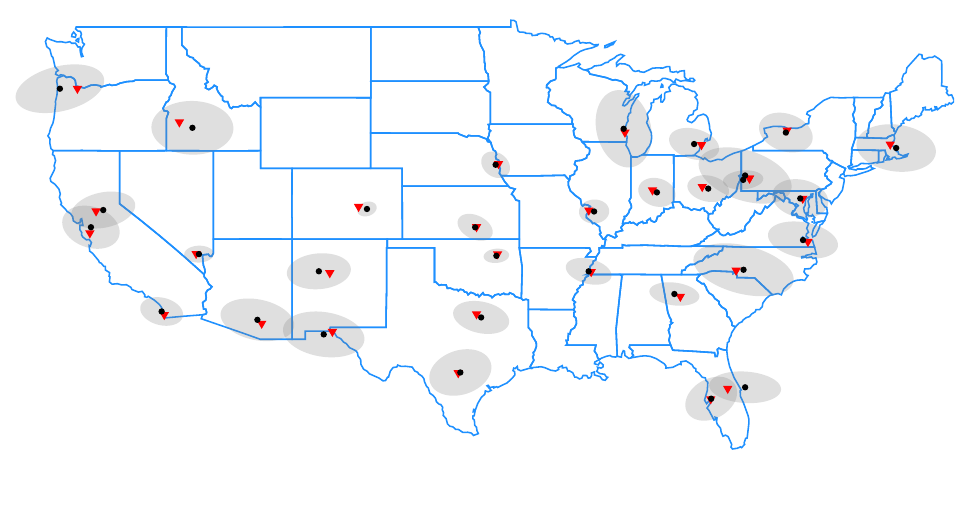}
    \vspace*{-1em}
    \caption{
        Classical multidimensional scaling (MDS) embedding of the noisy pairwise distance between 30 U.S. cities with multiplicative noise. The red points denote the true city locations $\X$, while the black points represent the MDS estimates $\hX$ after Procrustes alignment. The gray ellipsoids constitute the $90\%$ confidence set, which guarantees that each true location lies within its corresponding ellipsoid with probability $1-\alpha=0.9$.
    }
    \label{fig:confidence-set-example}
\end{figure}

We illustrate the performance of our proposed bootstrap procedures through simulations and numerical experiments in \cref{sec:experiments}. For example, \cref{fig:confidence-set-example} shows the typical output of the multiplier bootstrap on the noisy pairwise distances between 30 U.S. cities with multiplicative noise.  The gray ellipsoids representing a 90\% confidence set capture the true latent positions (in red) with high fidelity.

\subsection{Related Work}

Various inference strategies for MDS have been proposed in prior work. \cite{ramsay1977maximum,ramsay1978confidence,ramsay1982some} was among the first to address statistical inference in this context. In particular, he introduced a maximum likelihood framework assuming that $d_{ij} \sim \log{N}(\del_{ij}, \sigma_{ij}^2)$ where the observed dissimilarities are log-normally distributed around the true values. In this framework, it is assumed that multiple independent replicates of the dissimilarity matrix are available, i.e., $D^{(1)}, \ldots, D^{(N)} \in \Rnn$. 
Notably, this assumption circumvents the main difficulty of performing inference from a single dissimilarity matrix where the number of parameters to be estimated ($np$ in total for $X \in \Rnp$) grows with the sample size $n$. {This last setup is the one we consider, and we develop our inference without imposing any parametric assumptions on the noise distribution.} 

To the best of our knowledge, the setup we consider here was first formally studied by \cite{li2020central}, who consider three specific noise models in the noisy realizable setting. For each row $\hx_i \in \Rp$ of the MDS embedding, they derive a central limit theorem, establishing that $\hx_i$ is asymptotically normally distributed around its latent counterpart, $g(x_i)$, after a suitable rigid transformation. Our work differs in three key aspects. First, our emphasis is on \textit{uniform} confidence sets, which guarantee simultaneous coverage for all points in the configuration, whereas the results in \citep{li2020central} are pointwise. This strengthening is non-trivial, and requires new technical machinery from extreme value theory. Second, our results are quantitative and non-asymptotic, providing explicit rates of convergence to the limiting distributions in the Kolmogorov-Smirnov metric. Lastly, we work in the more general setting of  heteroscedastic noise models studied in \citep{vishwanath2025minimax}, which includes the noise models considered in \citep{li2020central}.
{We note that we use the finite-sample error bounds established in \citep{vishwanath2025minimax} to derive distributional convergence results.}

On the application side, several studies have proposed practical methods for constructing confidence regions for the output of MDS more generally. \cite{jacoby2014bootstrap} were the first to investigate the use of bootstrap resampling for MDS. Their method relies on generating bootstrap replicates by resampling the rows of $\X$, from which confidence ellipsoids are constructed. Their approach, while only applicable to the case where $X$ is available, is primarily ad-hoc and provides no formal coverage guarantees from a statistical standpoint {and seems to yield anti-conservative confidence sets. 
We note that our bootstrap procedures are fundamentally different operationally, and come with theoretical guarantees.} 
In a different direction, \cite{de2019pseudo} avoids the need for resampling entirely by constructing \textit{pseudo-confidence regions} using the Hessian of the MDS stress function---a tool which is typically used in stability analyses. This approach, however, is also ad-hoc, {and is developed without specifying the type of noise model being considered}. 
In the 1980s, \cite{de1986special} used the jackknife (equivalently, the leave-one-out) method in order to assess the stability of the MDS solutions. 
Interestingly, in the same way that the jackknife can be viewed as a precursor to the bootstrap \citep{wu1986jackknife}, the bootstrap approach we consider here is perhaps most similar in spirit to {the jackknife approach of} \cite{de1986special};  in particular, we generate replicates from the residuals as opposed to from the observed dissimilarities itself. 
Finally, \cite{nikitas2023assessing} conduct a comparative study of various methods for constructing confidence ellipsoids for MDS, including the methods discussed above. They employ qualitative criteria based on a ``visual inspection of plots'' and quantitative criteria based on examining cluster probabilities and stability measures resulting from the effect of adding a constant value to all dissimilarities \cite[Section~5.2]{nikitas2023assessing}. Their study, however, does not examine the actual coverage guarantees for the methods they consider.

\textbf{Organization.}\quad 
In \cref{sec:background}, we introduce the {setting, including a description of CMDS and a definition of the noise models that we consider}. We present our main distributional convergence results for CMDS embeddings in \cref{sec:convergence}. 
{In \cref{sec:bootstrap}, we introduce and provide theoretical guarantees for the bootstrap: the multiplier bootstrap in \cref{sec:multiplier-bootstrap} and the empirical bootstrap in \cref{sec:empirical-bootstrap}, the latter being analyzed under the special case of \iid{} noise}.  We illustrate the performance of our methods through numerical experiments in \cref{sec:experiments}. {\cref{sec:discussion} contains a brief discussion of our results.} The proofs for the main results are deferred to \cref{sec:proofs}, and the more technical details are relegated to \cref{sec:toolkit,sec:proof-auxiliary}.

\textbf{Notation.}\quad For $\onev \in \Rn$ and $J = \onev\onev\tr$ (the matrix of all 1s),   $H = I - J/n$ denotes the centering matrix. For $x \in \Rp$, $\norm{x}$ denotes the Euclidean norm (i.e., the $\ell_2$-norm). For $A \in \R^{m \times k}$, $\opnorm{A}, \ttinf{A}$ and $\frobenius{A}$ denote the $\ell_2$-operator norm, the $\ell\ttinft$-operator norm and the Frobenius norm of $A$, respectively. $\orth{p}$ denotes the group of $p\times p$ orthogonal matrices, and $\euc{p}$ the group of rigid transformations on $\Rp$.

We also use standard asymptotic notation: we write $a_n = O(b_n)$ (equiv. ${a_n \lesssim b_n}$) for two sequences $a_n, b_n$, if there exists $C > 0$ such that $\abs{a_n} \le C \abs{b_n}$ for sufficiently large $n$, and $a_n \asymp b_n$ if $a_n \lesssim b_n$ and $b_n \lesssim a_n$. Similarly, $a_n = o(b_n)$ if $\lim_n\abs{a_n/b_n} = 0$ and $a_n \sim b_n$ if ${\lim_n\abs{a_n/b_n-1} = o(1).}$ 
For a sequence of random variables $\xi_n$, we write $\xi_n = \Op(a_n)$ if there exists $C > 0$ such that $\pr(\abs{\xi_n/a_n} > C) \le 1/n$ for all $n > N_C$, and $\xi_n = o_p(1)$ if $\lim_n\pr(\abs{\xi_n/a_n} > C) = 0$ for all $C > 0$. 

For a real valued random variable $\xi$, $\norm{\xi}_{\psi_1}$ and $\norm{\xi}_{\psi_2}$ denote its sub-exponential and sub-Gaussian norms \citep[Chapter~2]{vershynin2018high}. For a random vector $\zeta \in \R^k$, $\norm{\zeta}_{\psi_p}\!\!:=\!\!\max_{\norm{x}=1}\norm{x\tr\zeta}_{\psi_p}$. A summary of additional notation introduced in the text is collected in \cref{tab:notation}.

\section{Background}
\label{sec:background}

In the realizable setting, the matrix $\Del$ is assumed to be a Euclidean dissimilarity matrix, i.e., $\del_{ij} = \norm{x_i - x_j}^2$, or, equivalently, in matrix form, 
$$
\Del = \diag(XX\tr)\onev\tr + \onev\diag(XX\tr)\tr - 2XX\tr,
$$ 
where $X \in \Rnp$ is the latent configuration. Throughout, $p < n$ is assumed to be fixed and known. 

A classical result due to \cite{schoenberg1935remarks} {(essentially in parallel with \citealp{young1938discussion})} establishes that $\Del$ is a Euclidean dissimilarity matrix if and only if the double-centering transformation $\Delc=-\half H\Del H$ appearing in line~1 of \cref{alg:cmds} is positive semi-definite. In fact, since $H\onev = \onev\tr H = \zerov$, it is easy to see that $\Delc = (H\X)(H\X)\tr$ corresponds to the Gram matrix of $HX$. Moreover, since we restrict our attention to the equivalence class of configurations up to rigid transformations, without loss of generality, we assume that the latent configuration $X$ is centered, i.e., $\onev\tr\X = 0$, from which it follows that $\Delc = XX\tr$.

Let the reduced rank-$p$ singular value decomposition of $X$ be given by
$$
X = U \L^{1/2}Q,
$$  
where $Q \in \orth{p}$, $\L = \diag(\lambda_1, \dots, \lambda_p)$, and $U {= [u_1 \cdots u_n]^\top} \in \Rnp$ satisfying $U\tr U = I$; the Gram matrix and scatter matrix of $X$ are, respectively, given by $XX\tr = U\L U\tr$ {and} $X\tr X = Q\tr\L Q.$
From lines~2~and~3 of \cref{alg:cmds}, it follows that classical multidimensional scaling with $\Del$ as input results in $\mds(\Del, p) = U\Lambda^{1/2}$ as the output. Therefore, the rotation $Q \in \orth{p}$ perfectly aligns $X$ to $U\L^{1/2}$ via the identity $XQ\tr = \hX$.

For the \textit{noisy realizable setting} in \cref{noisy-setting}, we are given $D = \Del + \Eps$, where $\Eps = (\eps_{ij})$ is a symmetric and hollow random matrix. Some examples of noise models which fit into this framework include: the \textit{additive} noise model, the \textit{multiplicative} noise model, and the \textit{log-normal} noise model,
\begin{align}
    d_{ij} = \del_{ij} + \xi_{ij}\qc{}\qquad
    d_{ij} = \del_{ij}(1 + \xi_{ij}),\qq{and}\qquad{}
    \log{d_{ij}} = \log{\del_{ij}} + \xi_{ij} ,\label{eq:noise-models}
\end{align}
where $(\xi_{ij})$ is an $n \times n$ symmetric and hollow random matrix with \iid{} entries. 

The resulting noise matrices $\Eps$, respectively, have entries:
\begin{align}
    \eps_{ij} = \xi_{ij}\qc{}\qquad
    \eps_{ij} = \del_{ij}\xi_{ij},\qq{and}\qquad{}
    \eps_{ij} = \del_{ij}(\exp(\xi_{ij}) - 1).
\end{align}
See Table~1 of \cite{vishwanath2025minimax} for other examples of noise models that fall within this framework. Let $\hU\hL\hU\tr$ be the rank-$p$ spectral decomposition of $\Dc = -\half HDH$. Then, the output of \cref{alg:cmds} applied to $D$ results in 
\begin{align}
    \hX = \mds(D, p) = \hU\hL^{1/2} \in \Rnp.
\end{align}
Unlike the noiseless case, in general, $\hX$ cannot be perfectly aligned to $\X$. A candidate for the optimal rigid transformation is obtained by solving the orthogonal Procrustes problem: 
\begin{align}
    \hQ = \argmin_{Q \in \orth{p}}\norm{\hU - UQ}^2_F.\label{eq:procrustes}
\end{align}
The matrix $\hQ\in \orth{p}$ solving \cref{eq:procrustes} admits a closed form solution based on the singular value decomposition of $\hU\tr U$. The resulting rigid transformation, $\hg$, aligning $X$ to $\hX$, is given by 
\begin{align}
    \hg(x) = \hP x \qq{where} \hP = \hQ\tr Q.\label{eq:hg}
\end{align}
\begin{remark}\label{rem:centering}
    Since $X$ is assumed to be centered, the optimal rigid transformation $\hg$ only has a rotation component and no  translation component; thus, $\hg(X) = X\hP\tr$ is always centered. 
    
    The map $\hg(X) = X\hP\tr = XQ\tr\hQ$ is a composition of two transformations: 
    \textup{(i)} $XQ\tr$ aligns $X$ to $U\Lambda^{1/2}$ as seen in the noiseless case, and \textup{(ii)} $(XQ\tr)\hQ$ then aligns $U\Lambda^{1/2}$ to $\hX$ via \cref{eq:procrustes}.
\end{remark}

Our main results {are based on} the following assumptions on the configuration $X$ and the noise~$\Eps$.

\cbox{black!5}{%
\begin{enumerate}[label=(\textbf{A}$_{\arabic*}$), ref=\textup{(\textbf{A}$_{\arabic*}$)}]
    \item\label{assumption:compact} For $\Rx > 0$ and $\kappa > 1$, the \textit{centered} configuration matrix $X = U\L^{1/2}Q$ is such that
    \begin{align}
        \ttinf{X} \le \Rx \qq{and} \frac{n}{\kappa^2} \le \lambda_p < \dots < \lambda_1 \le \kappa^2n.
    \end{align}
    \item\label{assumption:noise} The random matrix $\Eps = (\eps_{ij}) \in \Rnn$ is symmetric, hollow, and satisfies the following:
    \begin{enumerate}[label=(\roman*), ref=\ref{assumption:noise}\,\textup{(\roman*)}]
        \item\label{noise-2} {For $\Msigma > 0$, the $\eps_{ij}$ are uniformly $\Msigma$-sub-Exponential, i.e.,}
$$
\displaystyle\max_{i < j}\norm{\eps_{ij}}_{\psi_1} \le \Msigma.
$$
        \item\label{noise-1} $\qty{\eps_{ij}: i < j}$ are independent with $\E(\eps_{ij}) = 0$ and $\Var(\eps_{ij}) = \sigma_{ij}^2$.
        \item\label{noise-3} {For $\msigma > 0$,}
        \begin{align}
            \sum_{\!{\{k \in [n] : \sigma_{ik}^2 > 0\}}\!\!\!} \sigma_{ik}^2 u_k u_k\tr \succcurlyeq \msigma^2 I_p \qq{for all } i \in [n].\label{eq:noise-variance}
        \end{align}
    \end{enumerate}
\end{enumerate}
}

We make a few remarks about these assumptions. First, we note that \ref{assumption:compact} is standard in recent work on CMDS \citep{arias2020perturbation,li2020central,little2023analysis,vishwanath2025minimax}. In particular, {$\ttinf{X} = \max_i\norm{x_i}\le \Rx$ means that the configuration remains compactly supported, and the lower bound on $\lambda_p$ ensures that the configuration remains quantitatively full-dimensional, i.e., the point cloud $\{x_1, \dots, x_p\}$ spans the whole space $\R^p$ and does not become `infinitesimally thin' in the asymptotic limit}. 

\begin{remark}\label{rem:iid}
From Lemma~1 of \cite{vishwanath2025minimax} it follows that a random design where $x_1, \dots, x_n$ are generated \iid{} from some probability distribution $F$ supported on $\Rp$, satisfies \ref{assumption:compact} with high probability (as $n\to\infty$){, up to $o(1)$ additive terms in the constants,} when
    \begin{align}
        \diam(\supp(F)) \le \Rx \qq{and} \kappa^{-2} I_p \preccurlyeq \cov(F) \preccurlyeq \kappa^2 I_p.\label{eq:iid-assumption}
    \end{align}
\end{remark}

The assumptions on the noise $(\eps_{ij})$ are somewhat different from those in \citep{vishwanath2025minimax}. The sub-exponential assumption in \ref{noise-2} is an artefact of our proofs. This assumption can, in principle, be relaxed to requiring that $\E\abs{\eps_{ij}}^4 \le \Msigma^4$ at the price of more tedious truncation arguments in the proofs, which we do not pursue here. 

On the other hand, assumption \ref{noise-1} allows the results to be applicable for a broad class of noise models including noise models in \cref{eq:noise-models}. Note that from \ref{noise-2} \textit{\&} \citep[Proposition~2.7.1]{vershynin2018high}, we automatically also have that $\max_{i < j}\sigma_{ij}^2 \le 4\Msigma^2$. The zero-mean assumption $\E(\Eps) = \O$ can be trivially relaxed to the requirement that $H\E(\Eps)H = \O$, since the multidimensional scaling procedure operates only on the double-centered dissimilarities $\Dc = -\half HDH$. 

The lower bound in \cref{eq:noise-variance} cannot be relaxed in general. In particular, a necessary condition for \cref{eq:noise-variance} to hold is that $\#{\qty{k: \sigma_{ik}^2 > 0}} \ge p$ for every $i \in [n]$. In other words, for each $x_i$ we require at least $p$ observations in $\qty{d_{ik}: k \in [n]}$ to have non-zero variance in order to be able to construct a $p$-dimensional confidence set $C_{\alpha, i} \subset \Rp$ containing $x_i$. Moreover, since $\sum_{i}u_iu_i\tr = U\tr U = I$, a sufficient condition for \cref{eq:noise-variance} to hold is that half (or any other constant fraction $\ge p/n$) of the $\sigma_{ij}$ are bounded from below by $\msigma > 0$. For the noise models in \cref{eq:noise-models}, this is automatically satisfied for the additive noise model and for the multiplicative models in the random design setting of \cref{rem:iid}.

\section{Distributional convergence of the reconstruction error}
\label{sec:convergence}

Given the setup in \cref{sec:background} with configuration $X$ and noise $\Eps = (\eps_{ij})$, for each $i \in [n]$ let $\Sigma_i := \diag(\sigma_{i1}^2, \dots, \sigma_{in}^2)$ and $\Omega_i \in \Rpp$ be the matrix given by
\begin{align}
    \Omega_i := \frac{n}{4} \cdot ({X\tr X})^{-1} ({X\tr\Sigma_i X}) ({X\tr X})^{-1}.\label{eq:Omega-i}
\end{align}
The matrix $\Omega_i$ approximately captures the local covariance of each $\hx_i \in \Rp$ up to higher order terms. The condition in assumption \ref{noise-3} ensures that $\Om_i$ is positive definite for all $i \in [n]$. Heuristically, the noisy observations $d_{i, *} = \del_{i, *} + \epsilon_{i, *} \in \Rn$ can be viewed through the lens of linear regression; here $x_i \in \Rp$ are the ``unknown regression coefficients'' and $\hx_i$ is the estimated coefficient. In this analogy, the matrix $\Omega_i$ appearing  in \cref{eq:Omega-i} can be viewed as a rescaled analogue of White's correction for heteroscedasticity \citep{white1980heteroskedasticity}. While this analogy disregards the fact that the estimated $\h{x}_i$ are only identified up to rigid transformations, it provides some intuition for the appearance of $\Om_i$. A formal justification is provided in \cref{prop:decomposition}.

Let $G$ be a random variable following the Gumbel/Type-I extreme value distribution with c.d.f. 
\begin{align}
    \pr(G \le t) = \exp(-\exp(-t)).\label{eq:gumbel-cdf}
\end{align}
For two random variables $X$ and $Y$, with a slight abuse of notation, let $\ks(X, Y)$ denote the Kolmogorov-Smirnov metric between the distributions of $X$ and $Y$, given by
\begin{align}
    \ks(X, Y) \equiv \ks( \mathscr{L}(X), \mathscr{L}(Y) ) := \sup_{t \in \R} \abs\Big{\pr(X \le t) - \pr(Y \le t)}.\label{eq:ks-metric}
\end{align}
Our main result below establishes that, after suitable alignment and normalization, the maximum deviation of the estimated latent configuration $\hX$ from $\X$ converges to the Gumbel distribution in the $\ks$ metric.

\cbox{black!5}{%
\begin{theorem}\label{thm:main}
    Suppose $\D(X) = \Del(X) + \Eps$ satisfying \ref{assumption:compact} \& \ref{assumption:noise}, and let $\hX = \mds(D, p)$ be the output of classical multidimensional scaling. Let $\Om_i$ be given by \cref{eq:Omega-i}, and define
    \begin{align}
        T_n := \max_{i \in [n]} \sqrt{n}\norm{\Om_i^{-1/2}\qty\big(x_i -  \hg\inv(\hx_i))}.
    \end{align}
    where $\hg\inv(x) = Q\tr\hQ x$ is given in \cref{eq:hg}. Let $a_n, b_n > 0$ be two sequences given by
    \begin{align}
        b_n^2 = {2\log{n} + (p-2)\log\log{n} - 2\log{\Gamma(p/2)}} \qq{and} a_n = 1/b_n.\label{eq:anbn}
    \end{align}
    {
    Then, there exist constants $C > 0$ and $\const_1(\pars) > 0$ such that
    \begin{align}
        \ks\qty( \frac{T_n - b_n}{a_n}, G ) \lesssim \frac{\log\log{n}}{\log{n}} + \const_1(\pars) \frac{\log^3{n}}{\sqrt{n}} =: \rate_n.
        \label{eq:gumbel-convergence}
    \end{align}}
\end{theorem}
}
The result in \cref{thm:main} is non-asymptotic and applies to any $X$ satisfying \ref{assumption:compact}, and the notation $\lesssim$ in \cref{eq:gumbel-convergence} hides only absolute constants that do not depend on $n$ or the model parameters $\pars$. Note that the dominant term in the convergence rate $\rate_n$ is $O(\log\log{n}/\log{n})$, which is typical in extreme value convergence \citep[e.g.,][]{leadbetter2012extremes,hall1979rate}, and the higher-order $\log^3{n}/\sqrt{n}$ term is explicitly given because it appears again in the bootstrap results in \cref{sec:bootstrap}.

\begin{proof}[Proof Sketch]
        The core idea of the proof is to write $\sqrt{n}\Om_i^\minushalf (x_i - \hg\inv(\hx_i)) = Y_i + R_i$, where the dominant term $Y_i$ can be written as a normalized sum of independent random variables, ${Y_i = n^\minushalf\sum_{k}\eps_{ik}\theta_{ik}}$ and $R_i$ is a remainder term satisfying $\max_i\norm{R_i} = o_p(1/\log{n})$. Here, $\theta_{ij} \in \Rp$ is a deterministic vector for all $i, j \in [n]$ (see \cref{prop:decomposition}), and the contribution of $\max_i\norm{R_i}$ to the limiting distribution of $T_n$ is negligible and is handled by Slutsky's theorem. By the central limit theorem, $Y_i$ approximately follows a Gaussian distribution, and moreover, owing to \ref{noise-2}, the Cramér moderate deviation principle ensures that the tails of $\norm{Y_i}$ are captured by a $\chi^2(p)$ distribution up to vanishing relative error in the extreme value regime.
            
    To finish the proof, note that if $Y_1, \dots, Y_n$ were independent, then classical results from extreme value theory would imply convergence to the Gumbel distribution at the same rate as in \cref{eq:gumbel-convergence}. However, for each $i \neq j$, the random variables $Y_i, Y_j$ are not independent, owing to the presence of the common noise component $\eps_{ij}$. The key technical hurdle in the proof is to use the Chen-Stein Poisson approximation to show that this dependence does not affect the limiting distribution of $T_n$. The classical Poisson approximation result due to \citep{arratia1989two, arratia1990poisson} is useful when the dependency graph for the random variables is either sparse or exponentially decaying (e.g., $m-$dependent or $\psi-$mixing). On the other hand, the dependency graph for $\qty{Y_i: i \in [n]}$ here is fully connected, wherein each $Y_i$ depends on all other $Y_j$, $j \neq i$, albeit very weakly. We use the monotone coupling result of \citep{barbour1992poisson} to handle this dependence structure. The sharper rate in the second order term in \cref{eq:gumbel-convergence} is obtained by carefully analyzing the tail dependence of $Y_i$ and $Y_j$. To this end, we require a local comparison inequality for non-central Chi-squared random variables (c.f., Lemma~A.2 of \citealp{zhilova2020nonclassical}), which may be of independent interest (see \cref{lem:chisq-anticoncentration}). The proof of \cref{thm:main} is given in \cref{proof:thm:main}.
\end{proof}

{
In view of \cref{rem:iid}, if $x_1, \dots x_n$ are sampled \iid{} from a distribution $F$ on $\Rp$, then \cref{thm:main} implies the following simple corollary.

\cbox{black!5}{%
\begin{corollary}\label{cor:iid}
    Suppose $x_1, \dots, x_n \simiid{} F$ where $F$ is a distribution on $\R^p$ satisfying \cref{eq:iid-assumption}, and $D = \Del(X) + \Eps$ satisfying \ref{assumption:compact}--\ref{assumption:noise}. Let $\hX = \mds(D, p)$. Then, under the same setup as \cref{thm:main},
    \begin{align}
      \frac{T_n - b_n}{a_n} \stackrel{d}{\longrightarrow} G \qq{as} n \to \infty.\label{eq:cor-iid}
    \end{align}
\end{corollary}
}
The randomness underlying $T_n$ in \cref{eq:cor-iid} arises from both the randomness in $X$ and in $\Eps$. On the other hand, if $X_n \in \Rnp$ is a deterministic sequence of configurations satisfying \ref{assumption:compact} {(with fixed constants)} for every $n$ along the sequence $n \to \infty$, then the same result in \cref{eq:cor-iid} follows directly from \cref{thm:main}. 

We make a few remarks on the relation of \cref{thm:main,cor:iid} to existing results in literature. For a similar \iid{} setup as above, \cite{li2020central} show that for each \textit{fixed} $i \in [n]$, 
\begin{align}
    \sqrt{n}\Om_i^\minushalf(x_i - \hg\inv(\hx_i)) \stackrel{d}{\longrightarrow} N(0, I_p). 
\end{align}
The result in \cref{cor:iid} strengthens this to a uniform convergence result over all $i \in [n]$, i.e.,
\begin{align}
    \frac{\sqrt{n}\max_{i \in [n]}\norm{\Om_i^\minushalf(x_i - \hg\inv(\hx_i))} - b_n}{a_n} \stackrel{d}{\longrightarrow} G.
\end{align}
{With this, along with the fact that} $\max_i\Om_i \preccurlyeq ({\Msigma^2} {4\kappa^2}) I_p$ and $b_n = 1/a_n \asymp \sqrt{\log{n}}$, we obtain the following uniform bound on the reconstruction error:
\begin{align}
    \max_{i \in [n]}\norm{x_i - \hg\inv(\hx_i)} 
    = \Op\qty( \kappa \cdot \Msigma\sqrt{\frac{\log{n}}{n}} ),
\end{align}
which recovers the rate established in Theorem~3 of \cite{vishwanath2025minimax}.
}

We now turn our attention to constructing confidence sets for $X$. To this end, observe that the map $X \mapsto \Om_i(X)$ given in \cref{eq:Omega-i} is \textit{equivariant} under the action of $\orth{p}$, i.e., for any $O \in \orth{p}$ and the rigid transformation\footnote{Once again, we only consider the action of $\orth{p}$ since $X$ is assumed to be centered. More generally, it is easy to see that $X \mapsto \Om_i(HX)$ is \textit{invariant} to translations. Therefore the expression in \cref{eq:Omega-equivariance} holds for any rigid transformation.} $g(X) = XO\tr$, we have
\begin{align}
    \Om_i(g(X))  = O \, \Om_i(X) \, O\tr.\label{eq:Omega-equivariance}
\end{align}

Therefore, in order to characterize the local covariance information around each embedded point $\hx_i$, we need to account for the rigid transformation $\hg$ aligning $\hX$ with $X$. 

For the expression in \cref{eq:Omega-i}, the matrix $\Om_i$ captures the local covariance information in the frame of $X$. In order to construct confidence sets for each $\hx_i$, we need to transform this covariance to the frame of $\hX$ as per \cref{eq:Omega-equivariance}, i.e.,
\begin{align}
    \Om_i( \hg(X) ) = \hP\,\Om_i(X)\,\hP\tr.\label{eq:Om-hat}
\end{align}
We can then use \cref{thm:main} and invert the pivotal quantity to construct uniform confidence sets for the latent configuration $X$. Specifically, for $\alpha \in (0, 1)$ let $q_{1-\alpha} = -\log\log(1/(1-\alpha))$ be the $(1-\alpha)$-quantile of the Gumbel distribution, and let $\ellipse_{\alpha, i} \subset \Rp$ be the {ellipsoid} given by
\begin{align}
    \ellipse_{\alpha, i} := \qty{ y \in \Rp: \sqrt{n}\norm{\Om_i^\minushalf\,\hP\tr\, (y - \hx_i)} \le b_n + a_n q_{1-\alpha}}.\label{eq:confidence-set-i}
\end{align}
The following corollary shows that $\prod_{i=1}^n \ellipse_{\alpha, i}$ is a valid uniform confidence {set} for $X$.

\cbox{black!5}{%
\begin{corollary}\label{cor:confidence-pivotal}
    Consider the setup in \cref{thm:main}, and let {$\ellipse_{\alpha, i}$} be as given in \cref{eq:confidence-set-i}. Then,
    \begin{align}
        \sup_{\alpha \in (0, 1)}\abs{\pr\qty\Big(  \hg(x_i) \in \ellipse_{\alpha, i}, \; \forall i \in [n] ) - (1 - \alpha)} \lesssim \rate_n.
    \end{align}
\end{corollary}
}

In practice, the matrices $\qty{\Om_i: i \in [n]}$ are not known, and need to be estimated from the data. We may replace $\Om_i$ with any consistent estimator $\hOm_i$. {For the matrix of residuals 
\begin{align}
    E = (e_{ij}) := D - \Del(\hX)\label{eq:residuals}
\end{align}
$\h{\Sigma}_i := \diag(e_{i1}^2, \dots, e_{in}^2)$, a simple choice is the plug-in estimator:}
\begin{align}
    \hOm_i = \frac{n}{4} \cdot ({\hX\tr\hX})^{-1} ({\hX\tr\h{\Sigma}_i \hX}) ({\hX\tr\hX})^{-1}.\label{eq:hOm}
\end{align}
The resulting plug-in ellipsoids are given~by
\begin{align}
    \conf_{\alpha, i} = \qty{y \in \Rp: \sqrt{n}\norm{\hOm_i^\minushalf(y-\hx_i)} \le b_n + a_n q_{1-\alpha}}.\label{eq:confidence-set-i-hat}
\end{align}
The following result shows that $\conf_\alpha = \prod_{i=1}^n\conf_{\alpha, i}$ is also a valid uniform confidence set for $X$.

\cbox{black!5}{%
\begin{proposition}\label{prop:plug-in}
    Consider the setup in \cref{thm:main}. Let $\hOm_i$ be given by \cref{eq:hOm}, and let
    \begin{align}
        \hat T_n := \max_{i \in [n]} \sqrt{n}\norm{\hOm_i^{-1/2}(x_i - \hg\inv(\hx_i))}.\label{eq:tn-hat}
    \end{align}
    Then, for $a_n, b_n$ given in \cref{eq:anbn}, there exists $\const_2(\pars) > 0$ such that
    \begin{align}
        \ks\qty( \frac{\hat T_n - b_n}{a_n}, G ) \lesssim \rate_n + \const_2(\pars)\frac{\log^{3}{n}}{\sqrt{n}} =: \rate_n'
        .\label{eq:gumbel-convergence-hat}
    \end{align}
    where $\rate_n$ is the rate in \cref{eq:gumbel-convergence}. Moreover, for ${\conf}_{\alpha, i}$ given by \cref{eq:confidence-set-i-hat},
    \begin{align}
        \sup_{\alpha \in (0, 1)}\abs{\pr\qty\Big(  \hg(x_i) \in {\conf}_{\alpha, i}, \; \forall i \in [n] ) - (1 - \alpha)} \lesssim \rate_n'.
    \end{align}
\end{proposition}
}
Notably, since $\hOm_i$ in \cref{eq:hOm} is already capturing the covariance information in the frame of $\hX$, no additional transformations such as \cref{eq:Om-hat} are required to ensure valid coverage. We also note that other alternatives to the simple plug-in estimator above can be constructed by adapting the estimators which appear in the context of heteroscedasticity correction for regression (see, e.g., \citealp{long2000using} and the references therein).

\begin{figure}
    \centering
    \includegraphics[width=0.327\textwidth]{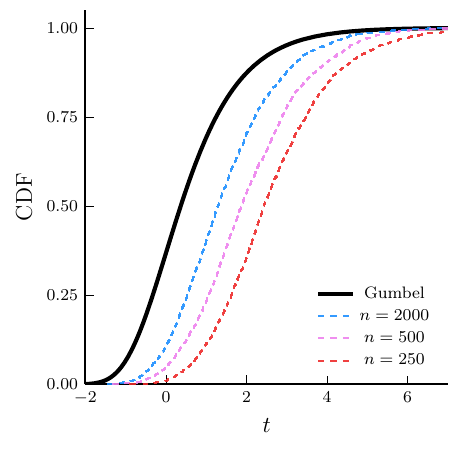}
    \includegraphics[width=0.327\textwidth]{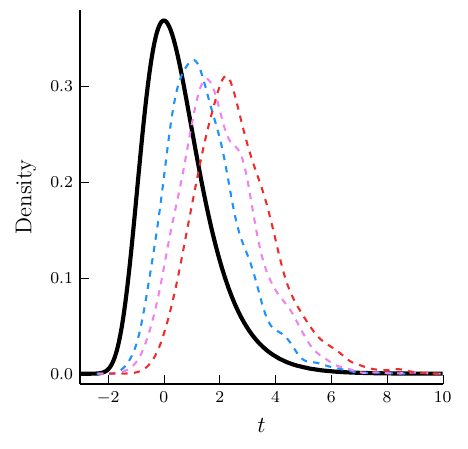}
    \includegraphics[width=0.327\textwidth]{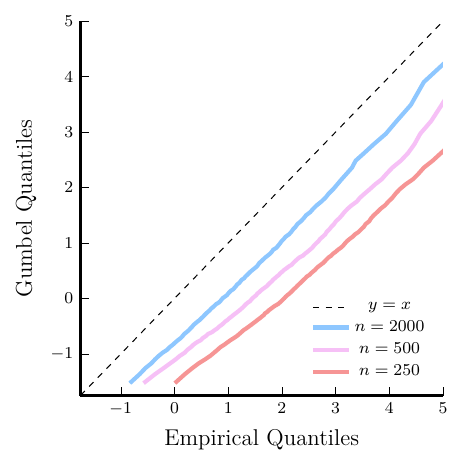}
    \caption{For $n \in \qty{250, 500, 2000}$ and $p = 5$,  latent configurations $X \in \Rnp$ from the same {distribution are generated}, and the noisy dissimilarities $D = \Del(X) + \Eps$ are generated under the additive noise model for ${\eps_{ij}\!\simiid{}\!N(0, 5)}$. (Left) The empirical c.d.f. of $(\hat{T}_n - b_n)/a_n$ is shown alongside the c.d.f. of the Gumbel distribution (Center) The kernel density estimates for the same data are compared against the p.d.f. of the Gumbel distribution. (Right) The QQ plot of the empirical quantiles vs. the Gumbel quantiles. Based on $2000$ Monte Carlo trials.}
    \label{fig:gumbel}
\end{figure}

\vspace*{-0.5em}
\section{Bootstrap Confidence Sets}
\vspace*{-0.5em}
\label{sec:bootstrap}

The main drawback in constructing confidence sets of the form \cref{eq:confidence-set-i} or \cref{eq:confidence-set-i-hat} is that the convergence to the Gumbel distribution (i.e., $\rate_n$ in \cref{eq:gumbel-convergence} and $\rate_n'$ in \cref{eq:gumbel-convergence-hat}) can be rather slow, and requires very large sample sizes in order to obtain reasonable coverage. Figure~\ref{fig:gumbel} illustrates how the empirical distribution of $(\hat{T}_n - b_n)/a_n$ compares to the Gumbel distribution for different values of $n$. In this section, we show that the bootstrap procedure can be used to construct valid confidence sets for $X$ in the noisy realizable setting. For a preview of the practical implications of the results in this section, see \cref{fig:boot-gumbel-pdf} in relation to \cref{fig:gumbel}.

\vspace*{-0.5em}
\subsection{Multiplier Bootstrap}
\vspace*{-0.5em}
\label{sec:multiplier-bootstrap}

\begin{algorithm}[t]
    \small
    \caption{\small Multiplier Bootstrap Confidence Sets for Noisy MDS}
    \label{alg:multiplier-bootstrap}
    \begin{algorithmic}[1]
        \Require Dissimilarity matrix $D \in \R^{n \times n}$, embedding dimension $p$,\\ number of bootstrap samples $B$, nominal level $\alpha \in (0, 1)$
        \State Compute $\hX \gets \mds(D, p)$ 
        \State Compute $E \gets D - \Del(\hX)$ and $\hOm_i$ using \cref{eq:hOm} for each $i \in [n]$
        \For{$b = 1$ to $B$} 
            \State Set $\Eps\b \gets R \circ E$ where $r_{ij} \simiid N(0, 1)$ for $i < j$ \Comment{Multiplier bootstrap}
            \State Generate noisy dissimilarities $D\b \gets \Del(\hX) + \Eps\b$
            \State Set $\hX\b \gets \mds(D\b, p)$ \Comment{Bootstrap embedding}
            \State Solve $\hP\b$ via orthogonal Procrustes analysis using \cref{eq:procrustes-bootstrap}
            \State Transform $\hgb\inv(\hX\b) = \hX\b\hP\b$ \Comment{Rigid transformation}
            \State $T\b_n(b) \gets \max_{i \in [n]} \sqrt{n}\norm{\hOm_i^{-1/2}(\hx_i - \hgb\inv(\hx_i\b))}$ \Comment{Bootstrap statistic}
        \EndFor
        \State Set $q_{1-\alpha}\b \gets$ the $(1-\alpha)$-quantile of $\qty{T\b_n(1), \dots, T\b_n(B)}$
        \State Compute the confidence \textit{ellipsoids} {$\conf_{\alpha, i}\b$} for each $i \in [n]$ using \cref{eq:confidence-set-i-bootstrap}.
        \State \Return Confidence sets $\conf_{\alpha}\b = \prod_{i=1}^n \conf_{\alpha, i}\b$
    \end{algorithmic}
\end{algorithm}

The multiplier bootstrap (also known as the \textit{wild} bootstrap) was originally formulated by \cite{wu1986jackknife}, and is based on the principle of externally randomizing the data to obtain a bootstrap sample. See, also, \cite{liu1988bootstrap, mammen1993bootstrap, shao2012jackknife} and the references therein for a comprehensive overview. We focus on the Gaussian multiplier bootstrap, which is arguably the most popular variant and is widely used in practice. The results below extend to other variants including \iid{} Rademacher random variables or Mammen's two-point distribution.

Let $R = (r_{ij}) \in \R^{n \times n}$ be a symmetric hollow matrix with $r_{ij} \simiid N(0, 1)$ for $i < j \in [n]$. For $\hX = \mds(D, p)$, let $\h{\Del} = \Del(\hX)$ be the pairwise Euclidean dissimilarities of $\hX$, and let {$E = (e_{ij})$} be the $n \times n$ symmetric hollow matrix of the residuals from \cref{eq:residuals}. Define $\Eps\b := R \circ E$ where $\eps\b_{ij} = r_{ij} e_{ij}$  for all $i < j$ be the externally randomized noise matrix, and let
\begin{align}
    D\b := \h{\Del} + \Eps\b \qq{and} \hX\b := \mds(D\b, p) \in \Rnp\label{eq:bootstrap-configuration}
\end{align}
denote the bootstrap dissimilarity matrix and bootstrap embedding of $D\b$, respectively. 

Conditionally on $\Eps$, $\hX=\hU\hL^{1/2}$ plays the role of the ``true'' configuration. For $\hX\b = \hU\b (\hL\b)^{1/2}$ obtained from the rank-$p$ spectral decomposition of $-\half H D\b H$, similar to \cref{eq:hg}, the optimal rigid transformation aligning $\hX$ to $\hX\b$ is simply
\begin{align}
    \hgb(x) = \hP\b x \qq{where} \hP\b = \argmin_{P \in \orth{p}}\norm{\hU\b - \hU P}^2_F,
    \label{eq:procrustes-bootstrap}
\end{align}
and $\hOm_i \in \Rpp$ plays the same role as $\Om_i$, i.e., it captures the covariance information of each $\hx_i\b$ in the frame of $\hX$ which generates the noisy dissimilarities~$\D\b$. \cref{alg:multiplier-bootstrap} summarizes the multiplier bootstrap procedure for constructing confidence sets for the latent configuration $X$.

The following result establishes the validity of the multiplier bootstrap procedure above by showing that, conditionally on $\Eps$, the distribution of $\hOm_i^{-1/2}(\hx_i - \hgb\inv(\hx\b_i))$ approximates the distribution of $\Om_i^{-1/2}(x_i - \hg\inv(\hx_i))$.

\cbox{black!5}{%
\begin{theorem}\label{thm:multiplier-bootstrap}
    Consider the setup in \cref{thm:main} with $D = \Del(\X) + \Eps$ under \ref{assumption:compact}\,\&\,\ref{assumption:noise}. Let $\hX = \mds(D, p)$ and $\hX\b = \mds(D\b, p)$ be as given in \cref{eq:bootstrap-configuration}, and define
    \begin{align}
        \h{T}_n := \max_{i \in [n]} \sqrt{n} \norm{\Om_i^{-1/2}(\hg(x_i) - \hx_i)}
        \qq{and} 
        T\b_n := \max_{i \in [n]} \sqrt{n} \norm{\hOm_i^{-1/2}(\hx_i - \hgb\inv(\hx\b_i))},\label{eq:tn-tnb}
    \end{align}
    where $\hg, \hgb$ are given in \cref{eq:hg} and \cref{eq:procrustes-bootstrap}, respectively. Then, with probability at least $1 - O(n^{-2})$ over the randomness of $\Eps$, we have
    \begin{align}
        \sup_{t \in \R}\abs{ \pr\qty\big(\h{T}_n \le t) - \pr\b\qty\big(T\b_n \le t) } \lesssim \const_1(\kappa, \Rx, \Msigma, \msigma) \frac{\log^5{n}}{\sqrt{n}} =: \rate\b_n,\label{eq:bootstrap-validity}
    \end{align}
    where $\pr\b(\cdot) = \pr(\cdot \mid \Eps)$ is the probability measure of $\Eps\b$ conditional on $\Eps$.
\end{theorem}
}

\begin{figure}[t!]
    \centering
    \includegraphics[width=0.327\textwidth]{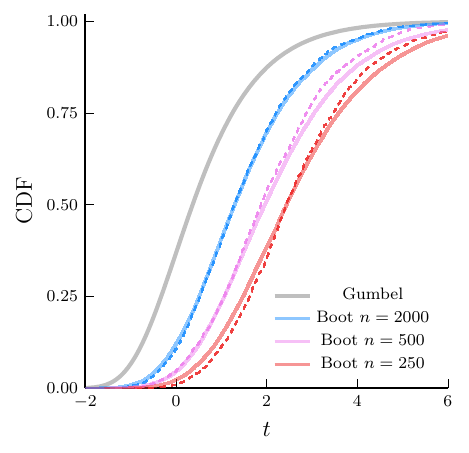}
    \includegraphics[width=0.327\textwidth]{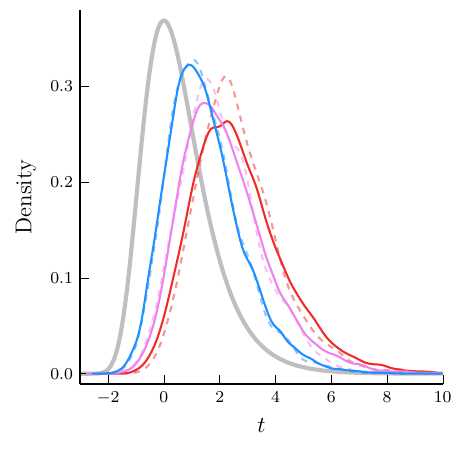}
    \includegraphics[width=0.327\textwidth]{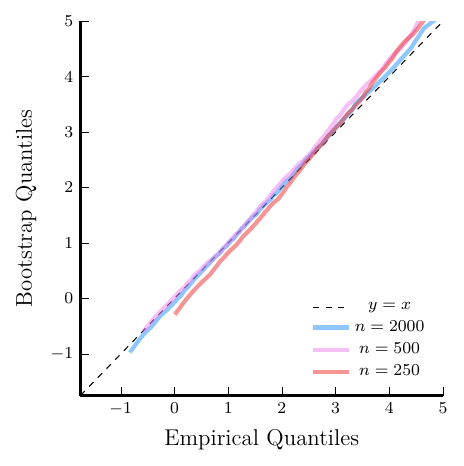}
    \caption{For the same data in \cref{fig:gumbel}, we perform  the multiplier bootstrap procedure using $B=4000$ replicates. (Left) The empirical c.d.f. of $(T_n\b - b_n)/a_n$ is compared to the c.d.f. of $(\hat{T}_n-b_n)/a_n$. (Center) The kernel density estimate based on the same bootstrap replicates is illustrated alongside the estimates from \cref{fig:gumbel}. The Gumbel c.d.f. and p.d.f. are shown in both figures for reference. (Right) The QQ plot of the empirical quantiles of $(T_n\b - b_n)/a_n$ vs. the empirical quantiles of $(\hat{T}_n - b_n)/a_n$. Based on $2000$ Monte Carlo trials.}
    \label{fig:boot-gumbel-pdf}
\end{figure}

\newpage
The proof of \cref{thm:multiplier-bootstrap} is given in \cref{proof:thm:multiplier-bootstrap}, and is based on intermediate approximations which appear in the proof of \cref{thm:main}. In contrast to the slow convergence to the Gumbel distribution in \cref{thm:main}, the bootstrap approximation in \cref{thm:multiplier-bootstrap} is substantially better, achieving nearly parametric rates up to logarithmic factors. For the same data from \cref{fig:gumbel}, the results in \cref{fig:boot-gumbel-pdf} show that the multiplier bootstrap estimation in \cref{thm:multiplier-bootstrap} is noticeably better in approximating the distribution of $\hat T_n$.

In comparison to the second term in the rate $\rate_n$ from \cref{eq:gumbel-convergence}, the rate in $\rate_n\b$ above has an extra $\log^2{n}$ factor (which arises from taking the maximum of $O(n)$ random variables with bounded $\psi_{1/2}$-Orlicz norm; see \cref{eq:bootstrap-params} and \ref{bound-6}). This may be an artifact of the proof technique, as we have not attempted to optimize the logarithmic factors in the convergence rate. 

The confidence set for $X$ can now be constructed using the bootstrap quantiles. For $\alpha \in (0, 1)$, let $q_{1-\alpha}\b$ denote the $(1-\alpha)$-quantile of the bootstrap statistic $T\b_n$, i.e.,
\begin{align}
    q_{1-\alpha}\b := \inf\qty{ t \in \R: \pr\b(T\b_n \le t) \ge 1 - \alpha }.
\end{align}
The resulting confidence set is the ellipsoid given by
\begin{align}
    \conf_{\alpha, i}\b := \qty{ y \in \Rp: \sqrt{n}\norm{\hOm_i^{-1/2}(y - \hx_i)} \le q\b_{1-\alpha}) }.\label{eq:confidence-set-i-bootstrap}
\end{align}
The coverage guarantee for $\prod_{i \in [n]}\conf_{\alpha, i}\b$ now follows from \cref{thm:multiplier-bootstrap}.

\cbox{black!5}{%
\begin{corollary}\label{cor:confidence-multiplier-bootstrap}
    Consider the setup in \cref{thm:multiplier-bootstrap}, and let $\conf_{\alpha, i}\b$ be given by \cref{eq:confidence-set-i-bootstrap}. Then, with probability greater than $1 - O(n^{-2})$ over the randomness of $\Eps$,
    \begin{align}
        \sup_{\alpha \in (0, 1)}\abs\Big{\pr\qty( \hg(x_i) \in \conf_{\alpha, i}\b, \; \forall i \in [n] ) - (1 - \alpha)} \lesssim \rate\b_n,
    \end{align}
    where $\rate\b_n$ is the convergence rate in \cref{eq:bootstrap-validity}.
\end{corollary}
}

In practice, $q_{1-\alpha}\b$ is approximated via Monte Carlo simulation, i.e., for $B$ draws of $(r_{ij})$, we can compute the bootstrap statistic $T\b_n(b)$ for $b = 1, \dots, B$ and approximate the quantile $q_{1-\alpha}\b$ as in line~10 of \cref{alg:multiplier-bootstrap}.

\subsection{Empirical Bootstrap}
\label{sec:empirical-bootstrap}

The empirical bootstrap procedure (also referred to as the \textit{nonparametric} or \textit{Efron's} bootstrap) was introduced by \cite{efron1979bootstrap}, and is arguably the most widely used bootstrap procedure in statistical estimation. While the multiplier bootstrap procedure above is valid in the heteroscedastic setting, the empirical bootstrap doesn't provide valid coverage guarantees in this setting. On the other hand, if $(\eps_{ij})$ is observed \iid{}, i.e., with equal variances, then the empirical bootstrap does provide valid coverage guarantees. 

The following result is a consequence of \cref{thm:main} and \cref{prop:plug-in}, and establishes a distributional convergence result when $(\eps_{ij})$ are \iid{}

\cbox{black!5}{%
\begin{proposition}\label{prop:iid}
    Under the conditions of \cref{thm:main}, assume that $\Eps = (\eps_{ij})$ are \iid{} with $\E(\eps_{ij})=0$, $\Var(\eps_{ij}) = \sigma^2$ and $\max_{i, j}\norm{\eps_{ij}}_{\psi_1} \le \varsigma$. For $\hX = \mds(D, p)$, define
    \begin{align}
        \tilde{T}_n := \frac{2\sqrt{n}}{\hat\sigma} \max_{i \in [n]}\Norm{\qty(\tfrac{\hX\tr\hX}{n})^{-1/2}( \hg(x_i) - \hx_i )},
    \end{align}
    where $\hat\sigma^2 = \binom{n}{2}\inv\sum_{i < j} (e_{ij} - \bar{e})^2$ is the sample variance of the residuals $(e_{ij})$ in \cref{eq:residuals}. Then, for $a_n, b_n > 0$ as given in \cref{eq:anbn},
    \begin{align}
        \ks\qty( \frac{\tilde{T}_n - b_n}{a_n}, G ) \lesssim \frac{\log\log{n}}{\log{n}} + \const_1(\kappa, \Rx, \sigma, \varsigma) \frac{\log^3{n}}{\sqrt{n}},\label{eq:gumbel-convergence-iid}
    \end{align}
    where $\const_1(\kappa, \Rx, \sigma, \varsigma)$ is the same constant as in \cref{thm:main} with $\varsigma$ and $\sigma$ in place of $\Msigma$ and $\msigma$.
\end{proposition}
}

We outline the empirical bootstrap procedure below. Let $E = D - \h{\Del}$ be the $n \times n$ matrix of residuals as in \cref{eq:residuals}. Let $\Eps\s = (\eps\s_{ij})$ be a symmetric hollow matrix where each $\eps\s_{ij}$ is an \iid{} draw from the empirical distribution of the centered residual matrix $(e_{ij} - \bar{e})$, i.e., 
\begin{align}
    \pr\s( \eps\s_{ij} = e_{kl} - \bar{e}) = \pr( \eps\s_{ij} = e_{kl} - \bar{e} \mid \Eps) = {\textstyle \binom{n}{2}}\inv \qq{for all} i < j \text{ and } k < l,\label{eq:empirical-bootstrap-measure}
\end{align}
In other words, the entries of $\Eps\s = (\eps\s_{ij})$ are obtained by sampling $\qty{e_{ij} - \bar{e}: i < j}$ with replacement. Let
\begin{align}
    D\s := \h{\Del} + \Eps\s \qq{and} \hX\s := \mds(D\s, p) \in \Rnp\label{eq:bootstrap-configuration-empirical}
\end{align}
be the bootstrap dissimilarity matrix and the bootstrap approximation of the latent configuration, respectively. Let $\hP\s \in \orth{p}$ be the Procrustes alignment given by
\begin{align}
    \hP\s = \argmin_{P \in \orth{p}}\norm{\hU\s - \hU P}^2_F\qq{and} \hgs(x) = \hP\s x.\label{eq:procrustes-bootstrap-empirical}
\end{align}
be the Frobenius-optimal rigid transformation. The resulting confidence set for $X$ is obtained similar to the multiplier bootstrap procedure in \cref{sec:multiplier-bootstrap}. 

For $\alpha \in (0, 1)$, let $
q\s_{1-\alpha} := \inf\qty\big{t \in \R: \pr\s(T\s_n \le t) \ge 1 - \alpha}
$ 
be the bootstrap quantile of $T\s_n$, and let $\conf\s_{\alpha, i}$ be the confidence ellipsoid for each $i \in [n]$ given by
\begin{align}
    \conf\s_{\alpha, i} := \qty\bigg{ y \in \R^p: \frac{\sqrt{n}}{\h{\sigma}}\Norm\Big{\qty\Big(\frac{\hX\tr\hX}{n})^{-1/2}(y - \hx_i)} \le q\s_{1-\alpha} }.\label{eq:confidence-set-i-empirical}
\end{align}
\cref{alg:empirical-bootstrap} in \cref{proof:thm:empirical-bootstrap} summarizes the empirical bootstrap procedure. 

\noindent Under the \iid{}\!\! assumption, the following result establishes the validity of the empirical bootstrap.

\cbox{black!5}{%
\begin{theorem}\label{thm:empirical-bootstrap}
    Consider the setup in \cref{prop:iid}. For $\hX = \mds(D, p)$, let $\X\s = \mds(D\s, p)$ be as given in \cref{eq:bootstrap-configuration-empirical}, and $\h{\sigma}^2 = \binom{n}{2}\inv\sum_{i < j}(e_{ij} - \bar{e})^2$. Let $\tilde{T}_n$ be as given in \cref{prop:iid}, and define
    \begin{align}
        T\s_n := \frac{2\sqrt{n}}{\hat\sigma} \max_{i \in [n]}\Norm{\qty(\tfrac{\hX\tr\hX}{n})^{-1/2}( \x_i - \hgs\inv(\hx\s_i) )},\label{eq:tn-tnstar}
    \end{align}
    where $\hg, \hgs$ are the rigid transformations given in \cref{eq:procrustes} and \cref{eq:procrustes-bootstrap-empirical}, respectively. 
    
    Then, with probability at least $1 - O(n^{-2})$ over the randomness of $\Eps$, we have
    \begin{align}
        \sup_{t \in \R}\abs{
        \pr\qty\big( {T_n} \le t ) - \pr\s\qty\big(  {T\s_n} \le t )
        } = \const_1(\kappa, \Rx, \sigma, \varsigma) \frac{\log^5{n}}{\sqrt{n}} =: \rate\s_n,\label{eq:bootstrap-validity-empirical}
    \end{align}
    where $\pr\s(\cdot) = \pr(\cdot | \Eps)$ is the empirical measure in \cref{eq:empirical-bootstrap-measure}; and, for $\conf\s_{\alpha, i}$ given in \cref{eq:confidence-set-i-empirical},
    \begin{align}
        \sup_{\alpha \in (0, 1)}\abs\Big{\pr\qty( \hg(x_i) \in \conf_{\alpha, i}\s, \; \forall i \in [n] ) - (1 - \alpha)} \lesssim \rate\s_n.\label{eq:bootstrap-validity-empirical-coverage}
    \end{align}
\end{theorem}
}

\begin{figure}[t!]
    \centering
    \subfigure{\includegraphics[width=0.83\linewidth]{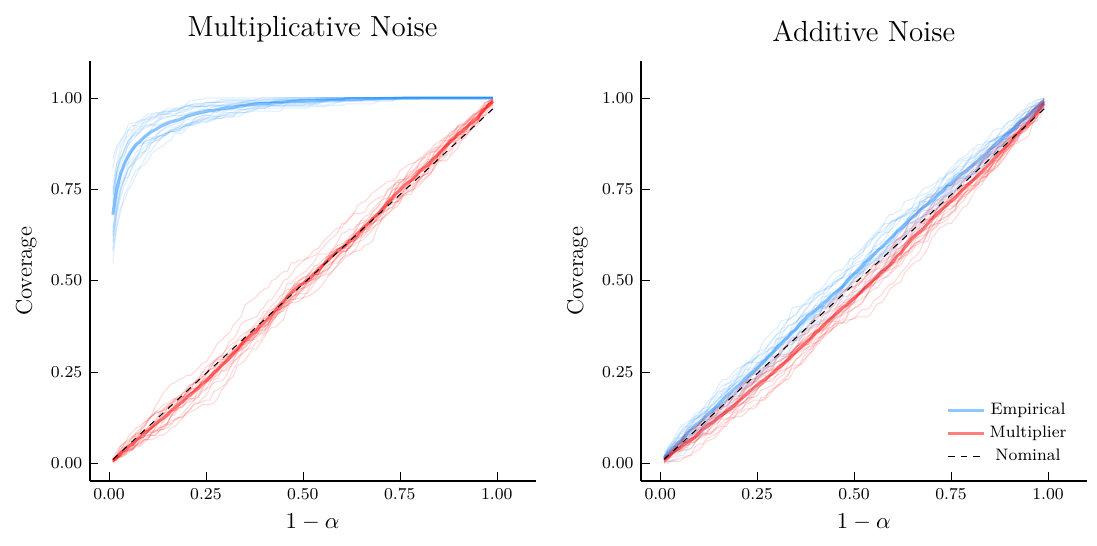}}
    \vspace*{-1em}
    \caption{\emph{Coverage probabilities for the multiplier bootstrap and the empirical bootstrap for different noise models.} For $N = 20$ different configurations, $X \in \R^{n \times 2}$, noisy dissimilarities $D$ are obtained using (left) multiplicative noise and (right) additive noise. Bootstrap confidence sets are computed using both the multiplier bootstrap and the empirical bootstrap procedures for a range of nominal levels $\alpha \in (0, 1)$, and the coverage probabilities are computed across $500$ Monte Carlo runs. Each of the $N$ thin lines correspond to the coverage probabilities obtained for a particular fixed configuration $X$, and the thick lines correspond to the average coverage across all configurations.
    }\label{fig:boot}
\end{figure}

\section{Numerical Experiments}
\label{sec:experiments}

We present some numerical experiments to illustrate the theoretical results in \cref{sec:convergence,sec:bootstrap}.

\textbf{Experiment 1.} (\textit{Multiplier vs. Empirical Bootstrap})\quad In the first experiment, for $n=500$ and $p=2$, we consider $N=20$ different configurations $X \in \Rnp$ which are all sampled  uniformly from an elliptical shape with eccentricity $2$. We consider two different noise models: (i) multiplicative noise where $\eps_{ij} \sim N(0, \sigma^2\del_{ij}^2)$ and (b) additive noise where $\eps_{ij} \simiid{} N(0, \sigma^2)$ for $i \neq j$.  In both cases, we fix $\sigma = 1.0$ and compute confidence sets using $B=500$ bootstrap replications using the multiplier bootstrap (\cref{alg:multiplier-bootstrap}) and the empirical bootstrap (\cref{alg:empirical-bootstrap}) procedures. \cref{fig:boot} plots the coverage probabilities for a range of $\alpha$ values computed across $500$ Monte Carlo trials. 

In the additive \iid{} noise setting, both the multiplier bootstrap and the empirical bootstrap yield valid coverage guarantees for the latent configuration $X$ as corroborated by \cref{thm:multiplier-bootstrap} and \cref{thm:empirical-bootstrap}. On the other hand, in the multiplicative noise setting, the empirical bootstrap doesn't provide valid confidence sets---the confidence sets are too conservative, leading to over-coverage. The multiplier bootstrap procedure, however, still provides valid coverage.

\begin{table}[t!]
\footnotesize
\centering
\caption{Coverage probabilities for different methods of constructing confidence sets for the setup in Experiment~2.}
\label{tab:coverage}
\resizebox{\linewidth}{!}{%
\begin{tabular}{ll*{9}{c}}
\toprule
\multirow{2}{*}{\textbf{Noise Model}} & \multirow{2}{*}{\textbf{Method}} & \multicolumn{9}{c}{\textbf{Nominal level }$(1-\alpha)$} \\
\cmidrule(lr){3-11}
& & $0.999$ & $0.99$ & $0.975$ & $0.95$ & $0.925$ & $0.9$ & $0.85$ & $0.8$ & $0.75$ \\
\midrule
\multirow{5}{*}{Additive}
    & Gaussian       & $0.997$ & $0.988$ & $0.967$ & $0.943$ & $0.915$ & $\mathbf{0.890}$ & $\mathbf{0.838}$ & $\mathbf{0.788}$ & $\mathbf{0.747}$ \\
    & Rademacher     & $0.997$ & $0.984$ & $0.968$ & $0.935$ & $0.906$ & $0.870$ & $0.813$ & $0.770$ & $0.712$ \\
    & Uniform        & $0.996$ & $0.986$ & $0.967$ & $0.943$ & {$\mathbf{0.916}$} & $0.883$ & $0.818$ & $0.779$ & $0.724$ \\
    & Empirical  & {$\mathbf{0.998}$} & {$\mathbf{0.989}$} & {$\mathbf{0.975}$} & {$\mathbf{0.955}$} & $0.940$ & $0.924$ & $0.867$ & $0.823$ & $0.795$ \\
    & Gumbel         & $0.985$ & $0.946$ & $0.911$ & $0.861$ & $0.821$ & $0.790$ & $0.730$ & $0.658$ & $0.603$ \\
    \midrule
    \multirow{5}{*}{Log-normal}
    & Gaussian       & {$\mathbf{0.998}$} & $0.992$ & $0.983$ & $0.954$ & $0.932$ & {$\mathbf{0.902}$} & {$\mathbf{0.851}$} & {$\mathbf{0.789}$} & {$\mathbf{0.707}$} \\
    & Rademacher     & $0.997$ & $0.988$ & {$\mathbf{0.980}$} & {$\mathbf{0.951}$} & {$\mathbf{0.924}$} & $0.895$ & $0.825$ & $0.738$ & $0.628$ \\
    & Uniform        & $0.997$ & {$\mathbf{0.991}$} & $0.983$ & $0.953$ & $0.930$ & $0.897$ & $0.837$ & $0.750$ & $0.659$ \\
    & Empirical  & {${1.000}$} & $1.000$ & $1.000$ & $1.000$ & $0.997$ & $0.994$ & $0.994$ & $0.953$ & $0.983$ \\
    & Gumbel         & $0.939$ & $0.847$ & $0.772$ & $0.690$ & $0.630$ & $0.580$ & $0.508$ & $0.423$ & $0.356$ \\
\bottomrule
\end{tabular}
}
\end{table}

\textbf{Experiment 2.}~(\textit{Comparison of multipliers in different noise settings})\quad As noted in \cref{sec:multiplier-bootstrap}, the multiplier bootstrap procedure is valid for a wide class of multipliers $(r_{ij})$ beyond the Gaussian multipliers considered in \cref{alg:multiplier-bootstrap}. The only requirement we have in our proofs is that $\E(r_{ij}) = 0$, $\Var(r_{ij})=1$ and $\max_{i, j}\norm{r_{ij}}_{\psi_2} < \infty$. 

In this experiment, we consider: (i) Gaussian multipliers, $r_{ij} \sim N(0, 1)$, (ii) Rademacher multipliers, $r_{ij} \sim \text{Ber}(\qty{+1, -1}; 1/2)$, and (iii) Uniform multipliers, $r_{ij} \sim \text{Unif}([- \sqrt{3}, \sqrt{3}])$. We also benchmark the performance of the multiplier bootstrap procedures against (iv)~the empirical bootstrap, and (v)~the extreme value approximation in \cref{prop:plug-in}. We take $X \in \R^{n \times 2}$ to be the locations (latitude/longitude) of $n=350$ largest cities in the U.S., and generate noisy dissimilarities using the additive and log-normal noise models described in~\cref{eq:noise-models}.

\cref{tab:coverage} reports the coverage probabilities for a range of nominal levels $(1-\alpha)$ computed across $1000$ Monte Carlo trials. For each bootstrap method, we generate $B=1000$ bootstrap replicates to compute the bootstrap quantiles. All the multiplier bootstrap methods yield valid coverage guarantees across both noise models. We also find that using the empirical bootstrap procedure, when valid (i.e., when the noise is additive), yields marginally better coverage in the extreme tails compared to the multiplier bootstrap procedures.

\begin{figure}[h!]
    \includegraphics[width=0.328\linewidth]{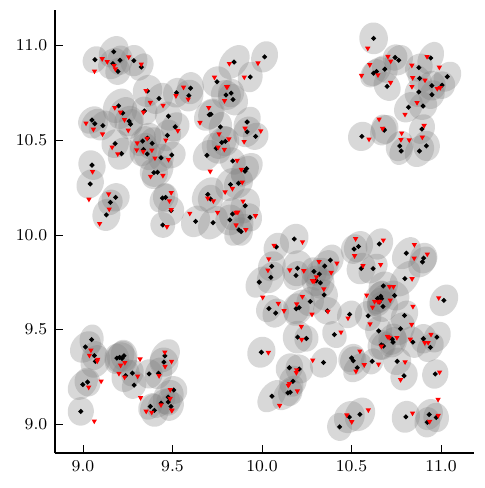}
    \includegraphics[width=0.328\linewidth]{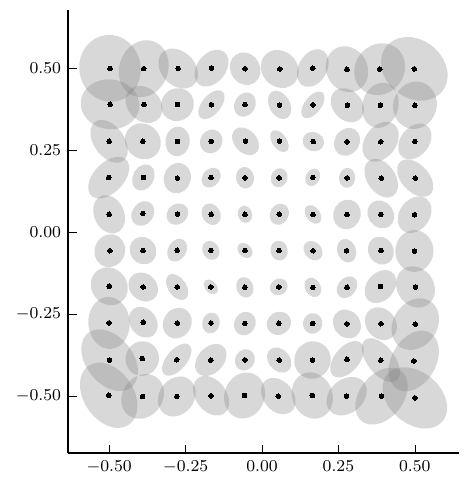}
    \includegraphics[width=0.328\linewidth]{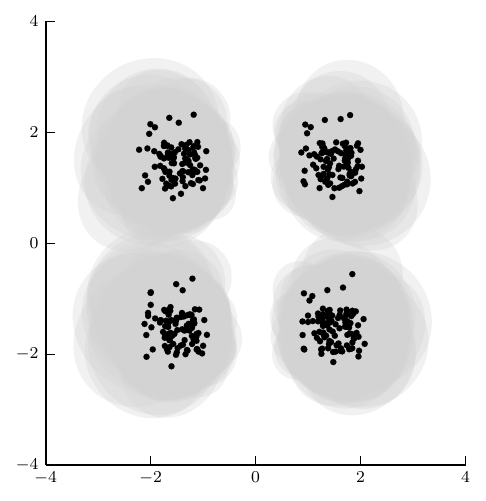}
    \subfigure[Additive vs. multiplicative noise\label{subfig:toy-a}]{\includegraphics[width=0.328\linewidth]{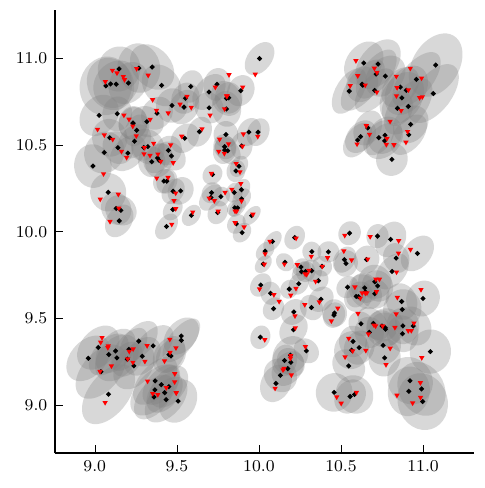}}\,
    \subfigure[Radial sum vs. difference noise\label{subfig:toy-b}]{\includegraphics[width=0.328\linewidth]{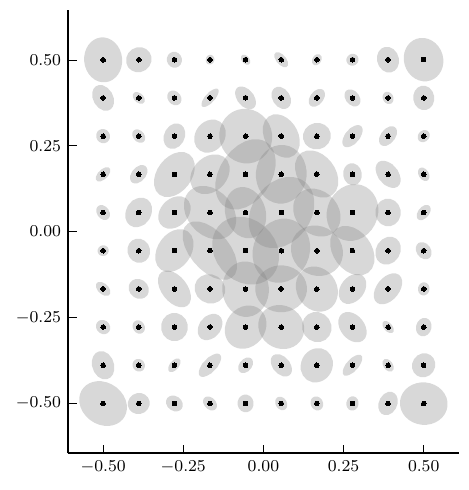}}
    \subfigure[Horizontal vs. vertical noise\label{subfig:toy-c}]{\includegraphics[width=0.328\linewidth]{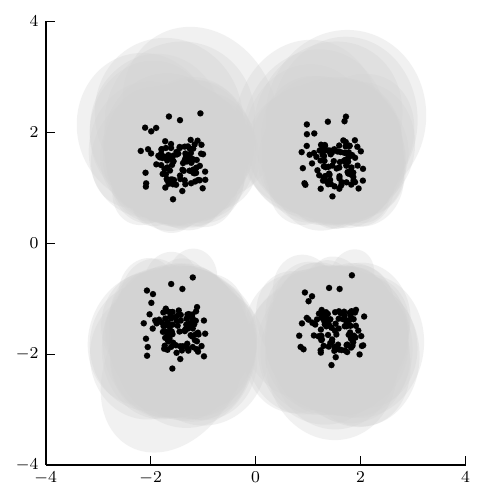}}
    \caption{
    \emph{Adaptivity of the multiplier bootstrap confidence sets to heteroscedasticity.} Points are sampled from a configuration $X \in \Rnp$ and various noise models, the embedding $\hX \in \Rnp$ is obtained via the classical MDS algorithm (black $\bullet$) and the confidence sets (grey ellipsoids) are computed using \cref{alg:multiplier-bootstrap}. (a) The noise is additive for in the figure on top and multiplicative in the figure below. The latent configuration is shown in red (\textcolor{red}{$\blacktriangledown$}). (b) For each pair of points, the noise variance depends on: (top) the sum of each point's squared norm and (bottom) the absolute difference of each point's squared norm. (c) The noise variance depends on: (top) the vertical pairwise distances and (bottom) the horizontal pairwise distances.
    }\label{fig:toy}
\end{figure}

\textsc{Experiment 3.}~(\textit{Adaptivity of the multiplier bootstrap to heteroscedasticity.})\quad In the final set of experiments, we demonstrate how the multiplier bootstrap confidence sets adapt to varying noise structures. In \cref{subfig:toy-a}, we fix a configuration of $n=150$ points sampled from a mixture of $4$ different squares in $\R^2$ (shown in red). We generate  noisy dissimilarities $D = \Del + \Eps$ under two noise models: (top) additive noise, $\eps_{ij} \sim N(0, \sigma^2)$ and (bottom) multiplicative noise, $\eps_{ij} \sim N(0, \sigma^2\del_{ij}^2)$ with $\sigma=0.4$. The confidence sets, computed using \cref{alg:multiplier-bootstrap}, are shown in grey. We note that both the shape and size of the resulting confidence sets adapt to the underlying noise. In particular, under additive noise, the sets around each $\hx_i$ are approximately spherical with similar radii, as expected from \cref{eq:confidence-set-i-empirical}. Under multiplicative noise, the sets become ellipsoidal and vary in size depending on the local noise level: points near the centroid have smaller variance (and thus smaller confidence sets), while points farther away from the centroid (e.g., those in the first and third quadrants) have larger variances, and therefore, larger confidence sets.

To further examine the adaptivity to heteroscedasticity, in \cref{subfig:toy-b} we consider a configuration of $n=100$ points uniformly placed on a square grid in $\R^2$. We generate noisy dissimilarities as follows: (top) $\eps_{ij} \sim N\qty\big(0, \norm{x_i}^2 + \norm{x_j}^2)$ and (bottom) $\eps_{ij} \sim N\qty\big(0, \abs\big{\norm{x_i}^2 - \norm{x_j}^2})$. The resulting confidence sets (shown in gray) capture the varying noise structure. In the first case, points farther from the origin have larger variances and thus larger confidence sets. In the second case, points with equal radii have similar noise variances. Moreover, the outermost ring---which contains the most points with identical radii---all have zero variance for their respective pairwise entries in the noisy dissimilarity matrix. On the other hand, points farther from this ring, i.e., points closer to the center or points at the corners of the grid, have the highest noise variances, as reflected in their confidence sets.

Beyond overall coverage, properly accounting for the noise in the dissimilarities can affect the inference from the embeddings. To illustrate this, in \cref{subfig:toy-c} we consider two noisy dissimilarity matrices with the same latent configuration but differing in their noise structures:
(top) variance depends only on vertical pairwise distances, and (bottom) variance depends only on horizontal pairwise distances. While the resulting embeddings in are visually hard to distinguish, their 95\% confidence sets show different patterns: there greater evidence for horizontal separation between the clusters on top, whereas, in the bottom figure, there is more evidence of vertical separation. Both are consistent with their respective noise structures.


\section{Discussion}
\label{sec:discussion}

\vspace*{-1em}

Our work places classical multidimensional scaling within a formal statistical framework. The distributional convergence results in \cref{sec:convergence} establishes the basis for constructing uniform confidence sets for the latent configuration, up to rigid transformations. The bootstrap procedures in \cref{sec:bootstrap} provide practical and efficient algorithms for constructing these confidence sets.

While our focus has been on constructing confidence sets, extending this framework to other inferential tasks may be of interest to practitioners, e.g., goodness-of-fit tests for the latent configuration, or two-sample tests for comparing the configurations underlying different dissimilarity matrices. Additionally, the theoretical guarantees obtained here apply when the noise is sufficiently regular, i.e., in the absence of (possibly adversarial) outliers or missing dissimilarities. Another practically relevant direction would be to develop an inferential framework for other MDS methods that are better able to handle (severe) outliers and/or missingness.

Our analysis considers the noisy realizable setting where the observed dissimilarities take the form: $d_{ij} = \norm{x_i - x_j}^2 + \eps_{ij}$, i.e., noise is added to the squared Euclidean distances between latent points lying in some low-dimensional subspace ($x_1, \dots, x_n \in \R^p$ for fixed $p < n$). An alternative and complementary framework considers the setting: $d_{ij} = \norm{y_i - y_j}^2$ where $y_i = Rx_i \in \R^m$ and $R \in \R^{m \times p}$ for $p \ll m$ is a random matrix which embeds the low-dimensional latent points into a higher dimensional space (see, e.g., \citealp{peterfreund2021multidimensional} and \citealp{little2023analysis}). Analyzing the statistical behavior of CMDS in this setting, particularly in high-dimensional regimes, is an interesting and open problem.

As noted in Section \ref{sec:introduction}, CMDS forms the basis for several embedding methods such as landmark MDS \citep{de2004sparse}, Isomap \citep{tenenbaum2000global}, and maximum variance unfolding \citep{weinberger2006introduction}, {and is often used in patch-based algorithms \citep[e.g.,][]{shang2004localization}}. Establishing similar results for these related methods is a promising direction for future work.


\clearpage
\section{Proofs}
\label{sec:proofs}

This section contains the proofs of the main results. In the interest of clarity, and to avoid notational clutter, throughout the proofs we will write $\const_{\square}$, $C_{\square}$, $c_{\square}$, $C'_{\square}$, etc., to denote constants $\const_{\square}(\pars)$, $C_{\square}(\pars)$, $c_{\square}(\pars)$, $C'_{\square}(\pars)$, etc. which depend only on the parameters $\pars$. Almost always $C, c > 0$ \textit{without any sub/super-scripts} are used to denote absolute constants. {Throughout the proofs, the notation $O(\dots)$ only suppresses constants possibly depending on $p$.}

We first present the following few lemmas which are used in the proofs. The first lemma is a well-known quantitative version of Slutsky's theorem and comes in handy for establishing the distributional convergence in the presence of relatively small remainder terms.

\cbox{black!5}
{%
\begin{lemma}\label{lem:slutsky}
    Let $S_n, T_n$ be {sequences} of random variables and $T$ a random variable such that
    \begin{align}
        \ks\qty( S_n, T ) = O(s_n) \qq{and} \pr\qty\Big( \abs{T_n - S_n} \ge u_n ) = O(r_n)
    \end{align}
    for some $C > 0$ and non-negative sequences {$u_n$} and $r_n, s_n = o(1)$. Then,
    \begin{align}
        \ks\qty(T_n, T) = O\qty\big( r_n + s_n + \omega_T(u_n) ),
    \end{align}
    where 
    $
    \omega_T(\eta) := \sup\qty\big{ \pr(t < T \le t+\eta): t \in \R}
    $
    is the modulus of continuity {of the c.d.f of $T$}. Moreover, if $T$ admits a p.d.f. uniformly bounded by $M > 0$, then $\omega_T(\epsilon) \le M\epsilon$.
\end{lemma}
}
Since the statement in this form was not available in standard references, the proof is provided in \cref{proof:lem:slutsky} for completeness. The next lemma characterizes the normalizing sequences $a_n, b_n$ in \cref{thm:main}, and the proof is deferred to \cref{proof:lem:chisq-tail}.

\cbox{black!5}
{%
\begin{lemma}\label{lem:chisq-tail}
    Let $Z \sim N(0, I_p)$ and $u_n(t) := a_nt + b_n$ for $a_n, b_n$ given in \cref{eq:anbn}. Then, for all $t \in \R$,
    \begin{align}
        \pr\qty\Big( \norm{Z} > u_n(t) ) = \frac{1}{n}e^{-t - t^2/2b_n^2} \qty( 1 + O\qty(\tfrac{\abs{t} + \log\log{n}}{\log{n}}) ).
    \end{align}
\end{lemma}
}

\subsection{Proof of \cref{thm:main}}
\label{proof:thm:main}

As noted in the proof sketch, we begin by writing $\Om_i^{-1/2}(x_i - \hg\inv(\hx_i))$ as normalized sum of independent random vectors plus a remainder term as follows.

\cbox{black!5}
{%
\begin{proposition}\label{prop:decomposition}
    Consider the setup in \cref{thm:main} where $X = U\L^{1/2} Q$, $\hX = \mds(D, p)$, $\hg$ is the rigid transformation given in \cref{eq:hg}, and $\Om_i$ is as given in \cref{eq:Omega-i}. Then, for each $i \in [n]$,
    \begin{align}
        \sqrt{n} \cdot \Om_i^{-1/2}(x_i - \hg\inv(\hx_i)) = Y_i + R_i,\label{eq:mds-decomposition}
    \end{align}
    where,
    \begin{align}
        Y_i := \frac{1}{\sqrt{n}} \sum_{k \in [n]} \eps_{ik} \theta_{ik} \qq{for} \theta_{ik} := \half\Om_i^\minushalf \qty(\tfrac{X\tr X}{n})\inv x_k,\label{eq:yi}
    \end{align}
    with
    \begin{align}
        ({\msigma^2}/{4\kappa^2}) I_p \preccurlyeq \Om_i \preccurlyeq ({\Msigma^2} {4\kappa^2}) I_p \qq{and}
        \max_{i, k}\norm{\theta_{ik}}  \le \frac{\kappa^3\Rx}{2\msigma} =: \const_0(\kappa,\Rx,\msigma),\label{eq:theta-norm}
    \end{align}
    and $R_{i}$ is a remainder term such that with probability greater than $1 - O(n^{-2})$,
    \begin{align}
        \max_i \norm{R_{i}} \lesssim {C}'_1(\pars)\sqrt{\frac{\log{n}}{n}}.\label{eq:remainder-term}
    \end{align}
\end{proposition}
}

From \cref{prop:decomposition}, it is clear that the $Y_i$ variables in \cref{eq:mds-decomposition} contribute to the dominant terms in $T_n$. To this end, let $\const_1 := (\const_0\Msigma)^2$ for $\const_0$ in \cref{prop:decomposition}, and let $M_n$ be defined as
\begin{align}
    M_n := \max_{i \in [n]} \norm{Y_i}.\label{eq:mn}
\end{align}
Using Slutsky's theorem in \cref{lem:slutsky}, we can restrict our attention to $M_n$ alone. Specifically, for $b_n \sim \sqrt{2\log{n}}$ and $a_n = 1/b_n$, from \cref{prop:decomposition} we have that with probability greater than $1 - O(n^{-2})$,
\begin{align}
    \abs{ \qty(\frac{T_n - b_n}{a_n}) - \qty(\frac{M_n - b_n}{a_n}) } \le \frac{1}{a_n} \max_{i \in [n]}\norm{R_i} \lesssim C'_1(\pars) \cdot \frac{\log{n}}{\sqrt{n}}.\label{eq:remainder-term-mn}
\end{align}
Also, {for the Gumbel distribution the p.d.f. satisfies} $f_G(t) \le e^{-1}$, and, therefore, for any $\epsilon > 0$,
\begin{align}
    \omega_G(\epsilon) := \sup\qty\big{ \pr(t \le G \le t+h): t \in \R, h \le \epsilon } \le \epsilon e^{-1} \le \epsilon.\label{eq:omega-gumbel}
\end{align}
If we can show that:
\begin{align}
    \sup_{t \in \R}\abs{\pr\qty( \frac{M_n - b_n}{a_n} \le t ) - \pr\qty( G \le t )} \lesssim \rate_n,\label{eq:claim-mn}
\end{align}
then the conclusion in \cref{eq:gumbel-convergence} follows from \cref{lem:slutsky} by combining \cref{eq:claim-mn} with the bound in \cref{eq:omega-gumbel} and the tail bound in \cref{eq:remainder-term-mn} and by noting that $C'_1\log{n}/\sqrt{n} = o(\rate_n)$. Therefore, the remainder of the proof is devoted to establishing the claim in \cref{eq:claim-mn}.

For $t \in \R$ and $a_n, b_n$ given in \cref{eq:anbn}, define $u_n(t) :=  a_n t + b_n$, and let
\begin{align}
    \lambda_n(t) := \sum_{i \in [n]} \pr\qty\Big( \norm{Y_i} > u_n(t) ).\label{eq:lambda-n}
\end{align}
Using the triangle inequality and by noting that $\pr(G \le t) = e^{-e^{-t}}$, we have
\begin{align}
    \sup_{t \in \R}\abs{\pr\qty( \frac{M_n - b_n}{a_n} \le t ) - \pr\qty( G \le t )} &\le \underbrace{\sup_{t \in \R}\abs{\pr\qty( {M_n} \le u_n(t) ) - e^{-\lambda_n(t)}}}_{=:\Circled{1}} + \underbrace{\sup_{t \in \R}\abs{e^{-\lambda_n(t)} - e^{-e^{-t}}}}_{=:\Circled{2}}.\label{eq:triangle-inequality}
\end{align}
The claim in \cref{eq:claim-mn} follows by establishing that
\begin{align}
    \Circled{1} \lesssim \const_1\frac{\log^3{n}}{\sqrt{n}} \qq{and} \Circled{2} \lesssim \frac{\log\log{n}}{\log{n}}.\label{eq:claim-1-2}
\end{align}
We prove these two bounds for \cref{eq:claim-1-2} in \cref{proof:claim-circled-1,proof:claim-circled-2}, respectively. To this end, the following lemma characterizes the tail behavior of the term $Y_i$ defined in \cref{eq:yi}.

\cbox{black!5}
{%
\begin{lemma}\label{lem:tail-bound}
    Let $Y_i$ be as given in \cref{eq:yi}, $u_n(t) := b_n + ta_n$ for $a_n, b_n$ given in \cref{eq:anbn}. Let $Z \sim N(0, I_p)$. There exists $\tau = \tau(\pars) > 0$ such that for all $i \in [n]$,
    \begin{align}
        \pr\qty\big(\norm{Y_i} > u_n(t) ) = 
        \begin{dcases}
            1 & \qq{if } t \in (-\infty, -b_n^2]\\
            \pr\qty\big(\norm{Z} > u_n(t)) \qty( 1 + O\qty(\tfrac{\const_1 \log^{3/2}{n}}{\sqrt{n}}))  & \qq{if } t \in (-b_n^2, \tau \log{n}]\\
            O(1/n^2) & \qq{if } t \in (\tau\log{n}, \infty).
        \end{dcases}
    \end{align}
\end{lemma}
}

{The $O(\dots)$ terms above do not depend on $t$. In particular, the $O(1/n^2)$ holds uniformly for all $t > \tau \log{n}$.} The proof of \cref{lem:tail-bound} is deferred to \cref{proof:lem:tail-bound}. Throughout, we also use the fact that for all $t \in (-b_n^2, \tau\log{n})$, by combining \cref{lem:tail-bound,lem:chisq-tail} we have $\lambda_n(t) \sim n\pr(\norm{Z} > u_n(t))$, or, equivalently,
\begin{align}\label{eq:lambda-n-k2}
    \lambda_n(t) = e^{-t - t^2/2b_n^2}\qty\big(1 + \zeta_n(t))
    \qq{where} 
    \abs{\zeta_n(t)} \lesssim \tfrac{\abs{t} + \log\log{n}}{\log{n}} + \tfrac{\const_1\log^{3/2}{n}}{\sqrt{n}}.
\end{align}

\subsubsection{Bound for $\Circled{1}$}
\label{proof:claim-circled-1}

Let $\tau > 0$ be as given in \cref{lem:tail-bound}, $\beta = 1/({2 + \sqrt{2}})$, and consider the following four intervals:
\begin{align}
    J_1 := \big(-\infty, -b_n^2\big]\qc{} J_2 := \big(-b_n^2, -\beta b_n^2\big]\qc{} J_3 := \big(-\beta b_n^2, \tau\log{n}\big]\qc{} J_4 := \big(\tau\log{n}, \infty\big).\label{eq:intervals}
\end{align}

\subsubsection*{\thesubsubsection~(i). $t \in J_1$.}

Since $u_n(t) = a_n t + b_n = t/b_n + b_n \le 0$ for $t \in J_1$, 
\begin{align}
    \pr(M_n \le u_n(t)) = 0\qq{and} \lambda_n(t) = \sum_{i \in [n]} \pr\qty( \norm{Y_i} > u_n(t) ) = n,\label{eq:lambda-n-j1}
\end{align}
we have
\begin{align}
    \sup_{t \in J_1}\abs\big{ \pr(M_n \le u_n(t)) - e^{-\lambda_n(t)} } = \sup_{t \in J_1}e^{-\lambda_n(t)} &\le e^{-n}.\label{eq:bound-circled-1-j1}
\end{align}

\subsubsection*{\thesubsubsection~(ii). $t \in J_4$.}

In terms of the upper-tail probability, we have
\begin{align}
    \sup_{t \in J_4}\abs\big{ \pr(M_n \le u_n(t)) - e^{-\lambda_n(t)} } 
    &= \sup_{t \in J_4}\abs\big{ \pr(M_n > u_n(t)) - (1-e^{-\lambda_n(t)}) }\label{eq:bound-circled-1-j4-1}
\end{align}
Note that when $t > \tau \log{n}$, \cref{lem:tail-bound} gives
\begin{align}
    \lambda_n(t) = \sum_{i\in [n]} \pr\qty(\norm{Y_i} > u_n(t)) = O(1/n),\label{eq:lambda-n-j4}
\end{align}
and, using a union bound,
$
\pr( M_n > u_n(t) ) \le n \cdot \max_i\pr( \norm{Y_i} > u_n(t) ) = O(1/n).
$ 
Using the triangle inequality in \cref{eq:bound-circled-1-j4-1} and the fact that $1 - e^{-z} \le z$ for $z \ge 0$, we get
\begin{align}
    \sup_{t \in J_4}\abs\big{ \pr(M_n \le u_n(t)) - e^{-\lambda_n(t)} } 
    \le \sup_{t \in J_4}\qty{ n\max_{i} \pr(\norm{Y_1} > u_n(t)) + \lambda_n(t) } =  O(1/n).\label{eq:bound-circled-1-j4}
\end{align}

\subsubsection*{\thesubsubsection~(iii). $t \in J_2, J_3$.}

We define some additional quantities. Let
\begin{align}
    B_i(t) := \mathbb{1}\qty\Big( \norm{Y_i} > u_n(t) ),\quad \pi_i(t) := \pr\qty\Big( \norm{Y_i} > u_n(t) ),\qq{and} W(t) := \sum_{i \in [n]} B_i(t).\label{eq:B-pi-W}
\end{align}
Note that $B_i(t) \sim \text{Ber}(\pi_i(t))$ for $i \in [n]$ and $\E(W(t)) = \lambda_n(t)$. We also need the following bound for $\Cov(B_i(t), B_j(t))$, which is the main technical hurdle in this proof.

\cbox{black!5}
{%
\begin{lemma}\label{lem:covariance-bound}
    {For any $t \in J_2 \cup J_3$}, let $B_i(t) := \mathbb{1}\qty{ \norm{Y_i} > u_n(t) }$ for $Y_i$ given in \cref{eq:yi}. Then, for all $i \neq j$,
    \begin{align}
        \abs{\Cov(B_i(t), B_j(t))} \lesssim \const_1 \cdot \frac{\log^3{n}}{n} \pr\qty\Big( \norm{Z} > u_n(t) )^2 + O(n^{-4}).\label{eq:covariance-bound}
    \end{align}
\end{lemma}
}
We again note that the $O(n^{-4})$ term above does not depend on $t$. The proof of \cref{lem:covariance-bound} is in \cref{proof:lem:covariance-bound}, and is based on a local comparison inequality for non-central Chi-squared random variables, which may be of independent interest.

For $t \in J_2 \cup J_3$, from \cref{eq:lambda-n-j1,eq:bound-circled-1-j4-1} note that $\lambda_n(t)$ decreases from $\lambda_n(t) = n$ when $t = -b_n^2$ to $\lambda_n(t) = O(1/n)$ when $t = \tau\log{n}$. At $t_n := -\beta b_n^2$, we have
\begin{align}
    -t_n - \frac{t_n^2}{2b_n^2} = \frac{b_n^2}{(2+\sqrt{2})} - \frac{b_n^2}{2(2+\sqrt{2})^2} = \frac{3 + 2\sqrt{2}}{4(3 + 2\sqrt{2})}b_n^2 = \frac{b_n^2}{4} \sim \half\log{n},
\end{align}
and, therefore, from \cref{lem:chisq-tail},
\begin{align}
    \lambda_n(t_n) = e^{-t_n - t_n^2/2b_n^2} (1 + \zeta_n(t_n)) \asymp \sqrt{n}.
    \label{eq:lambda-n-tn}
\end{align}
This implies that, $\lambda_n(t) \gtrsim \sqrt{n}$ {uniformly over $J_2$} and $\lambda_n(t) \lesssim \sqrt{n}$ {uniformly over $J_3$}.  We use two different results to bound $\Circled{1}$ based on the value of $\lambda_n(t)$.

\subsubsection*{\thesubsubsection~(iii--a). $t \in J_2$.}

Since $\qty\big{ M_n \le u_n(t) } = \qty\big{ W(t) = 0 } \subseteq \qty{\abs{W(t) - \E W(t)} \ge \E W(t)}$,
\begin{align}\label{eq:circled-1-W-bound}
    \abs\big{ \pr\qty({M_n} \le u_n(t)) - e^{-\lambda_n(t)} }
    \le \pr(W(t) = 0) + e^{-\lambda_n(t)}
    &\le \pr\qty\Big(\abs{W - \lambda_n(t)} \ge \lambda_n(t)) + e^{-\lambda_n(t)}\\
    &\le \frac{\Var(W(t))}{\lambda_n(t)^2} + e^{-\lambda_n(t)},\label{eq:bound-circled-1-j2-1}
\end{align}
where last line follows from an application of Chebyshev's inequality. Now, using \cref{lem:covariance-bound} and by noting that $\lambda_n(t) \sim n \pr(\norm{Z} > u_n(t))$ from \cref{lem:tail-bound}, we have
\begin{align}
    \Var\qty\Big(\sum_{i \in [n]} B_i(t))
    &\le \sum_{i \in [n]}\pi_i(t) + n^2 \cdot \max_{i, j}\abs{\Cov(B_i(t), B_j(t))} \\
    &\lesssim \lambda_n(t) + n^2 \cdot \qty(\const_1 \log^3{n} \cdot \frac{\pr(\norm{Z} > u_n(t))^2}{n} + {O(n^{-4})}) \\
    &\lesssim \lambda_n(t) + \frac{\const_1\log^3{n}}{n} \lambda_n(t)^2 + n^{-2}.\label{eq:var-wt}
\end{align}
Plugging this back into \cref{eq:bound-circled-1-j2-1} and using the fact that $\lambda_n(t) \gtrsim \sqrt{n}$ {uniformly} on $J_2$, we get
\begin{align}
    \sup_{t \in J_2}\abs\big{\pr(W(t) = 0) - e^{-\lambda_n(t)}} 
    &\lesssim \sup_{t \in J_2}\qty(\frac{1}{\lambda_n(t)} + \frac{\const_1 \log^3{n}}{n} +\frac{1}{n^2\lambda_n(t)^2} + e^{-\lambda_n(t)}) \lesssim \frac{1}{\sqrt{n}}.\label{eq:bound-circled-1-j2}
\end{align}

\subsubsection*{\thesubsubsection~(iii--b). $t \in J_3$.}

We use {a} Poisson approximation \citep{chen1975poisson,barbour1992poisson}. For $Y_i$ given in \cref{eq:yi} and by definition of $\Om_i$ in \cref{eq:Omega-i}, it is easy to verify that
\begin{align}
    \E(Y_i) = 0\qc{} \Var(Y_i) = I_p, \qq{and} \Cov(Y_i, Y_j) = \E(Y_i Y_j\tr) = \frac{\sigma_{ij}^2}{n}(\theta_{ij}\theta_{ji}\tr + \theta_{ji}\theta_{ij}\tr) \quad \forall i \neq j.\label{eq:yi-moments}
\end{align}
Therefore, $Y_i \not\indep Y_j$, and consequently $B_i(t) \not\indep B_j(t)$ for all $i \neq j$. The Poisson approximation {derived} in \citep{arratia1989two} is not useful in {the present} situation where the dependency graph {is} fully connected. We use the variant in \citep[Theorem~2.C]{barbour1992poisson}.

\cbox{black!5}
{%
\begin{lemma}\label{lem:poisson-approximation}
    Let $W(t) = \sum_{i \in [n]} B_i(t)$ where $B_i(t) = \mathbb{1}\qty{ \norm{Y_i} > u_n(t) }$ for $Y_i$ given in \cref{eq:yi} and $u_n(t) = a_n t + b_n$. Then, for $M_n = \max_{i \in [n]} \norm{Y_i}$ and for all $t \in \R$,
    \begin{align}
        \abs\Big{ \pr\qty\big({M_n} \le u_n(t)) - e^{-\lambda_n(t)} } 
        \le \frac{1 - e^{-\lambda_n(t)}}{\lambda_n(t)} \qty( \sum_{i \in [n]} \pi_i(t)^2 + \sum_{i \neq j} \abs{\Cov\qty\Big(B_i(t), B_j(t))} ). \label{eq:poisson-approximation}
    \end{align}
\end{lemma}
}
See \cref{proof:lem:poisson-approximation} for the proof of \cref{lem:poisson-approximation}. From \cref{lem:tail-bound,lem:chisq-tail,eq:lambda-n-k2}, observe that $\pi_i(t) \sim \pr(\norm{Z} > u_n(t)) \sim \frac{1}{n}\lambda_n(t)$ for all $t \in J_3$; therefore
\begin{align}
    \sum_{i \in [n]}\pi_i(t)^2 \sim \frac{\lambda_n(t)^2}{n}.\label{eq:pi2-bound}
\end{align}
Similar to the steps in \cref{eq:var-wt}, we obtain
\begin{align}
    \sum_{i \neq j}\abs{\Cov(B_i(t), B_j(t))} 
    \lesssim n^2 \cdot \const_1 \log^3{n} \frac{\pr(\norm{Z} > u_n(t))^2}{n} \;{\sim}\; \const_1 \log^3{n} \cdot \frac{\lambda_n(t)^2}{n}.\label{eq:covariance-bound-in-t}
\end{align}
Using \cref{eq:pi2-bound,eq:covariance-bound-in-t} in \cref{lem:poisson-approximation} and by noting that $\lambda_n(t) \lesssim \sqrt{n}$ {uniformly} on $J_3$ leads to
\begin{align}
    \sup_{t \in J_3}\abs\Big{ \pr\qty\big({M_n} \le u_n(t)) - e^{-\lambda_n(t)} } 
    &\lesssim \sup_{t \in J_3}\frac{1 - e^{-\lambda_n(t)}}{\lambda_n(t)} \qty( \frac{\lambda_n(t)^2}{n} + \const_1\log^3{n} \frac{\lambda_n(t)^2}{n} )\\
    &\le (1 + \const_1\log^3{n}) \cdot \sup_{t \in J_3}\frac{\lambda_n(t)}{n} \lesssim \const_1\frac{\log^3{n}}{\sqrt{n}}.\label{eq:bound-circled-1-j3}
\end{align}
Combining the bounds in \cref{eq:bound-circled-1-j1,eq:bound-circled-1-j4,eq:bound-circled-1-j2,eq:bound-circled-1-j3}, we have
\begin{align}
    {\Circled{1} = \sup_{t \in J_1 \cup J_2 \cup J_3 \cup J_4} \abs{ \pr\qty({M_n} \le u_n(t)) - e^{-\lambda_n(t)} } \le \const_1\frac{\log^3{n}}{\sqrt{n}}.}\label{eq:final-bound-circled-1}
\end{align}

\subsubsection{Bound for $\Circled{2}$}
\label{proof:claim-circled-2}

Similar to the bound for $\Circled{1}$, let {$t_n := {\log\log{n}}$} and consider the following three intervals:
\begin{align}
    K_1 := (-\infty, -t_n),\qq{} 
    K_2 := \qty[-t_n, t_n],\qq{}
    K_3 := (t_n,\infty).
    \label{eq:intervals-circled-2}
\end{align}
Throughout, we will also use the fact that for all $\abs{t} < b_n \asymp \sqrt{\log{n}}$, from \cref{eq:lambda-n-k2} we have
\begin{align}
    \lambda_n(t) = e^{-t}\qty\big( 1 + \eta_n(t) ) \qq{where} \abs{\eta_n(t)} \lesssim \abs{\zeta_n(t)} + \tfrac{t^2}{b_n^2}, \label{eq:lambda-n-k2-eta}
\end{align}

{
where $\lesssim$ above only suppresses absolute constants, and $\sup_{t \in K_2}\abs{\eta_n(t)} = o(1)$ uniformly.
}

\subsubsection*{\thesubsubsection~(i). $t \in K_1$.}

For $t < -t_n$, we have {$e^{-t} > e^{t_n}$}. Similarly, from \cref{eq:lambda-n-j1} we have $\lambda_n(t) = n$ for $t \le -b_n^2$, and when $t \in (-b_n^2, -t_n)$, from \cref{eq:lambda-n-k2,eq:lambda-n-k2-eta},
\begin{align}
    \lambda_n(t) = e^{-t-t^2/2b_n^2}(1 + \zeta_n(t)) > e^{t_n}\qty\big(1 - \abs{\eta_n(t_n)}).
\end{align}
Because $\eta_n(t_n) = o(1)$, for sufficiently large $n$ we have $\abs{\eta_n(t_n)} < \frac{3}{4}$ and $\min\qty{\lambda_n(t), e^{-t}} > \frac{1}{4}e^{t_n}$. Therefore, for all $t \in K_1$,
\begin{align}
    \abs\big{e^{-\lambda_n(t)} - e^{-e^{-t}}} &\le e^{-\lambda_n(t)} + e^{-e^{-t}} {\lesssim e^{-\min\qty{\lambda_n(t), e^{-t}}} \lesssim e^{-\frac{1}{4}e^{t_n}} = \frac{1}{n^{1/4}}}.
    \label{eq:bound-circled-2-k1}
\end{align}
\subsubsection*{\thesubsubsection~(ii). $t \in K_3$.}
We use the fact that $z \mapsto e^{-z}$ is $1$-Lipschitz for $z \ge 0$ to get
\begin{align}
    \sup_{t \in K_3}\abs\big{e^{-\lambda_n(t)} - e^{-e^{-t}}} \le \sup_{t \in K_3}\abs\big{\lambda_n(t) - e^{-t}}.
    \label{eq:lipschitz-bound-circled-2}
\end{align}
For $\tau$ given in \cref{lem:tail-bound}, we further split $K_3 = (t_n, \sqrt{\log{n}}] \cup (\sqrt{\log{n}}, \tau\log{n}] \cup [\tau\log{n}, \infty)$.

$\bullet$\quad For all $t > \tau \log{n}$, we have $e^{-t} \le 1/n^{\tau}$ and from \cref{eq:lambda-n-j4} we have $\lambda_n(t) = O(1/n)$. It follows that 
$$
\sup_{t > \log{n}}\abs{\lambda_n(t) - e^{-t}} = O\qty(n^{-1} \vee n^{-\tau}).
$$

$\bullet$\quad Similarly, for $\sqrt{\log{n}} < t \le \tau\log{n}$, we have $e^{-t} \le e^{-\sqrt{\log{n}}}$. From \cref{eq:lambda-n-k2}, we also have  
$$
\sup_{\sqrt{\log{n}} < t \le \tau\log{n}}\zeta_n(t) = O(1)
$$ 
from which it follows that $\lambda_n(t) \lesssim e^{-t - t^2/2b_n^2} \lesssim e^{-\sqrt{\log{n}}}$ uniformly for all $\sqrt{\log{n}} < t \le \tau\log{n}$. Therefore,
$$
\sup_{\sqrt{\log{n}} < t \le \tau\log{n}}\abs{\lambda_n(t) - e^{-t}} \lesssim e^{-\sqrt{\log{n}}} = o(1/\log{n}).
$$

$\bullet$\quad On the other hand, for $t \in (t_n, \sqrt{\log{n}}]$, using \cref{eq:lambda-n-k2-eta} and the bound for $\zeta_n(t)$ from \cref{eq:lambda-n-k2},
\begin{align}
    \abs{\lambda_n(t) - e^{-t}} = e^{-t}\abs{\eta_n(t)} \lesssim e^{-t}\qty( \frac{\abs{t} + t^2 + \log\log{n}}{\log{n}} + \frac{\const_1\log^{3/2}{n}}{\sqrt{n}}).
\end{align}
Using the fact that $e^{-t} \le 1$ and $t e^{-t} \le 1/e$ and $t^2e^{-t} \le ({2/e})^2$ for all $t \ge 0$, we obtain
\begin{align}
    \sup_{t \in K_3}\abs{e^{-\lambda_n(t)} - e^{-e^{-t}}} &\le \sup_{t \in (t_n, \sqrt{\log{n}}]} \abs{\lambda_n(t) - e^{-t}} \lesssim \frac{\log\log{n}}{\log{n}}+ \frac{\const_1\log^{3/2}{n}}{\sqrt{n}}.\label{eq:bound-circled-2-k3}
\end{align}

\subsubsection*{\thesubsubsection~(iii). $t \in K_2$.}

We need a tighter bound for this step. From the mean value theorem,
\begin{align}
    \abs\big{e^{-\lambda_n(t)} - e^{-e^{-t}}} \le e^{-\min\qty{\lambda_n(t), e^{-t}}} \cdot \abs{\lambda_n(t) - e^{-t}}.\label{eq:mean-value-theorem-circled-2}
\end{align}
From \cref{eq:lambda-n-k2-eta}, note that $\abs{\lambda_n(t) - e^{-t}} = e^{-t}\abs{\eta_n(t)}$ where $\sup_{t \in K_2}\eta_n(t) = o(1)$. For sufficiently large~$n$ we have $\abs{\eta_n(t)} < 3/4$ from which it follows that $\min\qty{\lambda_n(t), e^{-t}} > e^{-t}/4$ for all $t \in K_2$. Plugging this back into \cref{eq:mean-value-theorem-circled-2}, we get
\begin{align}
    \abs\big{e^{-\lambda_n(t)} - e^{-e^{-t}}} \le e^{-e^{-t}/4} \cdot e^{-t} \cdot \abs{\eta_n(t)}.
\end{align}
Note that $e^{-t} \cdot e^{-e^{-t}/4} = z(t)e^{-z(t)/4}$ for $z(t) = e^{-t}$. Using the fact that for $z \ge 0$ the function $f(z) = ze^{-z/4}$ has a maximum value of $4/e$ at $z=4$, we get
\begin{align}
    \sup_{t \in K_2}\abs\big{e^{-\lambda_n(t)} - e^{-e^{-t}}} &\le \frac{4}{e} \cdot \sup_{t \in K_2}\abs{\eta_n(t)} \lesssim \frac{\log\log{n}}{\log{n}} + \const_1\frac{\log^{3/2}{n}}{\sqrt{n}}.\label{eq:bound-circled-2-k2}
\end{align}
Combining the bounds in \cref{eq:bound-circled-2-k1,eq:bound-circled-2-k2,eq:bound-circled-2-k3}, we have
\begin{align}
    {\Circled{2} = \sup_{t \in K_1 \cup K_2 \cup K_3} \abs{e^{-\lambda_n(t)} - e^{-e^{-t}}} \lesssim \frac{\log\log{n}}{\log{n}} + \const_1\frac{\log^{3/2}{n}}{\sqrt{n}}.}
\end{align}
The desired bound for \cref{thm:main} now follows from \cref{eq:claim-1-2} and \cref{eq:claim-mn}.
\qed

\clearpage
\subsection*{Proof of \cref{cor:iid}}
\label{proof:cor:iid}

{
Let $\gamma_n := \kappa^2\Rx^2\sqrt{\log{n}/n}$, and for $x_1, \dots, x_n \simiid{} F$ let $\mathcal{A}$ be the event given by
\begin{align}
  \mathcal{A} := \qty{ 
    \ttinf{X} \le \Rx \qq{and}
    \frac{1}{\kappa(1 + \gamma_n))}  
    \le s_p\qty( \frac{HX}{\sqrt{n}})
    \le s_1\qty( \frac{HX}{\sqrt{n}})
    \le \kappa(1 + \gamma_n)
  }.
\end{align}
where $s_k(A)$ is the $k$-th largest singular value of $A$. By splitting the probability $\pr\qty( (T_n - b_n)/a_n \le t )$ conditionally on $\mathcal{A}$ and $\mathcal{A}^c$, we have
\begin{align}
  &\abs{\pr\qty(  \frac{T_n - b_n}{a_n} \le t) - \pr(G \le t)}\\ 
  &\qquad\qquad\le \pr(\mathcal{A}) \cdot \abs{ \pr\qty(  \frac{T_n - b_n}{a_n} \le t \mid \mathcal{A}) - \pr(G \le t) }
  + \pr(\mathcal{A}^c) \cdot \abs{ \pr\qty(  \frac{T_n - b_n}{a_n} \le t \mid \mathcal{A}^c) - \pr(G \le t) }\\
  &\qquad\qquad\le \abs{ \pr\qty(  \frac{T_n - b_n}{a_n} \le t \mid \mathcal{A}) - \pr(G \le t) } + \pr(\mathcal{A}^c).\label{eq:cor-iid-proof-1}
\end{align}
On the event $\mathcal{A}$, note that assumption \ref{assumption:compact} holds with the same $\Rx$ but $\kappa$ replaced by $\kappa_n = \kappa(1 + \gamma_n)$. Moreover, we further have that $\kappa_n \le 2\kappa$ for sufficiently large $n$. Thus, conditional on $\mathcal{A}$, we can apply \cref{thm:main} to obtain
\begin{align}
  \sup_{t \in \R}\abs{\pr\qty(\frac{T_n - b_n}{a_n} \le t \mid \mathcal{A}) - \pr(G \le t)} \le C \frac{\log\log{n}}{\log{n}} + \const_1(p, 2\kappa, \Rx, \msigma, \Msigma) \frac{\log^3{n}}{\sqrt{n}}.
\end{align}
From Lemma~1 of \citet{vishwanath2025minimax}, we also have $\pr(\mathcal{A}^c) = O(n^{-2})$ for sufficiently large $n > N_0$.  Plugging these bounds back into \cref{eq:cor-iid-proof-1}, we obtain
\begin{align}
  \lim_{n \to \infty}\sup_{t \in \R}\abs{\pr\qty(  \frac{T_n - b_n}{a_n} \le t) - \pr(G \le t)} = 0,
\end{align}
which implies the result in \cref{eq:cor-iid}.
\qed
}


\subsection{Proof of Corollary~\ref{cor:confidence-pivotal}}
\label{proof:cor:confidence-pivotal}

{
    For $a_n, b_n$ given in \cref{eq:anbn}, let 
    \begin{align}
        \alpha_n := 1 - e^{-e^{b_n/a_n}} \qq{such that} q_{1-\alpha_n} = -{b_n}/{a_n}.
    \end{align} 
    Note that ${\lim_n\alpha_n = 1}$. Therefore, for all practical values of the confidence level, $\alpha < \alpha_n$, we have $b_n + a_n q_{1-\alpha} > 0$, and, from the definition of $\ellipse_{\alpha, i}$ in \cref{eq:confidence-set-i}, it follows that $\ellipse_{\alpha, i} \neq \emptyset$ for all $\alpha < \alpha_n$. Since $\hg\inv(v) = \hP\tr v$ from \cref{eq:hg}, we additionally have $\norm{ \Om_i^{-1/2}\hP\tr( \hg(x_i) - \hx_i ) } = \norm{\Om_i^{-1/2}(x_i - \hg\inv(\hx_i))}$, which implies that
    \begin{align}
        \pr\qty\Big( \hg(x_i) \in \ellipse_{\alpha, i}, \forall i \in [n] ) = \pr\qty\Big( T_n \le b_n + a_n q_{1-\alpha} ).
    \end{align}
    Therefore, using the fact that $\pr(G\le q_{1-\alpha}) = 1-\alpha$ and from \cref{thm:main}, we have
    \begin{align}
    &\sup_{\substack{\alpha \in (0, 1)\\\alpha < \alpha_n}}\abs\Big{\pr\qty\Big( \hg(x_i) \in \ellipse_{\alpha, i}, \forall i \in [n] ) - (1 - \alpha)}
    = \sup_{\substack{\alpha \in (0, 1)\\\alpha < \alpha_n}}\abs{\pr\qty\Big(\frac{T_n - b_n}{a_n} \le q_{1-\alpha}) - \pr(G \le q_{1-\alpha})} \lesssim \rate_n.
    \end{align}
    For $\alpha \ge \alpha_n$, note that $b_n + a_nq_{1-\alpha} \le 0$ which implies that $\ellipse_{\alpha, i} = \emptyset$ and $\pr\qty\big( \hg(x_i) \in \ellipse_{\alpha, i}, \forall i \in [n] ) = 0$. It follows that
    \begin{align}
        \sup_{\substack{\alpha \in (0, 1)\\\alpha \ge \alpha_n}}\abs\Big{\pr\qty\Big( \hg(x_i) \in \ellipse_{\alpha, i}, \forall i \in [n] ) - (1 - \alpha)} \le 1-\alpha_n = \exp(-\exp(b_n/a_n)) \ll \rate_n,
    \end{align}
    since $b_n/a_n \asymp \log{n}$. Combining the bounds for the two cases above gives the desired result.
}
\qed


\subsection{Proof of \cref{prop:plug-in}}
\label{proof:prop:plug-in}

The proof is based on establishing a bound similar to Lemma~A.7 of \citep{spokoiny2015bootstrap}, which applies only to Gaussian random vectors and is, therefore, not directly applicable to our setting. Instead, we  use \cref{lem:slutsky} directly after establishing the following bound.
\cbox{black!5}{%
{
\begin{lemma}\label{lem:omega-error}
    Let $\Om_i, \hOm_i \in \Rpp$ be the matrices defined in \cref{eq:Omega-i} and \cref{eq:hOm}, respectively, and let $\hP$ be as given in \cref{eq:hg}. Then, with probability greater than $1 - O(n^{-2})$,
    \begin{align}
        \max_{i \in [n]}\Opnorm\Big{\Omega_i^{-1/2} \;\hP\tr\hat\Omega_i\hP\;\Omega_i^{-1/2} - I_p} \lesssim C'_2(\pars) \frac{\log^2{n}}{\sqrt{n}} \label{eq:plug-in-claim-1}
    \end{align}
\end{lemma}
}
}

The proof of \cref{lem:omega-error} is in \cref{proof:lem:omega-error}. Let $\Psi_i := \Omega_i^{-1/2} \;\hP\tr\hat\Omega_i\hP\;\Omega_i^{-1/2}$ be the matrix in \cref{eq:plug-in-claim-1} and note that for any $i \in [n]$ and $x, y \in \Rp$,
\begin{align}
    (\hg(x) - y)\tr\hOm_i\inv(\hg(x) - y) \
    &= \qty\big{\hP(x - \hP\tr y)}\tr\; \hOm_i\inv\; \qty\big{\hP(x - \hP\tr y)} \\
    &= (x - \hP\tr y)\tr(\hP\tr\hOm_i\inv\hP)(x - \hP\tr y) \\
    &= (x_i - \hg\inv(y))\tr \Om_i^{-1/2}\Psi_i\inv \Om_i^{-1/2}(x_i - \hg\inv(y))
\end{align}
from which it follows that
\begin{align}
    \Norm\big{\hOm_i^{-1/2}(\hg(x_i) - \hx_i)} = \Norm\big{{\Psi_i^{-1/2}}\Om_i^{-1/2}(x_i - \hg\inv(\hx_i))} \quad \forall i \in [n].\label{eq:tn-hat-alt}
\end{align}
From the definition of $T_n$ in \cref{thm:main} and using \cref{eq:tn-hat-alt} in the definition of $\hat T_n$ in \cref{eq:tn-hat}, we have
\begin{align}
    \abs{ \qty(\frac{\h{T}_n - b_n}{a_n}) - \qty(\frac{{T}_n - b_n}{a_n}) } 
    &= \frac{\sqrt{n}}{a_n} \cdot \abs{ \max_{i \in [n]}\Norm\big{\Psi_i^{-1/2}\Om_i^{-1/2}(x_i - \hg\inv(\hx_i))} - \max_{i \in [n]}\Norm\big{\Om_i^{-1/2}(x_i - \hg\inv(\hx_i))} }\\
    &\le \frac{\sqrt{n}}{a_n} \max_{i \in [n]}\Norm\big{ (I_p - \Psi_i^{-1/2})\Om_i^{-1/2}(x_i - \hg\inv(\hx_i)) }\\
    &\hspace{-3em}\lesssim {\sqrt{n\log{n}}} \cdot \max_{i \in [n]}\Opnorm\big{I_p - \Psi_i^{-1/2}} \cdot \max_{i \in [n]}\Opnorm\big{\Om_i^{-1/2}} \cdot \max_{i \in [n]}\Norm{x_i - \hg\inv(\hx_i)}.\label{eq:tn-hat-tn-diff}
\end{align}
where the second inequality follows by two applications of the reverse triangle inequality (one for the $\ell_{\infty}$-norm and one for the $\ell_2$-norm), and final inequality follows since in $a_n \sim 1/\sqrt{2\log{n}}$. From \ref{bound-2} and \cref{eq:omega-bounds}, with probability greater than $1 - O(n^{-2})$,
\begin{align}
    \max_{i \in [n]}\Opnorm\big{ \Om_i^{-1/2} } \le \frac{2\kappa}{\msigma} 
    \qq{and}
    \max_{i \in [n]}\norm{ x_i - \hg(\hx_i) } \lesssim c_2(\pars) \sqrt{\frac{\log{n}}{n}}. \label{eq:plug-in-bounds-2}
\end{align}
Lastly, let $z_n := C'_2 \log^2{n}/\sqrt{n}$ be the r.h.s. of \cref{eq:plug-in-claim-1}. On the event $\qty{\max_{i \in [n]}\opnorm{\Psi_i - I_p} \le z_n}$, which by \cref{lem:omega-error} holds with probability greater than $1 - O(n^{-2})$, for all $i \in [n]$ we have
\begin{align}
    (1-z_n)I_p \preccurlyeq \Psi_i \preccurlyeq (1+z_n)I_p \implies (1+z_n)^{-1/2} I_p \preccurlyeq \Psi_i^{-1/2} \preccurlyeq (1-z_n)^{-1/2} I_p.
\end{align}
Moreover, for sufficiently large $n$, $z_n < 1/2$ and it follows that 
\begin{align}
    \max_{i \in [n]}\Opnorm\big{ I_p - \Psi_i^{-1/2} } \le 1 - (1-z_n)^{-1/2}   \le 2z_n.\label{eq:plug-in-bounds-3}
\end{align}
Plugging in \cref{eq:plug-in-bounds-2,eq:plug-in-bounds-3} into \cref{eq:tn-hat-tn-diff} we get that with probability greater than $1 - O(n^{-2})$,
\begin{align}
    \abs{ \qty(\frac{\h{T}_n - b_n}{a_n}) - \qty(\frac{{T}_n - b_n}{a_n}) } \lesssim
    {\const}_2(\pars) \frac{\log^3{n}}{n},\label{eq:tn-hat-tn-bound-final}
\end{align}
for ${\const}_2 := {C_2' c_2 \kappa }/{\msigma}$. From \cref{thm:main}, we know that $\ks( (T_n-b_n)/a_n, G ) \lesssim \rate_n$; we can now use \cref{lem:slutsky} along with the modulus of continuity $\omega_G(\epsilon) \le \epsilon$ from \cref{eq:omega-gumbel} to get:
\begin{align}
    \ks\qty(\frac{\h{T}_n - b_n}{a_n}, G) 
    &\lesssim \rate_n + {\const}_2\frac{\log^{3}{n}}{n} + \frac{1}{n^2}\\
    &\lesssim \rate_n + {\const}_2\frac{\log^{3}{n}}{n},\label{eq:ks-hat-tn}
\end{align}
which completes the proof of \cref{eq:gumbel-convergence-hat}. For the plug-in confidence set ${\conf}_\alpha = \prod_{i \in [n]}{\conf}_{\alpha, i}$, note that
\begin{align}
    \pr\qty\Big( \hg(x_i) \in {\conf}_{\alpha, i}, \forall i \in [n] ) 
    \! =\! \pr\qty( \max_{i \in [n]}\sqrt{n}\norm{\hOm_i^{-1/2}(\hg(x_i) - \hx_i)} \le a_n q_{1-\alpha} + b_n)
    \! =\! \pr\qty( \h{T}_n \le b_n + a_nq_{1-\alpha} ).
\end{align}
The proof for the coverage guarantee is now identical to the proof of \cref{cor:confidence-pivotal} in \cref{proof:cor:confidence-pivotal}.
\qed

\subsection{Proof of \cref{thm:multiplier-bootstrap}}
\label{proof:thm:multiplier-bootstrap}

Conditional on $\Eps$, notice that $D\b = \h{\Del} + \Eps\b$ satisfies the noisy realizable setting defined in \cref{noisy-setting} where the ``\textit{true}'' latent configuration is $\hX$ and the noise terms are $\eps\b_{ij} = r_{ij}e_{ij}$ where $\Var(\eps\b_{ij} \mid \eps_{ij}) = e_{ij}^2$ and $\norm{r_{ij} \cdot e_{ij}}_{\psi_1} = \abs{e_{ij}} \cdot \norm{r_{ij}}_{\psi_1}$ {since $e_{ij}$ are treated as fixed.} Conditionally on $\Eps$, this setup satisfies assumptions \ref{assumption:compact}\,\textit{\&}\,\ref{assumption:noise} with $\h{\Rx} := \ttinf{\hX},$
\begin{align}
  \h{\kappa} := \frac{1}{\sqrt{n}}\max\qty\big{s_1(\hX), s_p(\hX)\inv},\quad
  \h{\msigma} := \min_{i}\lambda_p\qty(\sum_{k}e_{ij}^2\h{u}_i\h{u}_i\tr),\quad
  \h{\Msigma} := \norm{Z}_{\psi_{_1}} \cdot \max_{i < j} \abs{e_{ij}}
\end{align}
where $s_k(\hX)$ are the singular values of $\hX$ and $Z \sim N(0, 1)$. The matrix $\hOm_i$ plays the same role for $\hX\b$ as $\Om_i$ does for $\hX$. {Let $\Sigma = (\sigma_{ij}^2)$ be the matrix of variances for $\Eps$, and let $\hat\Sigma = (e_{ij}^2)$ be the matrix of variances for $\Eps\b$. } Using \labelcref{bound-1,bound-2,bound-6,bound-7}, with probability at least $1 - O(n^{-2})$ over the randomness of $\Eps$, we have
\begin{align}
    s_1(\hX) &\le s_1(\hg(X)) + \opnorm{\hX - \hg(X)} \le \kappa\sqrt{n} + c_1 \\[0.35em]
    s_p(\hX) &\ge s_p(\hg(X)) - \opnorm{\hX - \hg(X)} \ge \kappa\inv\sqrt{n} - c_1 \\[0.35em]
    \ttinf{\hX} &\le \ttinf{\hg(X)} + \ttinf{\hX - \hg(X)} \le \Rx + c_2\sqrt{\log{n}/n}\\[0.35em]
    \max_{i < j} e_{ij}^2 &\le \max_{i < j} \sigma_{ij}^2 + \maxnorm{\h{\Sigma} - \Sigma} \le \Msigma^2 + \Msigma^2 \log^2{n},\\[0.15em]
    \min_{i}\lambda_p({\hU\tr \h{\Sigma}_i \hU\tr})
    &\ge \min_{i}\lambda_p({U\tr \Sigma_i U}) - \max_{i \in [n]}\opnorm{\hU\tr \h{\Sigma}_i \hU - U\tr \Sigma_i U}
    \ge \msigma^2 - c'_8\frac{\log^2{n}}{\sqrt{n}},
\end{align}
{where, in the last inequality, we used: $\lambda_p(A_i) \ge \lambda_p(B_i) - \opnorm{A_i - B_i} \ge \lambda_p(B_i) - \max_{i \in [n]}\opnorm{A_i - B_i}$ for all $i \in [n]$ by an application of Weyl's inequality.} Thus, with probability at least $1 - O(n^{-2})$ and for all sufficiently large $n$ we have
\begin{align}
    \h{\kappa} &= n^{-1/2}\max\{ s_1(\hX), s_p(\hX)\inv \} \le 2\kappa\\[0.25em]
    \h{\Rx} &= \ttinf{\hX} \le 2\Rx\\[0.25em]
    \h{\Msigma} &= \sqrt{\max_{i < j} e_{ij}^2} \le \Msigma\log{n},\\[0em]
    \h{\msigma} &= \min_{i}\lambda_p({\hU\tr \h{\Sigma}_i \hU\tr}) \ge \half\msigma.\label{eq:bootstrap-params}
\end{align}

With these bounds in hand, consider the following three statistics:
\begin{align}
  \h{T}_n = \max_{i \in [n]}\norm{\hOm_i^{-1/2}(\hg(x_i) - &\hx_i)},\qq{}\qq{}
  T\b_n = \max_{i \in [n]}\norm{\hOm_i^{-1/2}(\hx_i - \hgb\inv(\hx\b_i))},\qq{and}\\
      &T_n = \max_{i \in [n]}\norm{\Om_i^{-1/2}(x_i - \hg\inv(\hx_i))}.
\end{align}
$T_n$ was analyzed in \cref{thm:main} and $\h{T}_n$ in \cref{prop:plug-in}. In what follows, we will prove:
\begin{align}
  \sup_{t \in \R}\abs{ \pr\qty\big( T_n \le t ) - \pr\b\qty\big(T\b_n \le t) } \lesssim \rate\b_n.\label{eq:bootstrap-claim-1}
\end{align}
The desired bound containing $\pr(\h{T}_n \le t)$ in place of $\pr(T_n \le t)$ will then follow from an identical argument used to prove \cref{prop:plug-in}, i.e., using exactly the same arguments leading up to \cref{eq:tn-hat-tn-bound-final} and then applying \cref{lem:slutsky} yields the desired result:
\begin{align}
  \sup_{t \in \R}\abs{ \pr\qty\big( \h{T}_n \le t ) - \pr\b\qty\big(T\b_n \le t) }
  &\le \sup_{t \in \R}\abs{ \pr\qty\big( {T}_n \le t ) - \pr\b\qty\big(T\b_n \le t) } + {\const}_2\frac{\log^{3}{n}}{n} \;\; \lesssim \;\; \rate\b_n.
\end{align}
The proof of the claim in \cref{eq:bootstrap-claim-1} is based on the bound established for term $\Circled{1}$ in the proof of \cref{thm:main}. 
The high-level procedure is fairly standard \citep[Chapter~2]{shao2012jackknife}, and the outline is as follows:
\begin{enumerate}[itemsep=-0.25em]
  \item Obtain an intermediate approximation for the distribution of $T_n$ using the bound for $\Circled{1}$.
  \item Conditionally on $\Eps$, obtain an analogous approximation for $T\b_n$.
  \item Combine the two intermediate approximations and uncondition on $\Eps$.
\end{enumerate}
Let us briefly recall the intermediate quantities which appeared in the proof of \cref{thm:main}.

\subsubsection*{Step 1. Intermediate bound for $T_n$.}

From \cref{prop:decomposition}, we had $\sqrt{n}\Om_i^{-1/2}(x_i - \hg\inv(\hx_i)) = Y_i + R_i$ where
\begin{align}
  Y_i = \frac{1}{\sqrt{n}}\sum_{k}\eps_{ik}\theta_{ik} 
  \qq{and}
  \pr\qty\Big( \textstyle{\max\limits_{i\in [n]}\norm{R_i} > C_1' \sqrt{\frac{\log{n}}{n}}} ) = O(n^{-2}),
\end{align}
for a constant $C'_1 \equiv C'_1(\pars)$.  Moreover, for $u_n(t) = a_n t + b_n$, 
\begin{align}
  M_n := \max_{i \in [n]}\norm{Y_i} \qc{} \lambda_n(t) := \sum_{i \in [n]} \pr\qty\big( \norm{Y_i} > u_n(t) ),\qq{and} \mathcal{F}_n(t) := e^{-\lambda_n(t)},
\end{align}
and from the bound on $\Circled{1}$ in \cref{eq:final-bound-circled-1}, we have
\begin{align}
  \sup_{t \in \R}\abs{ \pr\qty\big( M_n \le u_n(t) ) - \mathcal{F}_n(t) } \le \const_1\frac{\log^3{n}}{\sqrt{n}}.\label{eq:Mn-bound-bootstrap}
\end{align}
Here, $\mathcal{F}_n(t)$ is the c.d.f. of an intermediate tight sequence which was shown in \cref{proof:claim-circled-2} to converge to the c.d.f. of the Gumbel distribution. Notice that $\lambda_n(t)$ is a non-increasing function of $t$ and hence $\mathcal{F}_n(t)$ is non-decreasing and right-continuous. Now, using \cref{eq:Mn-bound-bootstrap} along with \cref{eq:remainder-term-mn} in \cref{lem:slutsky} gives
\begin{align}
  \sup_{t \in \R}\abs{ \pr\qty\big( T_n \le u_n(t) ) - \mathcal{F}_n(t) } &\lesssim \const_1\frac{\log^3{n}}{\sqrt{n}} + \omega_{n}\qty(C_1' \sqrt{\frac{\log{n}}{n}}),\label{eq:tn-intermediate}
\end{align}
where $\omega_{n}(\eta) := \sup\qty{ \mathcal{F}_n(t+\eta) - {\mathcal{F}_n}(\eta): t \in \R }$ is the modulus of continuity of $\mathcal{F}_n$. This gives us the desired intermediate approximation for $T_n$.

\subsubsection*{Step 2. Intermediate bound for $T\b_n$.}

We use $\h{\const}_{\square}, \h{C}_{\square}$ to denote the same constants as earlier but with ${\h{\Rx}, \h{\kappa}}$,~etc. in place of in place of $\Rx, \kappa$, etc. Conditional on $\Eps$, since $D\b = \h{\Del} + \Eps\b$ satisfies the assumptions \ref{assumption:compact}\,\textit{\&}\,\ref{assumption:noise} with parameters ${\h{\Rx}, \h{\kappa}}, \h{\msigma}, \h{\Msigma}$. Using \cref{prop:decomposition}, it follows that
$$
\sqrt{n}\hOm_i^{-1/2}(\hx_i - \hgb\inv(\hx\b_i)) = Y\b_i + R\b_i
$$ 
where, for $\theta\b_{ik} = \half\hOm_i^\minushalf \qty(\tfrac{\hX\tr \hX}{n})\inv \hx_k$, 
\begin{align}
  Y\b_i &:= \frac{1}{\sqrt{n}}\sum_{k}r_{ik}e_{ik}\theta\b_{ik},\qq{and}
  \pr\b\qty\Big( \textstyle{\max\limits_{i\in [n]}\norm{R_i\b} > \h{C}'_1 \sqrt{\frac{\log{n}}{n}}} ) = O(n^{-2}),
\end{align}
{where $\pr\b$ is with respect to the randomness in $R=(r_{ij})$ and $O(n^{-2})$ holds uniformly on $\Eps$.} For the same sequence $u_n(t) = a_n t + b_n$, we can similarly define
\begin{align}
  M\b_n &= \max_{i \in [n]}\norm{Y\b_i},\quad
  \lambda\b_n(t) = \sum_{i \in [n]} \pr\b\qty\big( \norm{Y\b_i} > u_n(t) ),\qq{and}
  \mathcal{F}\b_n(t) = e^{-\lambda\b_n(t)}.
\end{align}
Using the same bound for $\Circled{1}$ in \cref{eq:final-bound-circled-1} (applied conditionally on $\Eps$) followed by an application of \cref{lem:slutsky} (again, applied conditionally on $\Eps$) similar to \cref{eq:tn-intermediate} gives
\begin{align}
  \sup_{t \in \R}\abs{ \pr\b\qty\big( T\b_n \le u_n(t) ) - \mathcal{F}\b_n(t) } &\lesssim \h{\const}_1\frac{\log^3{n}}{\sqrt{n}} + \omega\b_{n}\qty(\h{C}'_1 \sqrt{\frac{\log{n}}{n}}),\label{eq:tnb-intermediate}
\end{align}
where $\omega\b_{n}(\eta) := \sup\qty{ \mathcal{F}\b_n(t+\eta) - \mathcal{F}\b(\eta): t \in \R }$ is the modulus of continuity of $\mathcal{F}\b_n$.

\subsubsection*{Step 3. Combining the intermediate approximations for $T_n$ and $T\b_n$.}

Using the intermediate approximations above and by noting that $t \mapsto u_n(t)$ is bijective, we can write \cref{eq:bootstrap-claim-1} as
\begin{align}
  &\sup_{t \in \R}\abs{ \pr\qty\Big( T_n \le t) - \pr\qty\Big( T\b_n \le t) }\\[0.5em]
  &\qquad \le \sup_{t \in \R}\abs{ \pr\qty\Big( T_n \le u_n(t)) - \pr\qty\Big( T\b_n \le u_n(t)) }\\[0.5em]
  &\qquad \le \sup_{t \in \R}\abs{ \pr\qty\Big( T_n \le u_n(t)) - \mathcal{F}_n(t) }
  + \sup_{t \in \R}\abs{ \mathcal{F}_n(t) - \mathcal{F}\b_n(t) }
  + \sup_{t \in \R}\abs{ \pr\qty\Big( T\b_n \le u_n(t)) - \mathcal{F}\b_n(t) }\\[0.5em]
  &\qquad \lesssim \sup_{t \in \R}\abs{ \mathcal{F}_n(t) - \mathcal{F}\b_n(t) } + \qty{\qty\big(\const_1 + \h{\const}_1)\frac{\log^3{n}}{\sqrt{n}} + \omega_n\qty(C'_1 \sqrt{\frac{\log{n}}{n}}) + \omega\b_n\qty(\h{C}'_1 \sqrt{\frac{\log{n}}{n}})}.\label{eq:tn-tnb-triangle}
\end{align}
In order to bound the {first} term, from \cref{lem:tail-bound} we have that,
\begin{align}
  \pr\qty(\norm{Y_i} > u_n(t)) = 
  \begin{dcases}
      1 & \qq{if } t \in (-\infty, -b_n^2]\\
      \pr\qty\big(\norm{Z} > u_n(t)) \qty( 1 + O\qty(\tfrac{\const_1 \log^{3/2}{n}}{\sqrt{n}}))  & \qq{if } t \in (-b_n^2, \tau \log{n}]\\
      O(1/n^2) & \qq{if } t \in (\tau\log{n}, \infty),
  \end{dcases}\label{eq:multiplier-tail-bound-yi}
\end{align}
and there exists $\hat\tau$ such that,
\begin{align}
  \pr\b\qty(\norm{Y_i\b} > u_n(t)) = 
  \begin{dcases}
      1 & \qq{if } t \in (-\infty, -b_n^2]\\
      \pr\qty\big(\norm{Z} > u_n(t)) \qty( 1 + O\qty(\tfrac{\h{\const}_1 \log^{3/2}{n}}{\sqrt{n}}))  & \qq{if } t \in (-b_n^2, \hat\tau \log{n}]\\
      O(1/n^2) & \qq{if } t \in (\hat\tau\log{n}, \infty).
  \end{dcases}
\end{align}
From \cref{proof:thm:main}, we have that $\const_1(\pars) = \const_0(\kappa, \Rx, \msigma)^2 \cdot \Msigma^2$. Using \cref{eq:bootstrap-params}, it follows that with probability at least $1 - O(n^{-2})$ over the randomness of $\Eps$,
$$
\h{\const}_1 \le \const_1\qty\big(2\kappa, 2\Rx, \Msigma\log{n}, \msigma/2) \lesssim \const_1(\pars)\log^2{n};
$$ 
similarly $\h{C}'_1 \lesssim C'_1\log^2{n}$ and for $\tau = 4\tilde{C}_p \Msigma \const_0$ from \cref{eq:bernstein-type-bound}, we have $\hat\tau \lesssim \tau\log{n}$. It follows that with probability at least $1 - O(n^{-2})$ over the randomness of $\Eps$,
\begin{align}
  \abs{ \pr\qty\big(\norm{Y_i} > u_n(t)) - \pr\b\qty\big(\norm{Y_i\b} > u_n(t)) } =
  \begin{dcases}
    0 & \text{if } t < -b_n^2\\
    \pr\qty\big( \norm{Z} > u_n(t) ) \cdot O\qty(\const_1\frac{\log^{7/2}{n}}{\sqrt{n}}) & \text{if } t \in (-b_n^2, \tau \log^2{n}]\\
    O(1/n^2) & \text{if } t > \tau\log^2{n}.\\
  \end{dcases}
\end{align}
Using the fact that $z \mapsto e^{-z}$ is $1$-Lipschitz for $z \ge 0$, we have
\begin{align}
  \abs{ \mathcal{F}_n(t) - \mathcal{F}\b_n(t) } 
  &= \abs{e^{-\lambda_n(t)} - e^{-\lambda\b_n(t)}} \\
  &\le \abs{\lambda_n(t) - \lambda\b_n(t)} \\
  &\le \sum_{i \in [n]} \abs{ \pr\qty(\norm{Y_i} > u_n(t)) - \pr\b\qty(\norm{Y_i\b} > u_n(t)) }.\label{eq:fn-fnb}
\end{align}
Note that \cref{eq:fn-fnb} is $0$ for $t < -b_n^2$ and is $O(1/n)$ for $t > \tau'\log{n}$. Since $\pr\qty(\norm{Z} > u_n(t)) \sim 1/n$ on $(-b_n^2, \tau\log^2{n})$ we get
\begin{align}
  \sup_{t \in \R}\abs{ \mathcal{F}_n(t) - \mathcal{F}\b_n(t) } \lesssim \const_1\frac{\log^{7/2}{n}}{\sqrt{n}}.\label{eq:fn-fnb-bound}
\end{align}

Plugging in the bound from \cref{eq:fn-fnb-bound} into \cref{eq:tn-tnb-triangle} gives
\begin{align}
    \sup_{t \in \R}&\abs{ \pr\qty\Big( T_n \le t) - \pr\qty\Big( T\b_n \le t) }
    \lesssim \const_1\frac{\log^{7/2}{n}}{\sqrt{n}} + \const_1\frac{\log^5{n}}{\sqrt{n}} + \omega_n\qty(C'_1 \sqrt{\frac{\log{n}}{n}}) + \omega\b_n\qty({C}'_1 \frac{\log^{5/2}{n}}{\sqrt{n}}).\label{eq:tn-tnb-final}
\end{align}

In order to complete the proof for \cref{eq:bootstrap-claim-1}, we need to establish a bound on $\omega_n(\eta)$ and $\omega\b_n(\eta)$. To this end, consider $\mathcal{F}_n(t) = \exp(-\lambda_n(t))$. Intuitively, $\omega_n(\eta)$ cannot be too different from $\omega_G(\eta)$ in \cref{eq:omega-gumbel} since $\mathcal{F}_n(t)$ converges to the Gumbel c.d.f. as $n \to \infty$. A multiplicative bound on $\omega_n(\eta)$ will suffice for our purposes.

Since we are evaluating $\omega_n(\eta)$ for the r.h.s. of \cref{eq:tn-tnb-final}, we may assume, without loss of generality, that $\eta < 1$ which holds for all sufficiently large $n$. From the discussion in Step 3 and \cref{eq:multiplier-tail-bound-yi} we have that $\mathcal{F}_n(t+\eta) - \mathcal{F}_n(t) = 0$ for all $t < -b_n^2$ and $\mathcal{F}_n(t+\eta) - \mathcal{F}_n(t) = O(\eta/{n})$ for all $t > \tau' \log{n}$. Therefore, it suffices to consider $t \in [-b_n^2, \tau' \log{n}]$. Let $\bar{F}_p(t) = \pr\qty( \norm{Z}^2 > u_n(t)^2 )$ denote the upper tail of a $\chi^2_p$ distribution and $f_p$ its density, and from \cref{eq:multiplier-tail-bound-yi} we have $\lambda_n(t) = n \cdot \bar{F}_p(u_n(t)^2) \cdot (1+s_n)$ for $s_n = O(\log^{3/2}{n}/\sqrt{n})$. Taking the derivative w.r.t. $t$ gives
\begin{align}
  \frac{d}{dt}\mathcal{F}_n(t) = -e^{-\lambda_n(t)} \cdot \frac{d}{dt}\lambda_n(t) 
  &= -e^{-\lambda_n(t)} \cdot n (1+ s_n) \cdot \frac{d}{dt}\bar{F}_p(u_n(t)^2)\\
  &= e^{-\lambda_n(t)} \cdot n(1 + s_n) \cdot f_p(u_n(t)^2) \cdot 2u_n(t)\frac{d}{dt}u_n(t).
\end{align}
Note that $\frac{d}{dt}u_n(t) = a_n = 1/b_n$ and $a_n u_n(t) = 1+ t/b_n^2$. Using the inverse Mills' ratio bound derived in~\cref{eq:mills-ratio}, i.e., $f_p(z) \le \half \bar{F}_p(z)$, we get
\begin{align}
  \frac{d}{dt}\mathcal{F}_n(t) \le {a_n u_n(t)} e^{-\lambda_n(t)} \cdot n(1+s_n) \cdot \bar{F}_p(u_n(t)^2) = (1 + t/b_n^2) \cdot \lambda_n(t)e^{-\lambda_n(t)}.
\end{align}
Using the fact that $ze^{-z} \le 1/e$ for all $z > 0$ and since $\abs{t}/b_n^2 \asymp 1/\sqrt{\log{n}}$ for all $t \in [-b_n^2, \tau' \log{n}]$, we have $\frac{d}{dt}\mathcal{F}_n(t) \le \frac{2}{e}$ {for all $n$ sufficiently large.} We therefore have that
\begin{align}
  \omega_n(\eta) \le \eta \cdot \sup_{t \in \R} \frac{d}{dt}\mathcal{F}_n(t) \le \frac{2}{e}\eta.\label{eq:omega-n-bound}
\end{align}
An identical analysis also gives $\omega\b_n(\eta) \le (2/e)\eta$. Finally, plugging the bound in \cref{eq:omega-n-bound} into \cref{eq:tn-tnb-final}, we obtain
\begin{align}
    \sup_{t \in \R}\abs{ \pr\qty\Big( T_n \le t) - \pr\qty\Big( T\b_n \le t) } &\lesssim  \const_1\frac{\log^5{n}}{\sqrt{n}},
\end{align}
with probability at least $1 - O(n^{-2})$ over the randomness of $\Eps$. This proves the claim in \cref{eq:bootstrap-claim-1}. The proof of the final result is then completed by the argument following \cref{eq:bootstrap-claim-1}.
\qed

\subsection{Proof of \cref{cor:confidence-multiplier-bootstrap}}
\label{proof:cor:confidence-multiplier-bootstrap}

Note that
\begin{align}
    \pr\qty\Big( \hg(x_i) \in \conf\b_{\alpha, i}, \; \forall i \in [n] ) 
    &= \pr\qty( (\hg(x_i) - \hx_i)\tr\hOm_i\inv(\hg(x_i) - \hx_i) \le \frac{(q\b_{1-\alpha})^2}{n},  \; \forall i \in [n])\\
    &= \pr\qty( \max_{i \in [n]}\sqrt{n}\norm{\hOm_i^{-1/2}(\hg(x_i) - \hx_i)} \le q\b_{1-\alpha}) = \pr\qty( \h{T}_n \le q\b_{1-\alpha} ).
\end{align}
Since $q\b_{1-\alpha}$ satisfies $1-\alpha \le \pr\b(T\b_n \le q\b_{1-\alpha})$, from \cref{thm:multiplier-bootstrap} it follows that for all $\alpha \in (0, 1)$ and with probability greater than $1 - O(n^{-2})$ over the randomness of $\Eps$,
\begin{align}
    \pr\qty\big( \h{T}_n \le q\b_{1-\alpha} ) \ge \pr\b\qty\big(T\b_n \le q\b_{1-\alpha}) - O(\rate\b_n) \ge (1 - \alpha) - O(\rate\b_n).\label{eq:coverage-bootstrap-lower}
\end{align}
On the other hand, for the reverse inequality, we need the intermediate bounds obtained in \cref{eq:tn-intermediate,eq:tnb-intermediate,eq:fn-fnb} from \cref{proof:thm:multiplier-bootstrap}. For any $\eta > 0$, $\alpha \in (0,1)$ and with probability greater than $1 - O(n^{-2})$ over the randomness of $\Eps$, we have
\begin{align}
    \pr\qty\big(\h{T}_n \le q\b_{1-\alpha}) 
    &\le \mathcal{F}_n( q\b_{1-\alpha} ) + O(\rate\b_n)\\
    &\le \mathcal{F}_n( q\b_{1-\alpha} - \eta ) + \omega_n(\eta) + O(\rate\b_n)\\
    &\le \mathcal{F}\b_n( q\b_{1-\alpha} - \eta ) + \omega_n(\eta) + O(2\rate\b_n)\\
    &\le \pr\b( q\b_{1-\alpha} - \eta ) + \omega_n(\eta) + O(3\rate\b_n) 
    \;\;<\;\; (1-\alpha) + \omega_n(\eta) + O(3\rate\b_n),
\end{align}
where the first, third and fourth inequalities follow from \cref{eq:tn-intermediate,eq:fn-fnb,eq:tnb-intermediate} respectively. The second inequality follows from the definition of the modulus of continuity $\omega_n(\eta)$ of $\mathcal{F}_n$, and the final inequality holds since $q\b_{1-\alpha} = \inf\{t \in \R: \pr\b(T\b_n \le t) \ge 1 - \alpha\}$, and, therefore, ${\pr\b(T\b_n \le q\b_{1-\alpha} - \eta) < 1 - \alpha}$ for all $\eta > 0$. 

Choosing $\eta = O(\rate\b_n)$ and using the bound $\omega_n(\eta) \le \eta$ from \cref{eq:omega-n-bound} gives
\begin{align}
    \pr\qty\big(\h{T}_n \le q\b_{1-\alpha}) - (1-\alpha) \lesssim {\rate\b_n}.\label{eq:coverage-bootstrap-upper}
\end{align}
Combining \cref{eq:coverage-bootstrap-lower,eq:coverage-bootstrap-upper} gives the desired result.
\qed


\subsection{Proof of \cref{prop:iid}}
\label{proof:prop:iid}

\cref{prop:iid} is a combination of \cref{thm:main} and \cref{prop:plug-in} to the setting where $(\eps_{ij})$ {are} \iid{} Specifically, note that since $\Sigma_{i} = \diag\qty(\sigma_{i1}^2, \dots, \sigma_{in}^2) = \sigma^2 I_p$, the matrix $\Om_i = \Om$ for all $i \in [n]$ where
\begin{align}
  \Om = \frac{\sigma^2}{4} \qty(\tfrac{X\tr X}{n})\inv,
  \qq{and}
  \Om^{-1/2}(x_i - \hg\inv(\hx_i)) = \frac{2}{\sigma}\qty(\tfrac{X\tr X}{n})^{-1/2}(x_i - \hg\inv(\hx_i)).\label{eq:omega-iid}
\end{align}
Let
\begin{align}
  T'_n := \max_{i \in n}\frac{2\sqrt{n}}{\sigma}\Norm{  \qty(\tfrac{X\tr X}{n})^{-1/2}(x_i - \hg\inv(\hx_i))}.
\end{align}
Applying \cref{thm:main} for $T'_n$ now gives
\begin{align}
  \ks\qty( \frac{T'_n - b_n}{a_n}, G ) \lesssim \frac{\log\log{n}}{\log{n}} + \const_1(\kappa, \Rx, \sigma, \varsigma) \frac{\log^3{n}}{\sqrt{n}}.\label{eq:tn-gumbel}
\end{align}
For $\tilde{T}_n$, as in the proof of \cref{prop:plug-in} in \cref{proof:prop:plug-in}, define 
\begin{align}
  \hOm = \frac{\hat\sigma^2}{4}\qty(\tfrac{\hX\tr\hX}{n})\inv \qq{and} \Psi := \Om^{-1/2}\;\hP\tr\hOm\hP\;\;\Om^{-1/2}\label{eq:psi-iid}
\end{align}
such that, using the same arguments leading up to \cref{eq:plug-in-bounds-2}, we have
\begin{align}
  \norm{\hOm^{-1/2}(\hg(x_i) - \hx_i)} = \norm{\Psi^{-1/2}\Om^{-1/2}(x_i - \hg\inv(\hx_i))} \qquad \forall\; i \in [n],\label{eq:tn-hat-alt-iid}
\end{align}
and, once again, by noting that $a_n \sim 1/\sqrt{2\log{n}}$, it follows that
\begin{align}
  \abs\big{ (\tilde{T}_n-b_n)/a_n - (T'_n-b_n)/a_n } 
  &\lesssim {\sqrt{n\log{n}}} \cdot \Opnorm\big{ I - \Psi^{-1/2} } \cdot \Opnorm\big{\Om^{-1/2}} \cdot \max_{i}\Norm{ x_i - \hg\inv(\hx_i) },\label{eq:tn-tnprime-diff}
\end{align}
From \ref{assumption:compact} and from \ref{bound-2}, with probability greater than $1 - O(n^{-2})$,
\begin{align}
  \Opnorm\big{\Om^{-1/2}} \le \frac{2\kappa}{\sigma}\qq{and} \max_{i}\Norm{ x_i - \hg\inv(\hx_i) } \lesssim c_2(\kappa, \Rx, \sigma, \varsigma) \cdot \sqrt{\frac{\log{n}}{n}}.\label{eq:omega-bounds}
\end{align}
Therefore, a bound on $\Opnorm\big{ \Psi - I }$ will give a bound on the r.h.s. of \cref{eq:tn-tnprime-diff} exactly as in the proof of \cref{prop:plug-in}.

\cbox{black!5}
{%
\begin{lemma}\label{lem:psi-bound-iid}
    Let $\Psi \in \Rpp$ be as defined in \cref{eq:psi-iid}. Then, with probability greater than $1 - O(n^{-2})$,
    \begin{align}
        \Opnorm\big{ I - \Psi^{-1/2} } \lesssim C_3'(\kappa, \Rx, \sigma, \varsigma) \cdot n^{-1/2}.
    \end{align}
\end{lemma}
}
The proof of \cref{lem:psi-bound-iid} is given in \cref{proof:lem:psi-bound-iid} below. The rest of the proof now follows exactly as in the proof of \cref{prop:plug-in} in \cref{proof:prop:plug-in}. Specifically, plugging in the bounds from \cref{eq:theta-norm} and \cref{lem:psi-bound-iid} into \cref{eq:tn-tnprime-diff}, with probability greater than $1 - O(n^{-2})$ we have
\begin{align}
  \abs\big{ (\tilde{T}_n-b_n)/a_n - (T'_n-b_n)/a_n } \lesssim \const_3(\kappa, \Rx, \sigma, \varsigma) \cdot \frac{\log{n}}{\sqrt{n}}.\label{eq:tn-tnprime-final}
\end{align}
where $\const_3 = C_3'c_2\kappa/\sigma$. Combining \cref{eq:tn-gumbel} with \cref{eq:tn-tnprime-final} and using \cref{lem:slutsky} with $\omega_G(\eta)\le \eta$,
\begin{align}
    \ks\qty(\frac{\tilde{T}_n - b_n}{a_n}, G) 
    &\lesssim \ks\qty( \frac{T'_n - b_n}{a_n}, G ) + {\const}_3\frac{\log{n}}{n} + \frac{1}{n^2}
    \lesssim \frac{\log\log{n}}{\log{n}} + \const_1 \frac{\log^3{n}}{\sqrt{n}} + {\const}_3\frac{\log{n}}{n},\label{eq:ks-hat-tn}
\end{align}
which gives the stated result by noting that $\const_3\log{n} \ll \const_1\log^3{n}$.
\qed

\begin{algorithm}[t]
    \small
    \caption{\small Nonparametric Bootstrap Confidence Sets for Noisy MDS with \iid{} Noise}
    \label{alg:empirical-bootstrap}
    \begin{algorithmic}[1]
        \Require Dissimilarity matrix $D \in \R^{n \times n}$, embedding dimension $p$,\\ number of bootstrap samples $B$, nominal level $\alpha \in (0, 1)$
        \Statex
        \State Compute $\hX \gets \mds(D, p)$ 
        \State Compute $E \gets D - \Del(\hX)$
        \State Set $\hOm \gets (\h{\sigma}/4) \cdot (\hX\tr\hX/n)^{-1}$
        \For{$b = 1$ to $B$} 
            \State Sample $\eps_{ij}\s$ from $E$ with replacement for $i < j$ \Comment{Nonparametric bootstrap}
            \State Generate noisy dissimilarities $D\s \gets \Del(\hX) + \Eps\s$ where $\Eps\s = (e\s_{ij})$
            \State Compute $\hX\s \gets \mds(D\s, p)$ \Comment{Bootstrap embedding}
            \State Solve $\hP\s$ via orthogonal Procrustes analysis using \cref{eq:procrustes-bootstrap}
            \State Transform $\hgs\inv(\hX\s) = \hX\s\hP\s$ \Comment{Rigid transformation}
            \State $T\s_n(b) \gets \max_{i \in [n]} \sqrt{n}\norm{\hOm^{-1/2}(\hx_i - \hgs\inv(\hx_i\s))}$ \Comment{Bootstrap statistic}
        \EndFor
        \State Set $q_{1-\alpha}\s \gets$ the $(1-\alpha)$-quantile of $\qty{T\s_n(1), \dots, T\s_n(B)}$
        \State Compute $\conf\s_{\alpha, i}$ for each $i \in [n]$ using \cref{eq:confidence-set-i-empirical}.
        \State \Return Confidence sets $\conf_{\alpha}\s = \prod_{i=1}^n \conf_{\alpha, i}\s$
    \end{algorithmic}
\end{algorithm}

\subsection{Proof of \cref{thm:empirical-bootstrap}}
\label{proof:thm:empirical-bootstrap}

As in the proof of \cref{prop:iid}, let $\Om$ and $\hOm$ be the matrices defined in \cref{eq:omega-iid} and \cref{eq:psi-iid}, respectively, i.e.,
\begin{align}
    \Om = \frac{\sigma^2}{4}\qty(\tfrac{X\tr X}{n})\inv, \qq{and} \hOm = \frac{\hat\sigma^2}{4}\qty(\tfrac{\hX\tr\hX}{n})\inv.
\end{align}
The proof is identical to the proof of \cref{thm:multiplier-bootstrap} in \cref{proof:thm:multiplier-bootstrap} with one modification in the last step, which we outline below. As in \cref{proof:thm:multiplier-bootstrap}, define
\begin{align}
    T'_n = \max_{i \in [n]}\norm{\Om^{-1/2}(x_i - \hg\inv(\hx_i))},\quad
    \tilde{T}_n = \max_{i \in [n]}\norm{\hOm^{-1/2}(\hg(x_i) - \hx_i)},\quad
    T\s_n = \max_{i \in [n]}\norm{\hOm^{-1/2}(\hx_i - \hgs\inv(\hx\s_i))}.
\end{align}
Conditional on $\Eps$, let $\pr\s := \binom{n}{2}\inv \sum_{i < j}\dirac\qty{e_{ij} - \bar{e}}$ be the empirical measure on $E$, where $\dirac\qty{e_{ij} - \bar{e}}$ is a Dirac measure at $(e_{ij}-\bar{e})$. For $\xi \sim \pr\s$, note that $\E\s(\xi) = 0$, $\Var\s(\xi) = \hat\sigma^2$ and
\begin{align}
    \norm{\xi}\s_{\psi_1} 
    &= \sup\qty{ t > 0: \E\s( \exp(\abs{\xi}/t) ) \le 2 }\\
    &= \sup\qty{ t > 0: {\binom{n}{2}\inv}\sum_{i < j}\exp( \abs{e_{ij} - \bar{e}} / t ) \le 2  } \le \max_{i < j}\abs{e_{ij} - \bar{e}}  < \infty,
\end{align}
and $\pr\s(\xi = e_{ij} -\bar{e}) = \binom{n}{2}\inv$. Therefore, the noise terms $\eps\s_{ij}$ are sampled \iid{} from the empirical measure $\pr\s$, and $D\s = \h{\Del} + \Eps\s$ satisfies the noisy realizable setting defined in \cref{noisy-setting} where the ``\textit{true}'' latent configuration is $\hX$. Conditionally on $\Eps$, this setup satisfies the assumptions of \cref{prop:iid} with
\begin{align}
    \h{\Rx} := \ttinf{\hX},\quad
    \h{\kappa} := \frac{1}{\sqrt{n}}\max\qty\big{s_1(\hX), s_p(\hX)\inv},\quad
    \h{\sigma}^2 := \frac{1}{n}\sum_{i < j}^n(e_{ij} - \bar{e})^2,\quad
    \h{\varsigma} := \max_{i < j} \abs{e_{ij} - \bar{e}},
\end{align}
where $s_k(\hX)$ are the singular values of $\hX$. Let $\h{\const}_1, C'_3$  denote the constant in \cref{prop:iid} with $\h{\Rx}, \h{\kappa}$, $\h{\sigma}, \h{\varsigma}$ in place of $\Rx, \kappa, \sigma, \varsigma$. Conditional on $\Eps$, $T\s_n$ mimics the distribution of $T'_n$. The arguments leading up to \cref{eq:tn-tnb-final} in \cref{proof:thm:multiplier-bootstrap} hold and gives
\begin{align}
    \sup_{t \in \R}\abs{ \pr\qty\Big( T'_n \le t) - \pr\s\qty\Big( T\s_n \le t) } &\lesssim \qty( \const_1 + \h{\const}_1 )\frac{\log^3{n}}{\sqrt{n}} + \omega_n\qty(C'_1 \sqrt{\frac{\log{n}}{n}}) + \omega\s_n\qty(\h{C}'_1 \sqrt{\frac{\log{n}}{n}}),\label{eq:tn-tns-final}
\end{align}
where $\omega_n(\eta), \omega\s_n(\eta) \le \eta$ follow from \cref{eq:omega-n-bound}. The only remaining step is to bound the $\h{\const}_1(\h{\kappa} ,\h{\omega}, \h{\sigma}, \h{\varsigma})$ by unconditioning on $\Eps$. From \cref{eq:bootstrap-params}, with probability greater than $1 - O(n^{-2})$ we have
\begin{align}
    \h{\Rx} &= \ttinf{\hX} \le 2\Rx\\[0.25em]
    \h{\kappa} &= n^{-1/2}\max\{ s_1(\hX), s_p(\hX)\inv \} \le 2\kappa\\[0.25em]
    \h{\Msigma} &= \max_{i < j} \abs{e_{ij}-\bar{e}} \le \Msigma\log{n},
\end{align}
and from \ref{bound-10}, with probability greater than $1 - O(n^{-2})$ and for sufficiently large $n$,
\begin{align}
    \h{\sigma}^2 \ge \sigma^2 - c_{10}\frac{\log{n}}{\sqrt{n}} \ge \frac{\sigma^2}{4}.
\end{align}
Plugging these bounds into $\h{\const}_1, \h{C}'_1$, we obtain
\begin{align}
    \sup_{t \in \R}\abs{ \pr\qty\Big( T'_n \le t) - \pr\qty\Big( T\b_n \le t) } &=  2\const_1\qty\big(2\kappa, 2\Rx, \sigma/2, \varsigma\log{n})\frac{\log^3{n}}{\sqrt{n}} + 2C'_1\qty\big(2\kappa, 2\Rx, \varsigma\log{n}, \sigma/2)\sqrt{\frac{\log{n}}{n}}.
\end{align}
Since $\const_1(\pars) = \const_0(\kappa, \Rx, \sigma)^2 \cdot \varsigma^2$ from \cref{proof:thm:main}, it follows that
\begin{align}
    \sup_{t \in \R}\abs{ \pr\qty\Big( T'_n \le t) - \pr\qty\Big( T\s_n \le t) } &=  2\const_1(\kappa, \Rx, \sigma, \varsigma)\frac{\log^5{n}}{\sqrt{n}},
\end{align}
The proof of \cref{eq:bootstrap-validity-empirical} is then completed by plugging the bound into \cref{eq:ks-hat-tn} with $\pr( T\s_n \le t)$ in place of $\pr(G \le t)$. The proof of the coverage guarantee for $\conf\s_{\alpha, i}$ in \cref{eq:bootstrap-validity-empirical} is identical to the proof for \cref{cor:confidence-multiplier-bootstrap} in \cref{proof:cor:confidence-multiplier-bootstrap}.
\qed



\bibliographystyle{abbrvnat}
\bibliography{./refs}

\begin{table}[t]
    \centering
    \small
    \caption{Summary of Notations}
    \label{tab:notation}
    \begin{tabularx}{\textwidth}{lX}
        \toprule[2pt]\\
        \multicolumn{2}{c}{\textsc{Linear Algebra}}\\[-0.3em]
        \multicolumn{2}{c}{\dotfill}\\
        $1_k, 0_k, e_k \in \R^k$ & vector of all ones, vector of all zeros, $k$th standard basis vector. \\
        $I_k, \O_k, J_k, H_k$ & $k \times k$ identity, all zeros, all ones, and the centering matrix $H_k := I_k - J_k/k$. \\
        $\orth{k}$ & group of orthogonal matrices in $\R^{k \times k}$. \\
        $\norm{x}$ & Euclidean $\ell_2$-norm of $x \in \R^k$. \\
        $\opnorm{A}, \frobenius{A}$ & $\ell_2$-operator norm and the Frobenius norm of $A \in \R^{k_1 \times k_2}$.\\
        $\ttinf{A}$ & $\ell\ttinft$-operator norm of $A \in \R^{k_1 \times k_2}$ where $\ttinf{A} = \max_{\norm{x}=1}\norm{Ax}_\infty$.\\
        $A_{i, *}, A_{*, j}$ & $i$th row and $j$th column of $A \in \R^{k_1 \times k_2}$. \\
        $\lambda_1(B) \ge \dots \ge \lambda_k(B)$ & non-increasing sequence of eigenvalues of $B \in \R^{k \times k}$.\\
        \hline\\[-0.4em]
        
        \multicolumn{2}{c}{\textsc{MDS}}\\[-0.3em]
        \multicolumn{2}{c}{\dotfill}\\
        $X = U\Lambda^{1/2}Q$ & Reduced rank-$p$ singular value decomposition of the \emph{centered} configuration $X$.\\ 
        $\Del, D$ & The noiseless and noisy dissimilarity matrices, respectively, where $D = \Del + \Eps$.\\
        $\Delc, \Dc$ & The double-centered dissimilarities $\Delc = -\half H\Del H$ and $\Dc = -\half HD H$.\\
        $\hU\hL\hU\tr$ & The rank-$p$ spectral decomposition of $\Dc = -\half HDH$.\\
        $\hX \in \Rnp$ & The output of the classical MDS algorithm where $\hX = \hU \hL^{1/2}$.\\
        $\hQ \in \orth{p}$ & The Frobenius-optimal Procrustes alignment of $\hU$ to $U$. See \cref{eq:procrustes}.\\
        $\hg \in \euc{p}, \hP \in \orth{p}$ & The optimal alignment of $X$ to $\hX$. $\hP = \hQ\tr Q$ and $\hg(v) = \hP v$. See \cref{eq:hg}.\\
        \hline\\
        
        \multicolumn{2}{c}{\textsc{Statistical}}\\[-0.3em]
        \multicolumn{2}{c}{\dotfill}\\
        $\mathscr{L}(\xi)$ & The law/probability distribution associated with the random variable $\xi$.\\
        $\norm{\xi}_{L^p}$ & $L^p$ norm of a real valued random variable $\xi$, i.e., $\norm{\xi}_{L^p} = (\E\abs{\xi}^p)^{1/p}$.\\
        $\norm{\xi}_{\psi_1}, \norm{\xi}_{\psi_2}$ & sub-exponential/Gaussian norm of $\xi$; $\norm{\xi}_{\psi_p} = \inf\{k > 0: \E\exp(|\xi/k|^p) \le 2\}$.\\
        & If $\supp(\xi) \subset \R^k$, then $\norm{\xi}_{\psi_p} = \max_{\norm{x}=1}\norm{x\tr\xi}_{\psi_p}$.\\
        $N(\mu, \Sigma)$ & Normal distribution with mean $\mu$ and covariance $\Sigma$.\\
        $\chi^2_k, \chi^2_k(\lambda)$ & \mbox{The central/non-central Chi-squared distribution with non-centrality $\lambda$.}\\[0.1em]
        $\ks$ & The Kolmogorov-Smirnov metric; $\ks(X, Y) = \textstyle{\sup_{t \in \R} |\pr(X \le t) - \pr(Y \le t)|}$\\
        $T_n, \h{T}_n$ & Standardized statistic for the deviation of $\hX$ from $X$ and its plug-in counterpart\\
        $T\b_n, T\s_n$ & Multiplier and empirical bootstrap versions of $T_n$\\
        \hline\\
        
        \multicolumn{2}{c}{\textsc{Asymptotic}}\\[-0.3em]
        \multicolumn{2}{c}{\dotfill}\\
        $a_n = O(b_n), a_n \lesssim b_n$ & there exists a constant $C > 0$ such that $\abs{a_n} \le C \abs{b_n}$ for all $n > N_C$.\\
        $a_n \asymp b_n$ & $a_n = O(b_n)$ and $b_n = O(a_n)$.\\
        $a_n = o(b_n), a_n \ll b_n$ & $\lim_{n \to \infty}\abs{a_n/b_n} = 0$.\\
        $a_n \sim b_n$ & $ \abs{a_n/b_n - 1} = o(1)$.\\
        $\xi_n = \Op(a_n)$ & there exists $C > 0$ such that $\pr(\abs{\xi_n/a_n} > C) \le 1/n$ for all $n > N_C$.\\
        $\xi_n = o_p(1)$ & $\lim_{n \to \infty}\pr(\abs{\xi_n/a_n} > C) = 0$ for all $C > 0$.\\
        
        \bottomrule[2pt]
    \end{tabularx}
\end{table}


\clearpage
\appendix


\section{Toolkit}
\label{sec:toolkit}

The first result establishes a Cramér-type moderate deviation bound for the $\ell_2$-norm of sums of independent random vectors.

{%
\begin{proposition}[Theorem~4.3~of~\citealp{fang2023p}]\label{prop:fang}
    Let $W = (1/\sqrt{n})\sum_{i=1}^n X_i \in \R^p$ where $X_1, \dots, X_n$ are independent random vectors, $\E(X_i) = 0$ for all $i$, and $\Var(W) = I_p$. Suppose $\norm{X_i}_{\psi_1} \le b$ for all $i \in [n]$. Let $Z \sim N(0, I_p)$. For any $p \ge 2$, and for {$S_n := {p^{1/4}b^2}/{\sqrt{n}},$}
    there exists $c > 0$ such that for sufficiently large $n$ 
    \begin{align}
        \abs{\frac{\pr\qty(\norm{W} > t)}{\pr\qty(\norm{Z} > t)} - 1} \lesssim S_n (1 + t) ( p\log{p} + \abs{\log{S_n}} + t^2 ) \qq{for all } 0 \le t \lesssim S_n^{-1/3}.
    \end{align}
\end{proposition}
}

In particular, for all $u_n \asymp \sqrt{\log{n}}$, it follows that
\begin{align}
    \frac{\pr\qty(\norm{W} > u_n)}{\pr\qty(\norm{Z} > u_n)} = 1 + O\qty(\frac{b^2\log^{3/2}{n}}{\sqrt{n}}).\label{eq:moderate-deviation-simplified}
\end{align}

{
The proof of \cref{thm:main} requires a Poisson approximation for sums of dependent Bernoulli random variables. We first define sequences of positively/negatively related Bernoulli random variables.

\begin{definition}[Definition~2.1.1 of \citealp{barbour1992poisson}]
    A collection $\qty{B_\alpha: \alpha \in \mathcal{I}^+}$ of Bernoulli random variables is said to be \textbf{positively related} if: for each $\alpha \in \mathcal{I}^+$ there exists a collection of random variables $\qty{B_{\beta, \alpha}: \beta \in \mathcal{I}^+ \setminus \qty{\alpha}}$ defined on the same probability space such that
    \begin{align}
        \mathscr{L}\qty( (B_{\beta, \alpha})_{\beta \in \mathcal{I}^+ \setminus \qty{\alpha}} ) = \mathscr{L}\qty( (B_\beta)_{\beta \in \mathcal{I}^+ \setminus \qty{\alpha}} | B_\alpha = 1 ) \qq{and} B_{\beta, \alpha} \ge B_\beta \quad\forall\;\beta \in \mathcal{I}^+ \setminus \qty{\alpha}.\label{eq:positive-relation}
    \end{align}
    Similarly, the collection $\qty{B_\alpha: \alpha \in \mathcal{I}^-}$ is said to be \textbf{negatively related} under the same conditions above but with $B_{\beta, \alpha} \le B_\alpha\; \forall \beta \in \mathcal{I}^-\setminus \qty{\alpha}$ in \cref{eq:positive-relation}.
\end{definition}
}
The following result from \cite{barbour1992poisson} provides a Poisson approximation for the sum of dependent Bernoulli random variables in the total-variation metric.

{%
\begin{proposition}[Theorem~2.C of \citealp{barbour1992poisson}]\label{prop:barbour}
    Suppose $W = \sum_{\alpha \in \mathcal{I}} B_\alpha$ where each $B_\alpha \sim \textup{Ber}(p_\alpha)$, and suppose that for each $\alpha \in \mathcal{I}$ there exists a partition $\qty{\mathcal{I}^-_\alpha, \mathcal{I}^+_\alpha, \mathcal{I}^0_\alpha}$ of $\mathcal{I} \setminus \qty{\alpha}$ such that $\qty{B_\alpha: \alpha \in \mathcal{I}^-_\alpha}$ are negatively related, $\qty{B_\alpha: \alpha \in \mathcal{I}^+_\alpha}$ are positively related. Then,
    \begin{align}
        \tv\qty\Big( \mathscr{L}(W), \textup{Poi}(\lambda) ) \le \frac{1 - e^{-\lambda}}{\lambda} 
        \Bigg( 
        &\sum_{\alpha \in \mathcal{I}}p_\alpha^2 + 
        \sum_{\alpha \in \mathcal{I}}\sum_{\beta \in \mathcal{I}^-_\alpha} \abs{ \Cov(B_\alpha, B_\beta) } \\
        &+ \sum_{\alpha \in \mathcal{I}}\sum_{\beta \in \mathcal{I}^+_\alpha} \Cov(B_\alpha, B_\beta) + \sum_{\alpha \in \mathcal{I}}\sum_{\beta \in \mathcal{I}^0_\alpha} {\big[}\E(B_\alpha B_\beta) + p_\alpha p_\beta {\big]}
        \Bigg),
    \end{align}
    where $\mathscr{L}(W)$ is the distribution of the random variable $W$, and $\lambda := \sum_{\alpha \in \mathcal{I}} p_\alpha$.
\end{proposition}
}

We also require the following well-known properties of the Orlicz $\psi_\alpha$-norms from \citet[Lemma~2.2.2]{van2000asymptotic}, \citet[Lemma~2.7.7]{vershynin2018high} and \citet[Eq.~3.5]{kuchibhotla2022moving}.

{%
\begin{proposition}\label{lem:orlicz-norms}
    Suppose $\xi$ is a random variable with $\norm{\xi}_{\psi_\alpha} < \infty$. Then, for any $t > 0$,
    \begin{align}
        \pr( \abs{\xi} > t ) \le 2\exp(-(t/\norm{\xi}_{\psi_\alpha})^\alpha).
    \end{align}
    Moreover, for (possibly dependent) random variables $\xi_1, \dots, \xi_N$,
    \begin{align}
        \norm{\max_{i \in [n]}\xi_i}_{\psi_\alpha} \lesssim \psi_\alpha\inv(N) \cdot \max_{i \in [n]}\norm{\xi_i}_{\psi_\alpha}. 
    \end{align}
    Lastly, for any $\alpha_i, \alpha_j > 0$ and $\xi_{i}, \xi_j$ such that $\norm{\xi_i}_{\psi_{\alpha_i}}, \norm{\xi_{j}}_{\psi_{\alpha_j}} < \infty$,
    \begin{align}
        \norm{\xi_i \cdot \xi_j}_{\psi_{\beta}} \le \norm{\xi_i}_{\psi_{\alpha_i}} \cdot \norm{\xi_j}_{\psi_{\alpha_j}}\qq{where} {1}/{\beta} = {1}/{\alpha_1} + {1}/{\alpha_2}.
    \end{align}
\end{proposition}
}

The next result is a concentration inequality for randomly weighted sums of fixed matrices. The result is a direct application of Theorem~3.2~and~Proposition~A.3 of \citet{kuchibhotla2022moving} along with an $\epsilon$-net argument \citep[Section~4.2.2]{vershynin2018high}. 

{%
\begin{proposition}\label{prop:matrix-bernstein}
    Let $\xi_1, \dots, \xi_n$ be a collection of zero-mean independent random variables with $\max_{i \in [n]}\norm{\xi_i}_{\psi_{\alpha}} \le K < \infty$ for some $\alpha \le 1$, and let $A_1, \dots, A_n \in \R^{p \times q}$ be fixed matrices. 
    Let
    \begin{align}
        \gamma^2 := \max\qty\Big{\Big\|{ \textstyle{\sum\limits_{i \in [n]} \E(\xi_i^2) A_i A_i\tr}}\Big\|_2, \Big\|{ \textstyle{\sum\limits_{i \in [n]} \E(\xi_i^2) A_i\tr A_i}\Big\|_2}} \qq{and} M := \max_{i \in [n]} \opnorm{A_i}.
    \end{align}
    Then, with probability at least $1 - 2e^{-t}$,
    \begin{align}
        \Big\|{\textstyle{\sum\limits_{i \in [n]} \xi_i A_i}}\Big\|_2 \lesssim \gamma\sqrt{t + p + q} + MK \cdot ((t + p + q)\log{n})^{1/\alpha}.
    \end{align}
\end{proposition}
}

\begingroup
\newcommand{\Uc}{\mathcal{U}}
\newcommand{\Vc}{\mathcal{V}}
\begin{proof}[Proof of \cref{prop:matrix-bernstein}]
    Let $S := \sum_i \xi_i A_i$, and define $\Uc, \Vc$ to be $1/4$-nets of the unit spheres $\mathbb{S}^{p-1}$ and $\mathbb{S}^{q-1}$, respectively, such that from Corollary~4.2.13 and Exercise~4.4.3 of \citet{vershynin2018high}, we have $\abs{\Uc} \le e^{p\log{9}}, \abs{\Vc} \le e^{q\log{9}}$, and
    \begin{align}
        \opnorm{S} \le 2\sup_{u \in \Uc, v \in \Vc}u\tr S v = 2 \sup_{u \in \Uc, v \in \Vc} S(u, v),\label{eq:e-net}
    \end{align}
    where $S(u, v) = u\tr S v$. Therefore, it suffices to bound the r.h.s. of \cref{eq:e-net}. For $u \in \Uc, v \in \Vc$, define 
    \begin{align}
        Z_i(u, v) := \xi_i \cdot u\tr A_i v \in \R\qc{} K_i(u, v) := \norm{Z_i(u, v)}_{\psi_\alpha}, \qq{and} \gamma^2(u, v) := \sum_{i \in [n]}\Var(Z_i(u, v)).
    \end{align}
    Since $K_i(u, v) = \abs{u\tr A_i v} \cdot \norm{\xi_i}_{\psi_{\alpha}}$ and $\Var(Z_i(u, v)) = \E(\xi_i^2) (u\tr A_i v)^2$, is straightforward to verify that
    \begin{align}
        \sup_{u \in \Uc, v\in\Vc}\max_{i}K_i(u, v) \le M K \qq{and} \sup_{u \in \Uc, v \in \Vc}\gamma^2(u, v) \le \gamma^2.
    \end{align}
    Using Theorem~3.2 of \citep{kuchibhotla2022moving}, we have that
    \begin{align}
        \norm{S(u, v)}_{\Psi_{\alpha, L_{n}(\alpha; u, v)}} \lesssim \gamma(u, v) \qq{for} L_n(\alpha; u, v) \asymp \frac{(\log{n})^{1/\alpha}}{\gamma(u, v)} \cdot \max_{i \in [n]} K_i(u, v),
    \end{align}
    where $\norm{\xi}_{\Psi_{\alpha, L}}$ denotes the GBO-norm defined in Definition~2.3 of their work. Using Proposition~A.3 of \citet{kuchibhotla2022moving}, we have that for all $t > 0$ and with probability at least $1 - 2e^{-t}$,
    \begin{align}
        \abs{S(u, v)} 
        &\lesssim \gamma(u, v) \qty( \sqrt{t} + L_n(\alpha; u, v) t^{1/\alpha} )\\
        &\lesssim \gamma(u, v) \sqrt{t} + (t\log{n})^{1/\alpha} \cdot \max_{i} K_i(u, v) \;\le\; \gamma \sqrt{t} + MK (t\log{n})^{1/\alpha}.
    \end{align}
    Taking a union bound over all $u \in \Uc, v \in \Vc$, it follows that
    \begin{align}
        \pr\qty( \sup_{u \in \Uc, v \in \Vc}\abs{S(u, v)} \gtrsim  \gamma \sqrt{t} + MK (t\log{n})^{1/\alpha}) \le 2e^{-(t - (p+q)\log{9})}.
    \end{align}
    Setting $t \mapsto t + (p+q)\log{9}$ in the above bound, we get that with probability at least $1 - 2e^{-t}$,
    \begin{align}
        \sup_{u \in \Uc, v \in \Vc}\abs{S(u, v)} \lesssim \gamma \sqrt{t + (p + q)} + MK \qty\Big(\qty\big(t + (p + q))\log{n})^{1/\alpha}.
    \end{align}
    The claim follows by plugging the above bound into \cref{eq:e-net}.
\end{proof}
\endgroup

\begin{remark}
    The sub-Gaussian type concentration for sums of heavy-tailed ($\psi_\alpha$ for $\alpha \le 1$) random variables in \cref{prop:matrix-bernstein} is a somewhat surprising consequence of the tail bounds associated with the Generalized Bernstein-Orlicz (GBO) norm \citep{kuchibhotla2022moving}. Since we assume that the dimensions $p, q$ are fixed throughout this work, we omit their dependence when \cref{prop:matrix-bernstein} is invoked in the proofs. 
\end{remark}



\section{Proofs for Auxiliary Lemmas}
\label{sec:proof-auxiliary}

\subsection{Proof of \cref{lem:slutsky}}
\label{proof:lem:slutsky}

For all $t \in \R$, we have
\begin{align}
    \pr\qty\big(T_n \le t) 
    &\le \pr\qty\big(  T_n \le t, \abs{T_n - S_n} \le C u_n ) + \pr\qty\big(  \abs{T_n - S_n} > C u_n )\\[5pt]
    &\le \pr\qty\big( S_n \le t + Cu_n ) + O(r_n)\\[5pt]
    &\le \pr\qty\big( T \le t + Cu_n ) + O(r_n + s_n)\\[5pt]
    &= \pr\qty\big( T \le t ) + \pr( t < T \le t + Cu_n ) + O(r_n + s_n)\\[5pt]
    &\le \pr\qty\big( T \le t ) + \omega_T(Cu_n) + O(r_n + s_n).
\end{align}
Similarly, starting with $\pr( S_n \le t - Cu_n ) \le \pr( T \le t ) + O(r_n)$, an identical argument gives,
\begin{align}
    \pr\qty(T_n \le t) \ge \pr\qty( T \le t ) - \omega_T(Cu_n/a_n) - O(r_n + s_n).
\end{align}
The claim follows by noting that the two bounds above hold for all $t \in \R$. Moreover, if $T$ admits a p.d.f. $f_T$, then
\begin{align}
    \sup_{\substack{t \in \R\\ 0\le h \le \eps}}\pr(t < T \le t + h) = \sup_{\substack{t \in \R\\ 0\le h \le \eps}}\int_{t}^{t+h} f_T(x) dx \le \sup_{0\le h \le \eps}\norm{f_T}_\infty h \le M\eps.
\end{align}
\qed

\subsection{Proof of \cref{lem:chisq-tail}}
\label{proof:lem:chisq-tail}

Since $\norm{Z}^2 \sim \chi_p^2$, a Chi-squared distribution with $p$ degrees of freedom, for all {$u \ge 0$},
\begin{align}
    \pr( \norm{Z} > u )  = \pr( \chi^2_{p} > u^2) = \frac{\Gamma(p/2, u^2/2)}{\Gamma(p/2)},\label{eq:chisq-tail}
\end{align}
where $\Gamma(z)$ is the gamma function and $\Gamma(\alpha, z)$ is the upper-incomplete gamma function which satisfies $\Gamma(\alpha, z) = z^{\alpha-1} e^{-z} (1 + O(1/z))$ as $z \to \infty$ \citep[\S 8.10.3]{olver2010nist}. 

For ${u_n(t) = a_n t + b_n \asymp \sqrt{\log{n}}}$,
\begin{align}
    \pr\qty( \norm{Z} > u_n(t) ) 
    &= \frac{1}{\Gamma(p/2)} \qty(\frac{u_n(t)^2}{2})^{p/2-1} \exp( -\frac{u_n(t)^2}{2} ) \cdot \qty(1 + O(1/\log{n})).
\end{align}
From the expression for $a_n, b_n$ in \cref{eq:anbn},
\begin{align}
    \frac{u_n(t)^2}{2} = \half\qty( t^2a_n^2 + 2a_nb_nt + b_n^2 ) = \frac{t^2}{2b_n^2} + t + \qty\Big(\log{n} + (p/2-1) \log\log{n} - \log\Gamma(p/2)). \label{eq:un-t-square}
\end{align}
and
\begin{align}
    \exp\qty(-\frac{u_n(t)^2}{2}) = \frac{e^{-t-t^2/2b_n^2}}{n} \cdot \frac{\Gamma(p/2)}{(\log{n})^{p/2-1}}.
\end{align}
Substituting this back into the expression in \cref{eq:chisq-tail}, we get
\begin{align}
    \pr\qty\big( \norm{Z} > u_n(t) ) 
    &= \frac{e^{-t-t^2/2b_n^2}}{n}\qty(\frac{u_n(t)^2}{2\log{n}})^{p/2-1} \cdot \qty(1 + O(1/\log{n})).
\end{align}
From \cref{eq:un-t-square} once again,
\begin{align}
    \qty(\frac{u_n(t)^2}{2\log{n}})^{p/2-1} = 1 + O\qty(\frac{\abs{t} + \log\log{n}}{\log{n}}),
\end{align}
which then gives the desired bound:
\begin{align}
    \pr\qty\big( \norm{Z} > u_n(t) ) = \frac{e^{-t}}{n}\qty(1 + O\qty(\frac{\abs{t} + \log\log{n}}{\log{n}})).\label{eq:unt-tail-bound}
\end{align}

\qed

\subsection{Proof of \cref{lem:tail-bound}}
\label{proof:lem:tail-bound}
Recall that $Y_i = \frac{1}{\sqrt{n}} \sum_{k \in [n]} \eps_{ik} \theta_{ik}$ from~\cref{prop:decomposition} where $\max_{i, k}\norm{\theta_{ik}} \le \const_0 \equiv \const_0(\kappa, \Rx, \msigma)$ from~\cref{eq:theta-norm}, and from \ref{noise-1} we also have
\begin{align}
    \max_{i \in [n]}\Norm\big{ \norm{\eps_{ik}\theta_{ik}} }_{\Psi_1} \le \const_0 \Msigma.
\end{align}
If $t \le -b_n^2$, note that $u_n(t) = t/b_n + b_n \le 0$, and, trivially, $\pr( \norm{Y_i} > u_n(t) ) = 1$. On the other hand, for any $\tau > 0$ and for all $-b_n^2 < t < \tau\log{n}$, we have $u_n(t) \asymp \sqrt{\log{n}}$; using \citet[Theorem~4.3]{fang2023p} (see \cref{prop:fang}), for $Z \sim N(0, I_p)$ and $\const_1 \equiv \const_1(\pars) := (\const_0\Msigma)^2$ we have
\begin{align}
    \pr\qty\big(\norm{Y_i} > u_n(t)) = \pr\qty\big(\norm{Z} > u_n(t)) \qty( 1 + O\qty( \frac{\const_1\log^{3/2}{n}}{\sqrt{n}} ) ).
\end{align}
For the last claim, we apply the matrix concentration inequality in \cref{prop:matrix-bernstein}. To this end, for $\alpha=1$ and $q=1$ and for each $i \in [n]$, let
\begin{align}
    K_i := \max_{k \in [n]}\norm{\eps_{ik}}_{\psi_1} \le \Msigma\qc{} 
    M_i := \max_{k \in [n]}\Norm{\frac{1}{\sqrt{n}}\theta_{ik}} \le \frac{\const_0}{\sqrt{n}},\qq{and}
    \gamma_i^2 := \sum_{k \in [n]} \E(\eps_{ik}^2) \frac{\norm{\theta_{ik}}^2}{n} \le \Msigma^2 \const_0^2.
\end{align}
From \cref{prop:matrix-bernstein}, for all $t > 0$ there exists a constant $\tilde{C}_p > 0$ depending on $p$ such that with probability at least $1 - 2e^{-t}$,
\begin{align}
    \norm{Y_i} \;\;\le\;\; \tilde{C}_p\qty(\gamma_i \sqrt{t} + M_i K_i ( t \log{n} )) \;\;\le\;\; \tilde{C}_p \Msigma \const_0\qty( \sqrt{t} + \frac{ t \log{n} }{\sqrt{n}}).\label{eq:bernstein-type-bound}
\end{align}
Define $\tau := 4\tilde{C}_p \Msigma \const_0$, and note that for all $t > \tau\log{n}$, we have
\begin{align}
    u_n(t) = \frac{t}{b_n} + b_n \ge \frac{t}{\sqrt{2\log{n}}} > 2 \tilde{C}_p \Msigma \const_0 \sqrt{2\log{n}} \ge \tilde{C}_p \Msigma \const_0\qty( \sqrt{2\log{n}} + \frac{ 2\log{n} }{\sqrt{n}}).
\end{align}
Therefore, from \cref{eq:bernstein-type-bound}, for all $t > \tau\log{n}$ we have
\begin{align}
    \pr\qty\Big( \norm{Y_i} > u_n(t) ) \le \pr\qty(\norm{Y_i} > \tilde{C}_p \Msigma \const_0\qty( \sqrt{2\log{n}} + \frac{ 2\log{n} }{\sqrt{n}}) ) \le 2e^{-2\log{n}} = O(1/n^2).
\end{align}
\qed

\subsection{Proof of \cref{lem:covariance-bound}}
\label{proof:lem:covariance-bound}

Since $t$ is fixed throughout, in the interest of clarity, we simply write $B_i \equiv B_i(t)$, $u_n \equiv u_n(t)$, etc., throughout. The proof of \cref{lem:covariance-bound} relies on the following local anti-concentration  inequality for the Chi-squared distribution.

\cbox{black!5}
{%
\begin{lemma}\label{lem:chisq-anticoncentration}
    Let $Z \sim N(0, I_p)$. Then, for any $x, \epsilon > 0$,
    \begin{align}\label{eq:chisq-anticoncentration}
        \pr\qty\Big( x < \norm{Z}^2 \le x+\epsilon) \le \frac{\epsilon}{2} \cdot {\pr\qty\Big( \norm{Z}^2 > x )}.
    \end{align}
    Moreover, for all $x, \epsilon > 0$ such that $x-\epsilon > p - 1$,
    \begin{align}\label{eq:chisq-comparison}
        \frac{\pr\qty( \norm{Z}^2 > x - \epsilon )}{\pr\qty( \norm{Z}^2 > x )} \le \frac{e^{\epsilon/2}}{1 - \frac{p-1}{x-\epsilon}}.
    \end{align}
\end{lemma}
}

\cref{lem:chisq-anticoncentration} is proved at the end of this section in \cref{proof:lem:chisq-anticoncentration}. We now proceed with the proof of \cref{lem:covariance-bound}. For each $i \neq j$, let $Y_{i \setminus j}$ be the sum of the individual terms in $Y_i$ excluding the $\eps_{ij} = \eps_{ji}$ term, i.e.,
\begin{align}
    Y_{i \setminus j} := \frac{1}{\sqrt{n}} \sum_{k \neq j} \eps_{ik} \theta_{ik}.\label{eq:yji}
\end{align}

\subsubsection*{Step 1. Simplifying $\Cov(B_i, B_j)$.}

Fix $i \neq j$, let $\xi := \eps_{ij}$ with $\Var(\xi) =\sigma_{ij}^2$, and let $v_1 \equiv v_{1, n} := \frac{1}{\sqrt{n}}\theta_{ij}, v_2 := \frac{1}{\sqrt{n}}\theta_{ji} \in \R^p$. From \cref{eq:yi},
\begin{align}
    Y_i = Y_{i \setminus j} + \xi v_1 \qq{and} Y_j = Y_{j \setminus i} + \xi v_2,
\end{align}

{
Since $\E(Y_i) = \E(\xi)=0$ and the $Y_i$'s are standardized, i.e., $\Var(Y_i) = I_p$ for all $i$, we also have that $\E(Y_{i \setminus j}) = 0$ and $\Sigma_{ij} := \Var(Y_{i \setminus j}) = I_p - \sigma_{ij}^2v_1v_1\tr.$
}
Using the fact that $\norm{\theta_{ij}} \le \const_0$ from \cref{eq:theta-norm}, we have 
\begin{align}\label{eq:sigma-ij-bound}
    \Norm\big{\norm{\eps_{ij}\theta_{ij}}}_{\psi_1} \le {\const_0\Msigma}
    \qq{and} 
    \opnorm{\Sigma_{ij} - I_p} \le \sigma_{ij}^2\norm{v}^2 \le \frac{\Msigma^2\const_0^2}{n} = \frac{\const_1}{n} =: s_n.
\end{align}

Lastly, from \ref{noise-1} we have that $Y_{i \setminus j} \indep Y_{j \setminus i} \indep \xi$. By definition of $B_i$'s in \cref{eq:B-pi-W}, $B_i, B_j$ are conditionally independent given $\xi$, i.e., $(B_i \indep B_j \mid \xi)$. Therefore, using the law of total covariance,
\begin{align}
    \Cov(B_i, B_j) = \Cov\qty\Big( \E(B_i \mid \xi), \E(B_j \mid \xi) ).\label{eq:law-of-total-covariance}
\end{align}
Define $g_1(\xi) := \E(B_i \mid \xi)$ and $g_2(\xi) := \E(B_j \mid \xi)$, i.e.,
\begin{align}
    g_1(\xi) := \pr\qty\Big( \norm{Y_{i \setminus j} + \xi v_1} > u_n \mid \xi ) \qq{and} g_2(\xi) := \pr\qty\Big( \norm{Y_{j \setminus i} + \xi v_2} > u_n \mid \xi ).
\end{align}
Let $\xi'$ be an \iid{} copy of $\xi$. Then, by Jensen's inequality,
\begin{align}
    \abs{\Cov(B_i, B_j)} 
    &\le \E\qty\Big[ \abs{g_1(\xi) - \E' g_1(\xi')} \cdot \abs{g_2(\xi) - \E' g_2(\xi')} ]\\
    &\le \E\E' \qty\Big[ \abs{g_1(\xi) - g_1(\xi')} \cdot \abs{g_2(\xi) - g_2(\xi')} ].
    \label{eq:covariance-jensen}
\end{align}
Moreover, from \ref{noise-2}, since $\xi, \xi'$ are $\Msigma$-sub-exponential there exists {$\tau_n = (1+o(1)) \cdot 4\Msigma\log{n}$}
such that for the event $A := \qty{\abs{\xi} \le \tau_n, \abs{\xi'} \le \tau_n}$ we have
\begin{align}
    \pr\qty(A^c) = \pr\qty( \qty{\abs{\xi} > \tau_n} \cup \qty{\abs{\xi'} > \tau_n} ) = \pr(\abs{\xi} > \tau_n) + \pr(\abs{\xi'} > \tau_n) = O(n^{-4}).
\end{align}
We can bound the right-hand side of \cref{eq:covariance-jensen} as follows:
\begin{align}
    \abs{\Cov(B_i, B_j)} &\le \E\E' \qty\Big[  \abs{g_1(\xi) - g_1(\xi')} \cdot \abs{g_2(\xi) - g_2(\xi')} \cdot \mathbb{1}_A ]\\
    &\qquad\qquad+\E\E' \qty\Big[  \abs{g_1(\xi) - g_1(\xi')} \cdot \abs{g_2(\xi) - g_2(\xi')} \cdot \mathbb{1}_{A^c} ].
\end{align}
Note that $\abs{g_1(\xi) - g_1(\xi')} \le 1$ and $\abs{g_2(\xi) - g_2(\xi')} \le 1$ for all $\xi, \xi'$, since $g_1, g_2$ are probabilities. It follows that
\begin{align}
    \abs{\Cov(B_i, B_j)} &\le \E\E' \qty\Big[  \abs{g_1(\xi) - g_1(\xi')} \cdot \abs{g_2(\xi) - g_2(\xi')} \cdot \mathbb{1}_A ] + \pr(A^c)\\
    &\le \E\E' \qty[  \abs{g_1(\xi) - g_1(\xi')} \cdot \abs{g_2(\xi) - g_2(\xi')} \cdot \mathbb{1}_A ] + O(n^{-4}).\label{eq:covariance-bound-1}
\end{align}
Thus, it suffices to bound $\abs{g_1(\xi) - g_1(\xi')}$ and $\abs{g_2(\xi) - g_2(\xi')}$ on the event $A$.

\subsubsection*{Step 2. Bounding $\abs{g_1(\xi) - g_1(\xi')}$ and $\abs{g_2(\xi) - g_2(\xi')}$ on $A$.}

Let $v \equiv v_{n} = v_{1, n}$. On the event $A$, define
\begin{align}
    S := \Sigma_{ij}^{-1/2}Y_{i \setminus j}, \quad r_n := \norm{\xi v} = \abs{\xi}\cdot\norm{v}, \quad r'_n := \norm{\xi' v} = \abs{\xi'}\cdot\norm{v}. 
\end{align}
Note that $\norm{v}\le \const_0/\sqrt{n}$ and $r_n, r_n' \le \tau_n \norm{v} \lesssim \const_0\Msigma \log{n}/\sqrt{n}$. By an application of the triangle inequality, we have
\begin{align}
    g_1(\xi)- g_1(\xi') &=  \pr\qty( \norm{Y_{i \setminus j} + \xi v} > u_n ) - \pr\qty( \norm{Y_{i \setminus j} + \xi' v} > u_n )\\  
    &\le \pr\qty( \norm{Y_{i \setminus j}} > u_n - r_n ) - \pr\qty(\norm{Y_{i \setminus j}} > u_n + r_n' ).
\end{align}
Similarly, we have $g_1(\xi')- g_1(\xi) \le \pr\qty( \norm{Y_{i \setminus j}} > u_n - r'_n ) - \pr\qty(\norm{Y_{i \setminus j}} > u_n + r_n )$. This implies that
\begin{align}
    \abs{g_1(\xi)- g_1(\xi') } &\le \pr\qty\Big(\norm{Y_{i \setminus j}} > u_n - r_n - r'_n ) - \pr\qty\Big(\norm{Y_{i \setminus j}} > u_n + r_n + r'_n ).\label{eq:anti-concentration-1}
\end{align}
Writing $Y_{i \setminus j} = \Sigma_{ij}^{1/2} S$, for all $z > 0$ we have
\begin{align}
    \pr\qty( \norm{S} \cdot \sqrt{\lambda_{\min}(\Sigma_{ij})} > z ) \le \pr\qty\Big( \norm{\Sigma_{ij}^{1/2}S} > z ) 
    \le \pr\qty\Big( \norm{S} \cdot \sqrt{\lambda_{\max}(\Sigma_{ij})} > z ).
\end{align}
Plugging this into \cref{eq:anti-concentration-1} and by noting that $\lambda_{\max}(\Sigma_{ij}) \le 1+s_n$ and $\lambda_{\min}(\Sigma_{ij}) \ge 1-s_n$ from \cref{eq:sigma-ij-bound},
\begin{align}
    \abs{ g_1(\xi)- g_1(\xi') } 
    &\le \pr\qty\Big( \norm{S} > \frac{u_n - r_n - r'_n}{\sqrt{1+s_n}} ) - \pr\qty\Big( \norm{S} > \frac{u_n + r_n + r'_n}{{\sqrt{1-s_n}}} ).
\end{align}
Note that since 
\begin{align}
    u_n \asymp \sqrt{\log{n}}\qc{} r_n, r_n' \lesssim \const_0\Msigma\frac{\log{n}}{\sqrt{n}},\qq{and}s_n \asymp \frac{1}{n},
\end{align}
there exists 
$$
\alpha_n = (1 + o(1)) \cdot \const_0\Msigma\log{n}/\sqrt{n}
$$ 
such that $(1+s_n)^{-1/2}(u_n - r_n - r'_n) \ge u_n - \alpha_n$ and $(1-s_n)^{-1/2}(u_n + r_n + r'_n) \le u_n + \alpha_n$, and, therefore,
\begin{align}
    \abs{ g_1(\xi)- g_1(\xi') } \le \pr\qty\Big( \norm{S} > {u_n - \alpha_n} ) - \pr\qty\Big( \norm{S} > {u_n + \alpha_n} ),
\end{align}
Observe that $S= \Sigma_{ij}^{-1/2}Y_{i \setminus j}$ is the normalized sum of independent random variables such that $\Var(S) = I_p$. Therefore, for $Z \sim N(0, I_p)$ we can now apply the Cramér-type moderate deviation bound in \cref{prop:fang} with $u_n \asymp \sqrt{\log{n}}$ to get

{
\begin{align}
    \abs{ g_1(\xi)- g_1(\xi') } \le \pr\qty\Big( \norm{Z} > u_n - \alpha_n )(1 + \eta_n) - \pr\qty\Big( \norm{Z} > u_n + \alpha_n )(1+\eta_n),\label{eq:anti-concentration-2}
\end{align}
}
where $\abs{\eta_n} \lesssim \const_0\Msigma\log^{3/2}{n}/\sqrt{n}$. 
{
Since $1 + \abs{\eta_n} \le 2$ for sufficiently large $n$, it follows that
\begin{align}
    \abs{ g_1(\xi)- g_1(\xi') } \le (1+\eta_n) \pr\qty\Big( u_n - \alpha_n < \norm{Z} \le u_n + \alpha_n ) \le 2\pr\qty\Big( u_n - \alpha_n < \norm{Z} \le u_n + \alpha_n ). \label{eq:anti-concentration-3}
\end{align}
}
Applying \cref{lem:chisq-anticoncentration} now gives
\begin{align}
    \pr\qty\Big( (u_n - \alpha_n)^2 < \norm{Z}^2 \le (u_n + \alpha_n)^2 ) 
    &\le \qty(\frac{(u_n + \alpha_n)^2 - (u_n - \alpha_n)^2}{2}) \cdot \pr\qty\big( \norm{Z}^2 > (u_n - \alpha_n)^2 )\\
    &= 2{\alpha_nu_n} \cdot \pr\qty\big( \norm{Z}^2 > u_n^2 - \epsilon_n ),\label{eq:chisq-anticoncentration-4}
\end{align}
where $\epsilon_n = 2u_n\alpha_n - \alpha_n^2 \asymp \log^{3/2}{n}/\sqrt{n}$. Since {$u_n - \alpha_n \gg 1$}, \cref{eq:chisq-comparison} applies and it follows that
\begin{align}
    \pr\qty\big( \norm{Z}^2 > u_n^2 - \epsilon_n ) 
    \lesssim \pr(\norm{Z}^2  > u_n^2).\label{eq:comparison-bound-1}
\end{align}
Plugging \cref{eq:comparison-bound-1} and \cref{eq:chisq-anticoncentration-4} back into \cref{eq:anti-concentration-3} and by noting that $u_n\alpha_n \lesssim \const_0\Msigma\frac{\log^{3/2}{n}}{\sqrt{n}}$,
\begin{align}
    \abs{ g_1(\xi)- g_1(\xi') } 
    &\le \const_0\Msigma\frac{\log^{3/2}{n}}{\sqrt{n}} \cdot \pr\qty\big( \norm{Z} > u_n ).
\end{align}
uniformly {on the event $A$}.  The same bound also holds for $\abs{g_2(\xi) - g_2(\xi')}$. Substituting this back into \cref{eq:covariance-bound-1} and by noting that $\const_1 = \const_0^2\Msigma^2$, we have
\begin{align}
    \abs{\Cov(B_i, B_j)} &\le \const_1\frac{\log^{3}{n}}{{n}} \cdot \pr\qty\big( \norm{Z} > u_n )^2 + O(n^{-4}),
\end{align}
which gives the desired bound in \cref{eq:covariance-bound}.
\qed

\subsubsection{Proof of \cref{lem:chisq-anticoncentration}}
\label{proof:lem:chisq-anticoncentration}
For $Z \sim N(0, I_p)$, let $f_{p}(x)$ denote the p.d.f. of $\norm{Z}^2 \sim \chi^2_p$ and $\bar{F}_p(x) = \pr\qty(\norm{Z}^2 > x)$ where
\begin{align}
    f_p(x) = \frac{(x/2)^{p/2-1}e^{-x/2}}{2\Gamma(p/2)} \qq{and} \bar{F}_p(x) = \frac{\Gamma(p/2, x/2)}{\Gamma(p/2)},
\end{align}
where $\Gamma(s, z)$ is the upper incomplete gamma function. Using the lower bound $\Gamma(s,z) \ge z^{s-1}e^{-z}$ \citep[\S8.10.2]{olver2010nist}, we have
\begin{align}
    \frac{f_p(x)}{\bar{F}_p(x)} = \half \cdot \frac{(x/2)^{p/2-1}e^{-x/2}}{\Gamma(p/2, x/2)} \le \half.\label{eq:mills-ratio}
\end{align}
The inverse of the ratio in \cref{eq:mills-ratio} is sometimes also called the Mills' ratio. Using the integral mean-value theorem for $f_{p}(x)$, there exists some $x' \in [x, x+\epsilon]$ such that
\begin{align}
    \pr\qty( x < \norm{Z}^2 \le x+\epsilon) &= \int_{x}^{x+\epsilon} f_p(t) dt = \epsilon \cdot f_p(x').
\end{align}
From \cref{eq:mills-ratio} and using the fact that $\bar{F}_p(x') = \pr(\norm{Z}^2 > x') \le \pr(\norm{Z}^2 > x)$, we have
\begin{align}
    \pr\qty( x < \norm{Z}^2 \le x+\epsilon) &= \epsilon \cdot f_p(x')  \le \frac{\epsilon}{2} \cdot \bar{F}_p(x') \le \frac{\epsilon}{2} \cdot \bar{F}_p(x). \FINEQ
\end{align}

For the claim in \cref{eq:chisq-comparison}, using the upper bound $\Gamma(s, z) \le z^{s-1}e^{-z} / (1 - (p-1)/z)$ in the numerator \citep[\S8.10.3]{olver2010nist}:
\begin{align}
    \frac{\pr\qty( \chi^2_p > x - \epsilon )}{\pr\qty( \chi^2_p > x )} = \frac{\Gamma(p/2, (x-\epsilon)/2)}{\Gamma(p/2, x/2)} \le \frac{(x-\epsilon)^{p/2-1}e^{-x/2 + \epsilon/2}}{x^{p/2-1}e^{-x/2} \cdot (1 - \frac{p-1}{x-\epsilon})} \le \qty(1-\frac{\epsilon}{x})^{p/2-1} \cdot \frac{e^{\epsilon/2}}{1 - \frac{p-1}{x-\epsilon}}.
\end{align}
The final bound follows by noting that $1 - \epsilon/x \le 1$.
\qed

\subsection{Proof of \cref{lem:poisson-approximation}}
\label{proof:lem:poisson-approximation}

In the interest of avoiding notational clutter, for fixed $t \in \R$, let $B_i \equiv B_i(t)$, $\lambda_n \equiv \lambda_n(t)$ and $u_n \equiv u_n(t)$. Let $\mathscr{L}(W)$ denote the distribution of $W$ and let $P_{\lambda_n} \sim \text{Poi}(\lambda_n)$ be a Poisson random variable with parameter $\lambda_n$. Since the event $\qty{M_n \le u_n} = \qty{W = 0}$, we have
\begin{align}
    \abs{ \pr\qty({M_n} \le u_n) - e^{-\lambda_n} } 
    &= \abs\Big{ \pr(W = 0) - \pr( P_{\lambda_n} = 0 ) }
    \le {\tv\qty\Big( \mathscr{L}(W), \text{Poi}(\lambda_n) )},\label{eq:poisson-approximation-1}
\end{align}
where $\tv(\cdot, \cdot)$ is the total variation metric. We aim to apply \cref{prop:barbour} to bound \cref{eq:poisson-approximation-1}.

To this end, for each $j \neq i$, let $Y_{j \setminus i}$ be the sum of the individual terms in $Y_j$ excluding the $\eps_{ij}=\eps_{ji}$ term as given in \cref{eq:yji}. Consider the random variable $B'_{ji}$ given as follows. If $B_i = 1$, then $B'_{ji} = B_j$; otherwise, if $B_i = 0$, then {draw}  only the $i$th row and the $i$th column of $\Eps$, i.e., $\qty{ \eps_{ik}=\eps_{ki}: k \in [n]\setminus \{i\} }$, until $B_i = 1$. Let $\Eps' = (\eps'_{ij})$ be the resulting matrix, and let $Y_j'$ be given by \cref{eq:yi} with $\Eps'$ in place of $\Eps$. It follows that
\begin{align}
    \qty{ B_{ji}: j \in [n] } \eqd \qty{ B_j: j \in [n]  \mid B_i = 1},
\end{align}
In other words, the $B_{ji}$'s follow the distribution of $B_j$ conditional on $B_i=1$. Moreover, note that $\eps'_{jk} = \eps_{jk}$ for all $j, k \neq i$. Therefore, $Y'_{j} := Y_{j \setminus i} + \frac{1}{\sqrt{n}}\eps_{ji}'\theta_{ji}$. Now, we define the sets $\mathcal{I}^+_i, \mathcal{I}^-_i \subset [n] \setminus \qty{i}$ as follows:
\begin{align}
    \mathcal{I}^+_i = \qty\Big{ j \in [n] \setminus \qty{i}: \eps'_{ij} > -\eps_{ij} + \sqrt{n}\tfrac{\theta_{ji}\tr Y_{j\setminus i}}{\norm{\theta_{ji}}^2}  } 
    \qq{and} 
    \mathcal{I}^-_i = \qty\Big{ j \in [n] \setminus \qty{i}: \eps'_{ij} \le -\eps_{ij} + \sqrt{n}\tfrac{\theta_{ji}\tr Y_{j\setminus i}}{\norm{\theta_{ji}}^2}  }.
\end{align}
Equivalently, note that:
\begin{align}
    \norm{Y'_j}^2 > \norm{Y_j}^2 &\iff \tfrac{\norm{\theta_{ji}^2}}{n}(\eps'_{ij})^2 + \tfrac{2 \theta_{ji}\tr Y_{j \setminus i}}{\sqrt{n}}\eps'_{ij} > \tfrac{\norm{\theta_{ji}^2}}{n}(\eps_{ij})^2 + \tfrac{2\theta_{ji}\tr Y_{j \setminus i}}{\sqrt{n}} \eps_{ij}\\
    &\iff \eps'_{ij} > - \eps_{ij} + \sqrt{n} \tfrac{2\theta_{ji}\tr Y_{j \setminus i}}{\norm{\theta_{ji}}^2};
\end{align}
a similar argument also holds for $j \in \mathcal{I}^-_i$, and we can equivalently write:
\begin{align}
    \mathcal{I}^+_i := \qty\big{ j \in [n] \setminus \qty{i}: \norm{Y'_j} > \norm{Y_j}  }, \text{ and } 
    \mathcal{I}^-_i := \qty\big{ j \in [n] \setminus \qty{i}: \norm{Y'_j} \le \norm{Y_j}  }.
\end{align}
In other words, $\mathcal{I}^+_i$ is the set of indices $j$ such that $\norm{Y'_j} > \norm{Y_j}$ (mutatis mutandis for $\mathcal{I}^-_i$), and are constructed purely based on the values of the resampled $\eps'_{ij}$s. It follows that
\begin{align}
    \begin{cases}
        B_{ji} > B_j & \text{if } j \in \mathcal{I}^+_i,\\
        B_{ji} \le B_j & \text{if } j \in \mathcal{I}^-_i;
    \end{cases}
\end{align}
therefore, from \citet[Definition~2.1.1]{barbour1992poisson}, for each $i$, $\qty{B_{ji}: j \in \mathcal{I}_i^{\pm}}$ is a monotone coupling and $\qty{\mathcal{I}^+_i, \mathcal{I}^-_i}$ is a partition of $[n]\setminus \qty{i}$ into positively and negatively related random variables. We can now apply \citet[Theorem~2.C]{barbour1992poisson} (see \cref{prop:barbour}) with $\mathcal{I}^0_i = \emptyset$ to get
\begin{align}
    \tv\qty\Big( \mathscr{L}(W), \ \text{Poi}(\lambda_n) ) 
    &\le \frac{1 - e^{-\lambda_n}}{\lambda_n} 
    \qty( 
    \sum_{i \in [n]}\pi_i^2 + 
    \sum_{i \in [n]}\sum_{j \in \mathcal{I}^-_i} \abs{ \Cov(B_i, B_j) } + \sum_{i \in [n]}\sum_{j \in \mathcal{I}^+_i} \Cov(B_i, B_j)
    )\\
    &\le \frac{1 - e^{-\lambda_n}}{\lambda_n}\qty(\sum_i\pi_i^2 + 
    \sum_{i \in [n]}\sum_{j \in [n]\setminus \qty{i}} \abs{ \Cov(B_i, B_j) }).\label{eq:tv-bound}
\end{align}
The final result now follows from \cref{eq:poisson-approximation-1}.
\qed

\subsection{Proof of \cref{lem:omega-error}}
\label{proof:lem:omega-error}

Note that since $X$ is assumed to be centered in \ref{assumption:compact}, we have $HX = X$ and the expression for $\Om_i$ and $\hOm_i$ are simplified to
\begin{align}
    \Om_i \equiv \Om_i(X) = \frac{n}{4} (X\tr X)\inv (X\tr \Sigma_i X) (X \tr X)\inv \qq{and} \hOm_i = \frac{n}{4} (\hX\tr \hX)\inv (\hX\tr \h{\Sigma}_i \hX) (\hX\tr \hX)\inv,\label{eq:Omega-i-hat}
\end{align}
and from \cref{eq:Omega-equivariance}, we also have
\begin{align}
    \hP \Om_i \hP\tr = \Om_i( X\hP\tr ) = \Om_i( \hg(X) ).
\end{align}

Let $A := (\hg(X)\tr \hg(X))\inv$, $\h{A} := (\hX\tr\hX)\inv$, $B_i := (\hg(X)\tr \Sigma_i \hg(X))$ and $\h{B}_i := (\hX\tr \h{\Sigma}_i \hX)$.
Then, we can write
\begin{align}
    \hOm_i - \hP\Om_i\hP\tr = \frac{n}{4}( \h{A}\h{B}_i\h{A} - AB_i A ) &= \frac{n}{4} \qty\Big( (\h{A} - A) \h{B}_i \h{A} + A (\h{B}_i - B_i) \h{A} + A B_i (\h{A} - A) ).\label{eq:omega-diff}
\end{align}
Since $\hg(X) = X\hP\tr$ is a rigid transformation, from \ref{assumption:compact} and \ref{noise-1} we have
\begin{align}
    \opnorm{A} = \opnorm{(X\tr X)\inv} \le \frac{\kappa}{n} \qq{and}
    \max_{i \in [n]}\opnorm{B_i} = \max_{i \in [n]}\opnorm{X\tr\Sigma_i X} \le {4\Msigma^2 \cdot \kappa^2 n}.\label{eq:omega-bounds-1}
\end{align}
From \ref{bound-4} and \ref{bound-7}, with probability greater than $1 - O(n^{-2})$ we also have
\begin{align}
    \opnorm{\h{A} - A} &= \Opnorm\big{(\hX\tr \hX)\inv - (\hg(X)\tr \hg(X))\inv} \le c_4 \cdot n^{-3/2}\\[0.5em]
    \max_{i \in [n]}\opnorm{\h{B}_i - B_i} &= \max_{i \in [n]}\Opnorm\big{(\hX\tr \h{\Sigma}_i \hX) - (\hg(X)\tr \Sigma_i \hg(X))} \le c_8 \cdot \log^2{n}\sqrt{n}.\label{eq:omega-bounds-2}
\end{align}
Plugging \cref{eq:omega-bounds-1,eq:omega-bounds-2} into \cref{eq:omega-diff} and using $\max_i\opnorm{\h{B}_i} \le \max_{i}\opnorm{B_i} + \max_{i}\opnorm{\h{B}_i - B_i}$, we get that with probability greater than $1 - O(n^{-2})$,
\begin{align}
    \max_{i \in [n]}\opnorm{\hOm_i - \hP\Om_i\hP\tr} 
    \;\;\lesssim\;\; \frac{n}{4} \qty( c_4 n^{-3/2} \cdot \frac{\kappa}{n} \cdot 4\Msigma^2 \kappa^2 n + \frac{\kappa^2}{n^2} \cdot c_8 \log^2{n}\sqrt{n}  )
    \;\;\lesssim\;\; \kappa^2c_8 \cdot \frac{\log^2{n}}{\sqrt{n}}.\label{eq:omega-diff-bound}
\end{align}
Finally, writing
\begin{align}
    {\Om_i^{-1/2} \; \hP\tr\hat\Om_i\hP\; \Om_i^{-1/2} - I_p} = \Om_i^{-1/2} \; \hP\tr \qty\Big(\h{\Om}_i - \hP \Om_i \hP\tr) \hP\; \Om_i^{-1/2},
\end{align}
we get
\begin{align}
    \max_{i \in [n]}\Opnorm\Big{\Om_i^{-1/2} \; \hP\tr\hat\Om_i\hP\; \Om_i^{-1/2} - I_p} \le \max_{i \in [n]}\Opnorm\Big{\Om_i^{-1/2}}^2 \cdot \max_{i \in [n]}\Opnorm\Big{\h{\Om}_i - \hP \Om_i \hP\tr}.\label{eq:omega-id-bound}
\end{align}
From \cref{eq:omega-diff-bound} and using the bound for $\Omega_i^{-1/2}$ from \cref{eq:theta-norm}, with probability at least $1 - O(n^{-2})$:
\begin{align}
    \max_{i \in [n]}\Opnorm\Big{\Om_i^{-1/2} \; \hP\tr\hat\Om_i\hP\; \Om_i^{-1/2} - I_p} \lesssim \frac{\kappa^2}{\msigma^2} \cdot \kappa^2c_8 \cdot  \frac{\log^2{n}}{\sqrt{n}}.\label{eq:omega-error-bound}
\end{align}
This completes the proof of \cref{lem:omega-error} by taking $C_2'(\pars) = \kappa^4 c_8 / \msigma^2$.
\qed

\subsection{Proof of \cref{lem:psi-bound-iid}}
\label{proof:lem:psi-bound-iid}

For $\Psi = \Om^{-1/2} \hP\tr \hOm \hP \Om^{-1/2}$, exactly as in the proof of \cref{lem:omega-error} above, from \cref{eq:omega-id-bound} we have
\begin{align}
    \opnorm{\Psi - I_p} &\le \opnorm{\Om^{-1/2}}^2 \cdot \opnorm{\h{\Om} - \hP \Om \hP\tr},\label{eq:psi-diff-bound}
\end{align}
where, from \cref{eq:omega-iid,eq:psi-iid},
\begin{align}
    \hOm = \frac{n}{4} \cdot \hat\sigma^2 \qty\big({\hX\tr\hX})\inv
    \qq{and}
    \hP\Om\hP\tr = \frac{n}{4} \cdot \sigma^2  \hP \qty( {X\tr X} )\inv \hP\tr = \frac{n}{4} \cdot \sigma^2 \qty({\hg(X)\tr\hg(X)})\inv.
\end{align}
Therefore, writing $A = (\hg(X)\tr \hg(X))\inv$ and $\h{A} = (\hX\tr \hX)\inv$, and using \ref{bound-4} and \ref{bound-10} it follows that with probability greater than $1 - O(n^{-2})$,
\begin{align}
    \opnorm{\hOm - \hP\Om\hP\tr} 
    &\le \frac{n}{4} \qty\Big( \opnorm{ (\h{A}-A)\hat\sigma^2 } + \opnorm{ A(\hat\sigma^2 - \sigma^2) } )\\
    &\lesssim \frac{n}{4}\qty( \hat\sigma^2 \cdot c_4 n^{-3/2} \;\;+\;\; \frac{\kappa}{n} \cdot c_{10} \tfrac{\log{n}}{n} )\\
    &\lesssim \frac{n}{4}\qty( \qty{ \sigma^2 + c_{10}\tfrac{\log{n}}{n} } \cdot c_4 n^{-3/2} \;\;+\;\; \frac{\kappa}{n} \cdot c_{10} \tfrac{\log{n}}{n} ) \;\;\lesssim\;\; c_4 \frac{\sigma^2}{\sqrt{n}}.\label{eq:omega-diff-bound-iid}
\end{align}
Using the fact that $\opnorm{(X\tr X/n)^{1/2}} \le \kappa$ it follows that $\opnorm{\Om^{-1/2}}^2 \le \kappa^2/4\sigma^2$; \cref{lem:psi-bound-iid} now gives:
\begin{align}
    \opnorm{\Psi - I_p} &\lesssim \frac{\kappa^2}{4\sigma^2} \cdot c_4 \frac{\sigma^2}{\sqrt{n}} \lesssim  \frac{\kappa^2 c_4}{\sqrt{n}}.\label{eq:psi-bound-final}
\end{align}
Let $z_n = {\kappa^2 c_4}/{\sqrt{n}}$ be the r.h.s. of \cref{eq:psi-bound-final}. On the event that the bound above holds, similar to \cref{eq:plug-in-bounds-3},
\begin{align}
    (1-z_n)I_p \preccurlyeq \Psi \preccurlyeq (1+z_n)I_p \implies (1+z_n)^{-1/2} I_p \preccurlyeq \Psi^{-1/2} \preccurlyeq (1-z_n)^{-1/2} I_p.
\end{align}
Moreover, for sufficiently large $n$, $z_n < 1/2$ and it follows that ${\opnorm{ I_p - \Psi^{-1/2} } \le 1 - (1-z_n)^{-1/2} \le 2z_n}$. Therefore, taking $C_3'(\kappa, \Rx, \sigma, \varsigma) = 2\kappa^2 c_4$, it follows that with probability greater than $1 - O(n^{-2})$,
\begin{align}
    \opnorm{ I_p - \Psi^{-1/2} } \le C_3'(\kappa, \Rx, \sigma, \varsigma) \cdot n^{-1/2}.
\end{align}


\section{Auxiliary Results for Classical Multidimensional Scaling}
\label{proof:auxiliary-mds-results}

The lemma below collects probabilistic bounds for various quantities which appear in other proofs.

\cbox{black!5}
{%
\begin{lemma}\label{lem:probabilistic-bounds}
    Consider the centered configuration $X = U \Lambda^{1/2} Q \in \R^{n \times p}$. Suppose $D = \Del + \Eps$ satisfying~\ref{assumption:compact}--\ref{assumption:noise}. Let $\hU \hL \hU\tr$ be the rank-$p$ spectral decomposition of $\Dc = -\tfrac{1}{2} H D H$ and
    \begin{align}
        \hX = \mds(D,p) = \hU \hL^{1/2}.
    \end{align}
    denote the output of classical multidimensional scaling. Let $\hQ \in \orth{p}$ be the Procrustes alignment from \eqref{eq:procrustes}, $\hP := \hQ\tr Q$, and denote the Frobenius-optimal rigid transformation of $X$ by $\hg(X) = X \hP\tr.$ Let $\h{\Del} = (\h{\del}_{ij})$ where $\h{\del}_{ij} = \norm{\hx_i - \hx_j}^2$, $E := D - \h{\Del}$
    denote the residual matrix. Let $\Sigma := (\sigma_{ij}^2)$ denote the matrix of noise variances, $\h{\Sigma} := (e_{ij}^2)$ and $\hat\sigma^2 = {\binom{n}{2}}\inv\sum_{i < j}(e_{ij} - \bar{e})^2$. For each $i \in [n]$, let $\Sigma_i := \diag(\sigma_{i1}^2, \dots, \sigma_{in}^2)$ and $\hat\Sigma_i = \diag(e_{i1}^2, \dots, e_{in}^2)$.\\
    
    Then, for sufficiently large $n$, with probability at least $1 - O(n^{-2})$ the following statements hold:
    \begin{enumerate}[label=\textup{(\roman*)}, ref=\cref{lem:probabilistic-bounds}\,(\roman*)]
        \item\label{bound-1} $\opnorm{\hX - \hg(X)} \lesssim c_1(\pars)$.
        \item\label{bound-2} $\ttinf{\hX - \hg(X)} \lesssim c_2(\pars)\sqrt{\log{n}/n}$.
        \item\label{bound-3} $\Opnorm\big{\hX\tr\hX - \hg(X)\tr\hg(X)} \lesssim c_3(\pars)\sqrt{n}$.
        \item\label{bound-4} $\Opnorm\big{(\hX\tr\hX)^{-1} - (\hg(X)\tr\hg(X))^{-1}} \lesssim c_4(\pars)n^{-3/2}$
        \item\label{bound-5} $\max_{i, j}\abs\big{\h{\del}_{ij} - {\del}_{ij}} \lesssim c_5(\pars) \sqrt{\log{n}/n}$.
        \item\label{bound-6} $\maxnorm{\h{\Sigma} - \Sigma} \lesssim \Msigma^2\log^2{n}.$
        \item\label{bound-7} For $i \in [n]$, $\h{\Sigma}_i := \diag(\hat\Sigma_{i, *})$ and $\Sigma_i := \diag(\Sigma_{i, *})$, there exists $c_7 \equiv c_7(\pars)$ such that
        \begin{align}
            \max_{i \in [n]}\Opnorm\big{X\tr\,(\h{\Sigma}_i - \Sigma_i)\,X} &\lesssim c_7 \cdot \sqrt{n \log{n}}\\
            \max_{i \in [n]}\Opnorm\big{U\tr\,(\h{\Sigma}_i - \Sigma_i)\,U} &\lesssim c_7 \cdot \sqrt{\frac{\log{n}}{n}}.
        \end{align}
        \item\label{bound-8} For the same setup as above, there exists $c_8 \equiv c_8(\pars)$ such that
        \begin{align}
            \max_{i \in [n]}\Opnorm\big{\hX\tr\,\h{\Sigma}_i\,\hX - \hg(X)\tr \,\Sigma_i\, \hg(X)} &\lesssim c_8 \log^2{n}\sqrt{n}\\
            \max_{i \in [n]}\Opnorm\big{\hU\tr\,\h{\Sigma}_i\,\hU - (U\hQ)\tr \,\Sigma_i\, (U\hQ)\tr} &\lesssim c_8 {\frac{\log^2{n}}{\sqrt{n}}}.
        \end{align}
        \item\label{bound-9} $\opnorm{U\tr\Eps U} \lesssim c_9(\pars) \cdot \sqrt{\log{n}}$.
        \item\label{bound-10} If $(\eps_{ij})$ are \iid{} with variance $\sigma^2$, then $\abs{\hat\sigma - \sigma} \lesssim c_{10} \log{n}/n.$
    \end{enumerate}
\end{lemma}
}

The proof of \cref{lem:probabilistic-bounds} is deferred to \cref{proof:lem:probabilistic-bounds}. 

\noindent The next result is a decomposition for $\hg(X)\!\!-\!\!\hX$ from.
{%
\begin{lemma}[Lemma~9 of \citealp{vishwanath2025minimax}]\label{lem:reconstruction-error}
    Under the same setup as \cref{lem:probabilistic-bounds},
    \begin{align}
        \hX - \hg(\X) = 
        &(\Dc - \Delc)U\L^{-1/2}\hQ\\
        &+(\Dc - \Delc)(\hU - U\hQ)\hL^{-1/2} + U\L (U\tr\hU - \hQ)\hL^{-1/2} + \Dc \U(\hP\hL^{-1/2}-\L^{-1/2}\hP).
    \end{align}
\end{lemma}
}

\subsection{Proof of \cref{prop:decomposition}}
\label{proof:prop:decomposition}

\begin{proof}
    Since $X$ is assumed to be centered in \ref{assumption:compact}, note that $HX = X$ and the expression for $\Om_i$ simplifies to 
    \begin{align}
        \Om_i = \frac{n}{4} (X\tr X)\inv (X\tr \Sigma_i X) (X \tr X)\inv.\label{eq:Omega-i-appendix-centered}
    \end{align}    
    Similarly, from \cref{rem:centering} this also implies that $\hg(X) = X\hP\tr$ and $\hg\inv(\hX) = \hX\hP$ for $\hP = \hQ\tr Q$. From \cref{lem:reconstruction-error}, we have
    \begin{align}
        X - \hX\hP = 
        &(\Delc - \Dc)U\L^{-1/2}\hQ\hP\\
        &+\underbrace{(\Delc - \Dc)(\hU - U\hQ)\hL^{-1/2}\hP}_{=:\zeta^{(1)}} + \underbrace{U\L (\hQ-U\tr\hU)\hL^{-1/2}\hP}_{=:\zeta^{(2)}} + \underbrace{\Dc \U(\L^{-1/2}\hP-\hP\hL^{-1/2})\hP}_{=:\zeta^{(3)}}.
    \end{align}
    For the first term, since $XQ\tr = U\L^{1/2}$ and $\hQ\hP = \hQ\hQ\tr Q = Q$, 
    \begin{align}
        \qty\Big(\Delc - \Dc)\qty(U\L^{-1/2})\hQ\hP 
        &= \qty(\half H\Eps H)\qty\Big(XQ\tr\L\inv)Q\\ 
        &= \half \qty(I - \frac{J}{n}) \qty\big(\Eps X) (X\tr X)\inv\\
        &= \frac{1}{2} \Eps\X \qty\big(X\tr\X)\inv - \underbrace{\frac{J}{2n}(\Eps X)(X\tr\X)\inv}_{=:\zeta^{(4)}},
    \end{align}
    where in the second line we used the fact that $H\X = X$ since $X$ is assumed to be centered, and $Q\tr \L\inv Q = (X\tr X)\inv$. For $\zeta := \zeta^{(1)} + \zeta^{(2)} + \zeta^{(3)} + \zeta^{(4)}$, we have
    \begin{align}
        X - \hX\hP = \frac{1}{2n} \Eps X \qty(\tfrac{X\tr\X}{n})\inv + \zeta.\label{eq:reconstruction-error-matrix-form}
    \end{align}
    From the intermediate calculations appearing in the proof of \citet[Theorem~3]{vishwanath2025minimax}, it can be shown that the residual matrix satisfies $\ttinf{\zeta} = o_p(n^{-1/2})$.
    \cbox{black!5}
    {%
    \begin{lemma}\label{lem:intermediate-quantities}
        For $\zeta := \zeta^{(1)} + \zeta^{(2)} + \zeta^{(3)} + \zeta^{(4)}$ defined above, with probability greater than $1 - O(n^{-2})$,
        \begin{align}
            \ttinf{\zeta} \lesssim c'(\pars)\frac{\sqrt{\log{n}}}{n},
        \end{align}
        where $c'(\pars) > 0$ is a term which depends only on $\pars$.
    \end{lemma}
    }
    The proof of \cref{lem:intermediate-quantities} is relegated to \cref{proof:lem:intermediate-quantities}. In other words, from \cref{eq:reconstruction-error-matrix-form,lem:intermediate-quantities}, for each $i \in [n]$
    \begin{align}
        \sqrt{n}\qty\Big(x_i - \hg\inv(\hx_i)) = \Upsilon_i + \sqrt{n}\zeta_i 
    \end{align}
    where
    \begin{align}
        \Upsilon_i := \frac{1}{\sqrt{n}} \sum_{k \in [n]} \eps_{ik} \cdot \tfrac{1}{2} \qty( \tfrac{X\tr X}{n} )\inv x_k \qq{and}
        \max_{i \in [n]}\norm{\zeta_i} = \Op\qty(\tfrac{c'(\pars)\log{n}}{n}).
    \end{align}
    In the display above, we used the fact that each row of $\hX\hP \in \Rnp$ is $\hP\tr \hx_i = \hg\inv(\hx_i)$. Since $\E(\epsilon_{ik}) = 0$, it follows that $\E(\Upsilon_i) = 0_p$ and
    \begin{align}
        \Var(\Upsilon_i) = \E(\Upsilon_i\Upsilon_i\tr) 
        &= \frac{1}{n}\sum_{k \in [n]}\sigma_{ik}^2 \cdot \tfrac{1}{2}\qty(\tfrac{X\tr X}{n})\inv x_k x_k\tr \cdot \tfrac{1}{2}\qty(\tfrac{X\tr X}{n})\inv\\ 
        &= \frac{n}{4} \qty\big( {X\tr X})\inv \qty\big( X\tr \Sigma_i X ) \qty\big({X\tr X})\inv = {\Om_i},\label{eq:variance-yi-appendix}
    \end{align}
    where $\Sigma_i := \diag(\sigma_{i1}^2, \dots, \sigma_{in}^2)$. As noted in \cref{eq:Omega-i-appendix-centered}, this is the expression for $\Om_i$ when $X$ is centered. If $X$ is not centered, the same analysis above gives $\Om_i$ as in \cref{eq:Omega-i}. Therefore, for $Y_i := \Om_i^{-1/2}\Upsilon_i$ and $R_i := \sqrt{n}\Om_i^{-1/2}\zeta_i$,
    \begin{align}
        \sqrt{n}\Om_i^{-1/2}(x_i - \hg\inv(\hx_i)) = Y_i + R_i.
    \end{align}
    Moreover, from \cref{eq:variance-yi-appendix}, we also have that $\E(Y_i)= 0$ and $\Var(Y_i) = I_p$. Rewriting $X = U \Lambda^{1/2} Q\tr$, it also follows that
    \begin{align}
    \Om_i
    &= \frac{n}{4} (\Lambda^{-1/2}Q)\tr  \qty\big(U\tr S_i U) (\Lambda^{-1/2} Q);
\end{align}
from the assumption that $\min_{i \in [n]} U^TS_iU \succcurlyeq \msigma^2 I_p$ in \ref{noise-3}, this implies that for all $i \in [n]$
\begin{align}
    ({\msigma^2}/{4\kappa^2}) I_p \preccurlyeq \Om_i \preccurlyeq ({\Msigma^2} {4\kappa^2}) I_p\qc{}\qquad ({2}/{\Msigma\kappa})I_p \preccurlyeq {\Om_i^{-1/2}} \preccurlyeq (2\kappa/\msigma)I_p,\label{eq:omega-bounds}
\end{align}
and
\begin{align}
    \max_{i, k}\norm{\theta_{ik}} 
    \le \frac{1}{2}\max_{i}\opnorm{\Om_i^{-1/2}} \cdot \opnorm{(\tfrac{X\tr X}{n})\inv} \cdot \ttinf{X}
    \le \frac{\kappa^3\Rx}{2\msigma} =: \const_0(\kappa, \Rx, \sigma),
\end{align}
which proves \cref{eq:theta-norm}. Lastly, note that $\ttinf{\zeta} \lesssim {c'(\kappa, \Rx, \sigma)\cdot\sqrt{\log{n}}}/{n}$ implies that
\begin{align}
    \max_{i \in [n]}\norm{R_i} \le \max_{i \in [n]}\opnorm{\Om_i^{-1/2}} \cdot \sqrt{n}\ttinf{\zeta} \lesssim  C'_1(\kappa, \Rx, \sigma) \frac{\log{n}}{\sqrt{n}} 
\end{align}
with probability greater than $1-O(n^{-2})$ and for $C'_1(\pars) := c'(\pars) \cdot \kappa/\msigma$.
\end{proof}

\subsection{Proof of \cref{lem:probabilistic-bounds}}
\label{proof:lem:probabilistic-bounds}

\begin{proof}
    In the interest of clarity, throughout we will write $\tX := \hg(X)$
    and write $c_{\square}$ to denote constants $c_{\square}(\pars)$ which depend on $\pars$. 
    
    Note that since $s_k(\tX) = s_k(X)$ for all $k \in [p]$ since $\tX = X \hP\tr$ for $\hP \in \orth{p}$. \ref{bound-1} and \ref{bound-2} are Theorem~2 and Theorem~3 from \cite{vishwanath2025minimax}, respectively, with $c_1 = \kappa \Msigma$ and $c_2 = \Msigma \kappa^2(\kappa + \Rx)$.

    \noindent $\bullet$ \textit{Proof of \ref{bound-3}.}\quad Using the triangle inequality, we have
    \begin{align}
        \Opnorm{\hX\tr\hX - \hg(X)\tr\hg(X)} &\le \Opnorm\big{(\tX - \hX)\tr\hX} + \Opnorm\big{\tX\tr(\hX - \tX)} \le \opnorm{\hX - \tX} \qty( \opnorm{\hX} + \opnorm{\tX} ).\label{eq:bound-3-decomposition}
    \end{align}
    From \ref{bound-1} and $\opnorm{\tX} \le \kappa\sqrt{n}$ from \ref{assumption:compact}, and using $\opnorm{\hX} \le \opnorm{\tX} + \opnorm{\hX - \tX}$, we obtain
    \begin{align}
        \Opnorm{\hX\tr\hX - \hg(X)\tr\hg(X)} 
        &\lesssim c_1\qty( 2\kappa\sqrt{n} + c_1 )\\ 
        &\lesssim 2c_1\kappa\sqrt{n} =: c_3\sqrt{n}.\FINEQ
    \end{align}

    \noindent\textit{Proof of \ref{bound-4}.}\quad By rewriting the difference, we have
    \begin{align}
        \Opnorm{  (\hX\tr\hX)^{-1} - (\tX\tr\tX)^{-1}} 
        &= \Opnorm{ (\tX\tr\tX)^{-1} \cdot \qty( \tX\tr\tX - \hX\tr\hX ) \cdot (\hX\tr\hX)^{-1} }\\
        &\le \frac{\opnorm{ \tX\tr\tX - \hX\tr\hX }}{s_p(\tX\tr\tX) s_p(\hX\tr\hX)}.
    \end{align}
    Note that ${n}/\kappa^2 \le \lambda_p(\tX\tr\tX)$, and for sufficiently large $n$ such that $n/2\kappa^2 > \sqrt{n}c_3$, and,
    \begin{align}
        s_p(\hX\tr\hX) \ge s_p(\tX\tr\tX) - \opnorm{\hX\tr\hX - \tX\tr\tX} \ge \frac{n}{\kappa^2} - \sqrt{n}c_3 \gtrsim \frac{n}{\kappa^2}
    \end{align}
    by an application of Weyl's inequality followed by \ref{bound-3}. It follows that
    \begin{align}
        \Opnorm{  (\hX\tr\hX)^{-1} - (\tX\tr\tX)^{-1}} \lesssim \frac{c_3\sqrt{n}}{n^2/\kappa^4} \lesssim \frac{\kappa^4c_3}{n^{3/2}} =: \frac{c_4}{n^{3/2}}.\FINEQ
    \end{align}

    \noindent $\bullet$ \textit{Proof of \ref{bound-5}.}\quad Since $\del_{ij} = \norm{x_i - x_j}^2 = \norm{\tx_i - \tx_j}^2$, we have
    \begin{align}
        \h{\del}_{ij} - {\del}_{ij} = \norm{\h{x}_i - \h{x}_j}^2 - \norm{x_i - x_j}^2 = \norm{ \tx_i - \tx_j + (\h{x}_i - \tx_i) - (\h{x}_j - \tx_j) }^2 - \norm{\tx_i - \tx_j}^2.
    \end{align}
    By expanding the square and using the Cauchy-Schwarz inequality, we obtain
    \begin{align}
        \abs{\h{\del}_{ij} - {\del}_{ij}} \le \norm{\h{x}_i - \tx_i}^2 + \norm{\h{x}_j - \tx_j}^2 + 2\norm{\tx_i - \tx_j}\qty\Big(\norm{\h{x}_i - \tx_i} + \norm{\h{x}_j - \tx_j}).
    \end{align}
    From \ref{assumption:compact}, note that $\max_{ij}\norm{\tx_i - \tx_j} \le 2\ttinf{X} \le 2\Rx$. Using \ref{bound-2}, we now have
    \begin{align}
        \max_{i, j}\abs{\h{\del}_{ij} - {\del}_{ij}} 
        &\le \max_{i \in [n]}\norm{\h{x}_i - \tx_i}^2 + \max_{j \in [n]}\norm{\h{x}_j - \tx_j}^2 + 2\Rx\max_{i, j}\qty(\norm{\h{x}_i - \tx_i} + \norm{\h{x}_j - \tx_j})\\
        &\lesssim 2c_2^2 \frac{\log{n}}{n} + 2\Rx c_2\sqrt{\frac{\log{n}}{n}} \lesssim c_5\sqrt{\frac{\log{n}}{n}},
    \end{align}
    for $c_5 = 2\Rx c_2$.\FIN

    \noindent $\bullet$ \textit{Proof of \ref{bound-6}.}\quad Using the definition of $E = D - \h{\Del}$, for each $i \neq j \in [n]$ we have
    \begin{align}
        \h{\sigma}_{ij}^2 - \sigma_{ij}^2 
        = e_{ij}^2- \sigma_{ij}^2 = (d_{ij} - \h{\del}_{ij})^2- \sigma_{ij}^2 
        &= \qty(\del_{ij} + \eps_{ij} - \h{\del}_{ij})^2- \sigma_{ij}^2\\ 
        &= \qty(\eps_{ij}^2- \sigma_{ij}^2) + 2\eps_{ij}(\hat\del_{ij} - \del_{ij}) + \qty(\hat\del_{ij} - \del_{ij})^2.\label{eq:bound-6-decomposition}
    \end{align}
    From \ref{bound-5}, we have
    \begin{align}
        \max_{i, j}\abs{\h{\del}_{ij} - \del_{ij}} \lesssim c_5\sqrt{\frac{\log{n}}{n}} \qq{and} \max_{i, j}(\hat\del_{ij} - \del_{ij})^2 \lesssim c_5^2\frac{\log{n}}{n}.\label{eq:bound-6-1}
    \end{align}
    Since $\eps_{ij}$ are uniformly $\Msigma$-sub-exponential, from \cref{lem:orlicz-norms}, for all $i, j \in [n]$ and $t > 0$
    \begin{align}
        \pr\qty( \abs{\eps_{ij}} \ge t ) \le 2e^{-t/\Msigma}.
    \end{align}
    Taking a union bound over all $i < j \in [n]$, we have
    \begin{align}
        \pr\qty( \max_{i < j} \abs{\eps_{ij}} \ge t ) \le \sum_{i < j}\pr\qty\Big(\abs{\eps_{ij}} \ge t) \le n^2 \cdot 2e^{-t/\Msigma}.
    \end{align}
    Setting $t = 4\Msigma\log{n}$, it follows that with probability at least $1 -2n^{-2}$,
    \begin{align}
        \max_{i < j} \abs{\eps_{ij}} \le 4\Msigma\log{n}.\label{eq:bound-6-2}
    \end{align}
    Finally, from \cref{lem:orlicz-norms} for $\alpha = 1/2$ and for all $i, j \in [n]$,
    \begin{align}
        \norm{\eps_{ij}^2 - \sigma_{ij}^2}_{\psi_\alpha} \le \norm{ \eps_{ij} + \sigma_{ij} }_{\psi_1} \cdot \norm{ \eps_{ij} - \sigma_{ij} }_{\psi_1} \le 2C\Msigma^2,\label{eq:bound-6-3-norm}
    \end{align}
    and using a similar argument as above, we have that for all $t > 0$
    \begin{align}
        \pr\qty( \abs{\eps_{ij}^2 - \sigma_{ij}^2} \ge t ) \le  2e^{-(t/(2C\Msigma^2))^\alpha}.
    \end{align}
    Taking a union bound over all $i < j \in [n]$ and setting $t = 32C\Msigma^2\log^2{n}$, it follows that with probability at least $1 - 2n^{-2}$,
    \begin{align}
        \max_{i < j} \abs{\eps_{ij}^2 - \sigma_{ij}^2} \lesssim \Msigma^2\log^2{n}.\label{eq:bound-6-3}
    \end{align}
    Combining \cref{eq:bound-6-1,eq:bound-6-2,eq:bound-6-3}, we have that with probability at least $1 - O(n^{-2})$,
    \begin{align}
        \max_{i, j} \abs{\h{\sigma}_{ij}^2 - \sigma_{ij}^2} \lesssim \qty(\Msigma^2\log^2{n}) + \qty(\frac{\Msigma c_5\log^{3/2}}{\sqrt{n}}) + \qty(\frac{c_5^2\log{n}}{n}) \lesssim c_6\log^2{n}
    \end{align}
    for $c_6 = \Msigma^2$.
    \FIN

    \noindent $\bullet$ \textit{Proof of \ref{bound-7}.}\quad Using the decomposition in \cref{eq:bound-6-decomposition} we have
    \begin{align}
        X\tr(\h{\Sigma}_i - \Sigma_i) X 
        &= \sum_{j \in [n]} (\h{\sigma}_{ij}^2 - \sigma_{ij}^2) x_jx_j\tr\\ 
        &= \sum_{j \in [n]}(\eps_{ij}^2 - \sigma_{ij}^2)x_jx_j\tr + 2\sum_{j \in [n]}\eps_{ij}(\hat\del_{ij} - \del_{ij})x_jx_j\tr + \sum_{j \in [n]}(\h{\del}_{ij} - \del_{ij})^2x_jx_j\tr.\label{eq:bound-7-decomposition}
    \end{align}
    We use \cref{prop:matrix-bernstein} to bound the first two terms. For $\alpha=1/2$ and $q=p$ and from \cref{eq:bound-6-3-norm}, we have
    \begin{align}
        K=\max_{j \in [n]}\norm{\eps_{ij}^2 - \sigma^2_{ij}}_{\psi_\alpha} \lesssim \Msigma^2\qc{}
        M = \max_{j \in [n]}\norm{x_jx_j\tr} = \max_{j \in [n]}\norm{x_j}^2 \le \Rx^2
    \end{align}
    and using \citet[Proposition~2.7.1]{vershynin2018high},
    \begin{align}
        \gamma^2 = \Opnorm\Big{\textstyle{\sum\limits_{j \in [n]} \E(\eps_{ij}^2 - \sigma^2_{ij})^2 (x_jx_j\tr)^2}} \lesssim 4\norm{\eps_{ij}^2 - \sigma^2_{ij}}^2_{\psi_\alpha} \cdot \max_{j}\norm{x_j}^2  \cdot\opnorm{X\tr X} \lesssim n \cdot \Msigma^4\Rx^2 \kappa^2.
    \end{align}
    Setting $t = 3\log{n}$ in \cref{prop:matrix-bernstein}, it follows that with probability at least $1 - O(n^{-3})$,
    \begin{align}
        \opnorm{\sum_{j \in [n]}(\eps_{ij}^2 - \sigma_{ij}^2)x_jx_j\tr} \lesssim \Msigma^2\Rx\kappa\sqrt{n \log{n}}.
    \end{align}
    Taking a union bound over all $i \in [n]$, we have that with probability at least $1 - O(n^{-2})$,
    \begin{align}
        \max_{i \in [n]}\Opnorm\Big{\sum_{j \in [n]}(\eps_{ij}^2 - \sigma_{ij}^2)x_jx_j\tr} \lesssim \Msigma^2\Rx\kappa\sqrt{n \log{n}}.\label{eq:bound-7-term-1}
    \end{align}
    Similarly, for the second term in \cref{eq:bound-7-decomposition}, for the same $M \le \Rx^2$ and with $\alpha=1$,
    \begin{align}
        K = \max_{j in [n]}\norm{\eps_{ij}}_{\psi_1} \le \Msigma\qq{and}
        \gamma^2 = \Opnorm\Big{\textstyle{\sum\limits_{j \in [n]} \E(\eps_{ij}^2) (x_jx_j\tr)^2}} \le n \cdot \Msigma^2\Rx^2\kappa^2,
    \end{align}
    from \cref{prop:matrix-bernstein} and using \ref{bound-5} it follows that with probability at least $1 - O(n^{-2})$,
    \begin{align}
        \max_{i \in [n]}\Opnorm\Big{\sum_{j \in [n]}\eps_{ij}(\hat\del_{ij} - \del_{ij})x_jx_j\tr} \le \max_{i \in [n]}\Opnorm\Big{\sum_{j \in [n]}\eps_{ij}x_jx_j\tr} \cdot \max_{i < j}\abs{\hat\del_{ij} - \del_{ij}} \lesssim \Msigma\Rx\kappa \cdot c_5 \log{n}.
    \end{align}
    Using \cref{eq:bound-6-1} once again: with probability at least $1 - O(n^{-2})$,
    \begin{align}
        \max_i\Opnorm\Big{ \textstyle{\sum\limits_{j \in [n]}}(\h{\del}_{ij} - \del_{ij})^2 x_j x_j\tr} \le \opnorm{X\tr X} \cdot \max_{i < j}(\h{\del}_{ij} - \del_{ij})^2 \lesssim \kappa^2 c_5^2 \log{n}.
    \end{align}
    Using the triangle inequality in \cref{eq:bound-7-decomposition}, plugging in the bounds above and taking $c_7 := \Msigma^2\Rx\kappa$, we get:
    \begin{align}
        \max_i\Opnorm\Big{X\tr(\h{\Sigma}_i - \Sigma_i) X} 
        &\lesssim c_7\sqrt{n \log{n}} \FINEQ
    \end{align}
    The proof for $\opnorm{U\tr(\h{\Sigma}_i - \Sigma_i)U}$ nearly identical. Similar to \cref{eq:bound-7-decomposition}, we have
    \begin{align}
        &\max_{i}\Opnorm{U\tr (\h{\Sigma}_i - \Sigma_i) U}\\ 
        &\qquad\qquad\le \max_{i}\Opnorm\Big{ \sum_j (\eps_{ij}^2 - \sigma_{ij}^2) u_j u_j\tr} + \max_{i}\Opnorm\Big{ \sum_j \eps_{ij} (\hat\del_{ij} - \del_{ij}) u_j u_j\tr} + \max_{i}\Opnorm\Big{ \sum_j (\h{\del}_{ij} - \del_{ij})^2 u_j u_j\tr}.
    \end{align}
    Note that $U\tr U=I_p$ and since $X = U\Lambda^{1/2}Q$, we also have $\ttinf{U} \le \ttinf{X} \opnorm{\Lambda^{-1/2}} \le \Rx\kappa/\sqrt{n}$. Therefore, the only adjustments needed are:
    \begin{align}
        M := \max_{j \in [n]}\norm{u_j}^2 \le \frac{\Rx^2\kappa^2}{n} \qq{and} \gamma^2 := \Opnorm\Big{\textstyle{\sum\limits_{j \in [n]} \E(\eps_{ij}^2 - \sigma^2_{ij})^2 (u_j u_j\tr)^2}} \lesssim \Msigma^4 \cdot \ttinf{U}^2 \cdot \opnorm{U} \le \frac{\Msigma^4\Rx^2\kappa^2}{n}.
    \end{align}
    Following the proof from above now leads to the following bounds with probability at least ${1 - O(n^{-2})}$:
    \begin{align}
        \max_{i \in [n]}\Opnorm\Big{\textstyle\sum_{j \in [n]}(\eps_{ij}^2 - \sigma_{ij}^2) u_j u_j\tr} 
        &\lesssim \Msigma^2\Rx\kappa\sqrt{{\log{n}}/{n}}\\
        \max_{i \in [n]}\Opnorm\Big{\textstyle\sum_{j \in [n]}\eps_{ij}(\hat\del_{ij} - \del_{ij}) u_j u_j\tr} 
        &\lesssim \Msigma\Rx\kappa\sqrt{{\log{n}}/{n}} \cdot \max_{i < j}\abs\big{\h{\del}_{ij} - \del_{ij}}
        \lesssim \Msigma\Rx\kappa \cdot c_5 {\log{n}}/{{n}}\\
        \max_{i \in [n]}\Opnorm{\textstyle\sum_{j \in [n]}(\h{\del}_{ij} - \del_{ij})^2 u_j u_j\tr}
        &\lesssim \max_{i < j}(\h{\del}_{ij} - \del_{ij})^2 \cdot \opnorm{U\tr U} \lesssim c_5^2 \frac{\log{n}}{n}.
    \end{align}
    Combining the bounds above gives the desired result:
    \begin{align}
        \max_i\Opnorm\big{U\tr(\h{\Sigma}_i - \Sigma_i) U} 
        &\lesssim c_7 \sqrt{\frac{\log{n}}{n}}. \FINEQ
    \end{align}

    \noindent $\bullet$ \textit{Proof of \ref{bound-8}.}\quad Rewriting the difference similar to the procedure in the proof for \ref{bound-3} gives:
    \begin{align}
        {\hX\tr\,\h{\Sigma}_i\,\hX - \tX\tr \,\Sigma_i\, \tX} = (\hX - \tX)\tr\,\h{\Sigma}_i\,\hX + \tX\tr\,\h{\Sigma}_i\,(\hX - \tX) + \tX\tr\,(\h{\Sigma}_i - \Sigma_i)\,\tX.\label{eq:bound-8-decomposition}
    \end{align}
    where $\opnorm{\tX} \le \kappa\sqrt{n}$ and with probability at least $1 - O(n^{-2})$, we have
    \begin{align}
        &\opnorm{\hX - \tX} \lesssim c_1\\[5pt]
        &\opnorm{\hX} \lesssim \opnorm{\tX} + \opnorm{\hX - \tX} \lesssim \kappa \sqrt{n}\\[5pt]
        &\max_i\opnorm{\hat\Sigma_i} \lesssim \opnorm{\Sigma_i} + \maxnorm{\hat\Sigma - \Sigma} \lesssim \Msigma^2\log^2{n} \\[5pt]
        &\max_i\opnorm{\tX\tr\,(\hat\Sigma_i - \Sigma_i)\,\tX} \lesssim c_7\sqrt{n\log{n}}.
    \end{align}
    Using the triangle inequality in \cref{eq:bound-8-decomposition} and plugging in the bounds above, we get with probability $1 - O(n^{-2})$,
    \begin{align}
        \max_i\Opnorm\big{\hX\tr\h{\Sigma}_i\hX - \tX\tr \Sigma_i \tX} \lesssim c_8\log^2{n}\sqrt{n}, \FINEQ
    \end{align}
    for $c_8=c_1 \Msigma^2 \kappa$. The proof is identical for the second claim as well. Note that
    \begin{align}
        \opnorm{\hU\tr\,\h{\Sigma}_i\,\hU - (U\hQ)\tr \,\Sigma_i\, (U\hQ)} = \opnorm{(\hU - U\hQ)\tr\,\h{\Sigma}_i\,\hU + (U\hQ)\tr\,\h{\Sigma}_i\,(\hU - U\hQ) + (U\hQ)\tr\,(\h{\Sigma}_i - \Sigma_i)\,(U\hQ)},
    \end{align}
    where $\opnorm{U\hQ} = \opnorm{\hU} = 1$, and from \citet[Lemma~14]{vishwanath2025minimax}, we have
    \begin{align}
        \opnorm{\hU - U\hQ} &\lesssim c'/\sqrt{n}\\
        \max_i\opnorm{\hat\Sigma_i} &\lesssim \Msigma^2\log^2{n}\\
        \max_i\opnorm{(U\hQ)\tr\,(\hat\Sigma_i - \Sigma_i)\,(U\hQ)} &\lesssim c_7\sqrt{\frac{\log{n}}{n}}
    \end{align}
    with probability at least $1 - O(n^{-2})$. Plugging in the bounds above, for $c_8' = \Msigma^2c'$, we get
    \begin{align}
        \max_i\Opnorm\big{\hU\tr\h{\Sigma}_i\hU - (U\hQ)\tr \Sigma_i (U\hQ)} \lesssim c'_{8}{\frac{\log^2{n}}{\sqrt{n}}}. \FINEQ
    \end{align}

    \noindent $\bullet$ \textit{Proof of \ref{bound-9}.}\quad Since $(\eps_{ij}) \in \Rnn$ is symmetric, we have
    \begin{align}
        U\tr\Eps U = \sum_{i < j} \eps_{ij} (u_i u_j\tr + u_ju_i\tr) \in \Rpp
    \end{align}
    is the sum of ${n \choose 2}$ independent matrices with sub-exponential entries. Similar to \ref{bound-7}, we will again use \cref{prop:matrix-bernstein} to bound the operator norm. To this end, we have, 
    \begin{align}
        K = \max_{i < j}\norm{\eps_{ij}}_{\psi_1} \le \Msigma
        \qc{}
        M = \max_{i < j}\Opnorm{u_i u_j\tr + u_ju_i\tr} \le 2\max_{i < j}\Opnorm{u_iu_j\tr} = 2\max_{i < j}\norm{u_i}\norm{u_j} \le 2\Rx^2\kappa^2/n,
    \end{align}
    and
    \begin{align}
        \gamma^2 = \Opnorm\Big{\textstyle\sum\limits_{i < j} \E(\eps_{ij}^2) (u_i u_j\tr + u_ju_i\tr)^2} \le {{n \choose 2}} \cdot \Msigma^2 \cdot 4\max\limits_{i < j}\opnorm{(u_i u_j\tr)^2} \lesssim \Msigma^2\Rx^4\kappa^4.
    \end{align}
    Using \cref{prop:matrix-bernstein}, it follows that for all $t > 0$, with probability at least $1 - 2e^{-t}$,
    \begin{align}
        \opnorm{U\tr \Eps U} \lesssim \Msigma\Rx^2\kappa^2\sqrt{t} + \frac{\Msigma\Rx^2\kappa^2}{n} t\log{n}.
    \end{align}
    Setting $t = 2\log{n}$, it follows that with probability at least $1 - O(n^{-2})$,
    \begin{align}
        \opnorm{U\tr \Eps U} \lesssim \Msigma\Rx^2\kappa^2\sqrt{\log{n}} =: c_{9}\sqrt{\log{n}}.
    \end{align}

    \noindent $\bullet$ \textit{Proof of \ref{bound-10}.}\quad Let $N:= \binom{n}{2} \asymp n^2$. From the decomposition in \cref{eq:bound-6-decomposition}, we have
    \begin{align}
    \abs{\hat\sigma^2 - \sigma^2} 
    &\le N\inv \abs\Big{\sum_{i < j} (\eps_{ij}^2 - \sigma^2) + 2 \sum_{i < j} \eps_{ij}(\hat\del_{ij} - \del_{ij}) + \sum_{i < j} (\hat\del_{ij} - \del_{ij})^2} + \bar{e}^2.\label{eq:sigma-decomposition-1}
    \end{align}
    For the first two terms, we use \cref{prop:matrix-bernstein} with $\alpha=1/2$, $p=q=1$. Specifically, for $\xi_{ij} = N\inv (\eps_{ij}^2-\sigma^2_{ij})$, it follows that for all $t > 0$ and 
    with probability at least $1 - 2e^{-t}$,
    \begin{align}
    \abs\Big{N\inv\sum\limits_{i < j} (\eps_{ij}^2 - \sigma^2)} \lesssim \gamma \sqrt{t} + MK (t\log{n})^2,
    \end{align}
    where, for $\alpha=1/2$, $A_i=1$, $M=1$,
    \begin{align}
    K = \max\limits_{i < j}\norm{\xi_{ij}}_{\psi_{\alpha}} \lesssim N\inv\norm{\eps_{ij}^2 - \sigma^2}_{\psi_{\alpha}} \lesssim N\inv\sigma^2,
    \qq{and}
    \gamma^2 = \sum\limits_{i < j}\E(\xi_{ij}^2) \lesssim N\inv \sigma^4.
    \end{align}
    Setting $t = 2\log{n}$ in \cref{prop:matrix-bernstein}, we get that with probability at least $1 - 2n^{-2}$,
    \begin{align}
    \abs\Big{N\inv\sum\limits_{i < j} (\eps_{ij}^2 - \sigma^2)} \lesssim \sigma^2N^{-1/2} \sqrt{\log{n}} \lesssim \sigma^2\frac{\sqrt{\log{n}}}{n}.
    \end{align}
    A similar analysis for the second term using \cref{prop:matrix-bernstein} with $\xi_{ij} = N\inv\eps_{ij}$, $\alpha=1$, $A_i=1$, $M=1$ gives:
    \begin{align}
    \abs\Big{N\inv\sum\limits_{i < j} \eps_{ij}} \lesssim \sigma N^{-1/2}\sqrt{\log{n}} \lesssim \sigma\frac{\sqrt{\log{n}}}{n}
    \end{align}
    with probability at least $1 - 2n^{-2}$. From \ref{bound-5}, we also have that with probability at least $1 - O(n^{-2})$,
    \begin{align}
    \max_{i < j}\abs\big{\h{\del}_{ij} - \del_{ij}} \lesssim c_5 \sqrt{\log{n}/n} \qq{and} \max_{i < j}\abs\big{\h{\del}_{ij} - \del_{ij}}^2 \lesssim c_5^2 \log{n}/n.
    \end{align}
    An identical analysis also follows for $\bar{e}$: with probability greater than $1 - O(n^{-2})$,
    \begin{align}
        \bar{e} = N\inv \sum_{i < j}e_{ij} = N\inv \sum_{i < j}(\del_{ij} - \h{\del}_{ij}) + N\inv \sum_{i < j}\eps_{ij} \;\;\lesssim\;\; c_5\sqrt{\frac{\log{n}}{n}} + \sigma\frac{\sqrt{\log{n}}}{n},
    \end{align}
    which implies that $\bar{e} \lesssim \sqrt{\log{n}/n}$, and, therefore, $\bar{e}^2 \lesssim \log{n}/n$. Plugging these bounds into \cref{eq:sigma-decomposition-1}, we get that with probability at least $1 - O(n^{-2})$,
    \begin{align}
    \abs{\hat\sigma^2 - \sigma^2} \lesssim \qty(\sigma\frac{\sqrt{\log{n}}}{n}) + \qty(c_5 \sqrt{\frac{\log{n}}{n}} \cdot \frac{\sigma\log{n}}{n}) + c_5^2\frac{\log{n}}{n} \lesssim c_{10}\frac{\log{n}}{n}.
    \end{align}
    for $c_{10} := c_5^2$.
\end{proof}

\subsection{Proof of \cref{lem:intermediate-quantities}}
\label{proof:lem:intermediate-quantities}

Note that from assumption \ref{noise-2} and Proposition~2.7.1 of \cite{vershynin2018high}, ${(\E|\eps_{ij}|^4)^{1/4} \le 4C \Msigma}$ for some absolute constant $C > 0$. Therefore, any appearance of $\sigma$ in \citep{vishwanath2025minimax} (which we refer to as \vac{} henceforth) can be replaced with $\Msigma$.

From the proof of Theorem~3 in Section~7.6 and from Lemma~14 of \vac{} it was already shown that on the event $\qty{ \opnorm{\Delc - \Dc} \lesssim \sqrt{n}\Msigma }$, which happens with probability greater than $1 - O(n^{-2})$, it also holds that:
\begin{align}
    \ttinf{\zeta^{(2)}} \lesssim \frac{c'_2}{n},
    &\qq{and}(n/\kappa^2) I_p \preccurlyeq  \hL \preccurlyeq (n\kappa^2) I_p,\\
    \ttinf{\Dc U} \lesssim c'_3 \sqrt{n}\qc{}
    \opnorm{\hU - U\hQ} \lesssim \frac{c'_4}{\sqrt{n}}\qc{}
    &\opnorm{U\tr\hU - \hQ} \lesssim \frac{c'_5}{n}\qc{}
    \opnorm{U\tr\hU - \hQ} \lesssim \frac{c'_5}{n}.\label{eq:intermediate-quantities-1}
\end{align}

$\bullet$ To bound $\ttinf{\zeta^{(3)}}$, we need a slightly stronger bound than that established in Eq.~(37) of \vac{}. From Eq.(88) of their work and the discussion immediately following it, note that
\begin{align}
    \opnorm{\hQ\hL^{-1/2} -\L^{-1/2}\hQ} &\lesssim \frac{c'_6}{n^{3/2}}\opnorm{\hQ\hL - \L\hQ}\\ 
    &\lesssim \frac{c'_6}{n^{3/2}}\qty\Big( \opnorm{\hQ - U\tr\hU} \cdot (\opnorm{\hL} + \opnorm{\L}) + \opnorm{U\tr(\Dc - \Delc)U} ).\label{eq:intermediate-quantities-2}
\end{align}
Using \cref{eq:intermediate-quantities-1}, it follows that
$
    \opnorm{\hQ - U\tr\hU} \cdot (\opnorm{\hL} + \opnorm{\L}) \lesssim ({c'_5}/{n}) \cdot 2\kappa^2n \lesssim \kappa^2 c'_5,
$
and using \ref{bound-9} instead:
\begin{align}
    \opnorm{U\tr(\Dc - \Delc)U}  \lesssim c'_7 \sqrt{\log{n}},
\end{align}
where, in the first inequality, we used the fact that $\opnorm{\cdot}$ is unitarily invariant and that $HX = X$. Plugging these into \cref{eq:intermediate-quantities-2} and using the bound on $\Dc U$ from \cref{eq:intermediate-quantities-1}, we get
\begin{align}
    \ttinf{\zeta^{(3)}} 
    &\lesssim \ttinf{\Dc U} \cdot \opnorm{\hQ\hL^{-1/2} -\L^{-1/2}\hQ}\\ 
    &\lesssim {c_3'}\sqrt{n} \cdot \frac{c'_6}{n^{3/2}} \qty(c_5'\kappa^2 + c_7' \sqrt{\log{n}}) \lesssim c_3' c_6' c_9\frac{\sqrt{\log{n}}}{n}.\label{eq:zeta3}
\end{align}

$\bullet$ For $\zeta^{(1)}$, in Eq.~(34) of \vac{} a bound on $\opnorm{\hU - U\hQ}$ was used for $\ttinf{\zeta^{(1)}}$. The result is improved if we use a bound on $\ttinf{\hU - U\hQ}$ instead. To this end, we follow the proof of Theorem~4.7 from \cite{cape2019two}. We are somewhat terse since the proof below is nearly identical. In particular, using the decomposition in \citep[Corollary~3.3]{cape2019two} followed by an application of the triangle inequality, we have
\begin{align}
    \ttinf{\hU - U\hQ} 
    &\lesssim \ttinf{(I - UU\tr){(\Dc - \Delc)}{U \hQ\hL\inv}}
    \\
    &\quad + \ttinf{{(I - UU\tr)}{(\Dc - \Delc)(\hU - U\hQ)\hL\inv}}
    \tag{$\lesssim c''_2 \cdot \sqrt{n} \cdot n^{-1/2} \cdot n^{-1}$}\\
    &\quad + \ttinf{{(I - UU\tr)}{\Delc}(\hU - UU\tr\hU)\hL\inv}
    \tag{$=0$}\\
    &\quad + \ttinf{U(U\tr\hU - \hQ)}
    \tag{$\lesssim c''_3 n^{-1/2} \cdot n^{-1}$},
\end{align}
where we used \cref{eq:intermediate-quantities-1} for the second term, $\ttinf{U} \le \Rx/\kappa\sqrt{n}$ in the fourth term, and $\Delc(\hU- UU\tr\hU) = U\L(U\tr\hU - U\hU)= 0$ in the third term; see, also, Section 6.10 of \citep{cape2019two} where this term is zero. For the first term, writing 
$$
(\Dc - \Delc)U = (\Dc - \Delc)U\L^{1/2} Q\tr Q \L^{-1/2} = (\Dc - \Delc)X Q\L^{-1/2},
$$
and using the bound from Proposition~3 of \vac{} for $\ttinf{(\Dc - \Delc)X}$ (see, also, p.22), we get
\begin{align}
    \ttinf{(I - UU\tr){(\Dc - \Delc)}{U \hQ\hL\inv}} 
    &\le \ttinf{(\Dc - \Delc)X} \cdot \opnorm{Q\L^{-1/2}\hQ\hL\inv}\\
    &\lesssim c_2\sqrt{n{\log{n}}} \cdot \frac{\kappa^{3}}{n^{3/2}} =: c''_1\frac{\sqrt{\log{n}}}{n}.
\end{align}
Plugging in these bounds back into $\ttinf{\hU - U\hQ}$ we get
\begin{align}
    \ttinf{\zeta^{(1)}} = \ttinf{(\Delc - \Dc)(\hU - U\hQ)\hL^{-1/2}} \lesssim \Msigma\sqrt{n} \cdot c''_1 \frac{\sqrt{\log{n}}}{n} \cdot \frac{\kappa}{\sqrt{n}} =: c^{(1)} \cdot\frac{\sqrt{\log{n}}}{n}.\label{eq:zeta1}
\end{align}

$\bullet$ For $\zeta^{(4)}$, for $J\Eps X = \onev \onev\tr\Eps\X$, we have
\begin{align}
    \ttinf{\zeta^{(4)}} \le \frac{1}{n}\ttinf{\onev\onev\tr\Eps X}\opnorm{(X\tr\X)\inv} \lesssim \frac{1}{n} \cdot \norm{\onev\tr\Eps X} \cdot \frac{\kappa^2}{n},\label{eq:zeta4-1}
\end{align}
where we used $\ttinf{\onev v^T} = \max_{i \in [n]}\norm{v} = \norm{v}$ for any $v \in \Rp$. Note that
\begin{align}
    \onev\tr\Eps X = \sum_{i, j \in [n]}\eps_{ij}x_j = \sum_{i < j}\eps_{ij}(x_i + x_j) \in \R^p,
\end{align}
is the sum of independent sub-exponential vectors with 
\begin{align}
    \norm{x_i + x_j} \le 2\Rx \qq{and} \Norm\big{\norm{\eps_{ij}(x_i + x_j)}_{_2}}_{\psi_1} \le 2\Rx\Msigma.
\end{align}
A straightforward application of \cref{prop:matrix-bernstein}; see, also, the proof of \ref{bound-9}, gives: with probability at least $1 - O(n^{-2})$,
\begin{align}
    \norm{\onev\tr\Eps X} \lesssim \Msigma\Rx n\sqrt{\log{n}},
\end{align}
and, plugging this back into \cref{eq:zeta4-1},
\begin{align}
    \ttinf{\zeta^{(4)}} \lesssim \frac{\Msigma\Rx n\sqrt{\log{n}}}{n^2} =: c^{(4)}\frac{\sqrt{\log{n}}}{n}.\label{eq:zeta4}
\end{align}

Combining \cref{eq:intermediate-quantities-1,eq:zeta3,eq:zeta1,eq:zeta4}, it follows that with probability greater than $1 - O(n^{-2})$,
\begin{align}
    \ttinf{\zeta} = \ttinf{\zeta^{(1)} + \zeta^{(2)} + \zeta^{(3)} + \zeta^{(4)}} \lesssim c' \frac{\sqrt{\log{n}}}{n},
\end{align}
for $c' = \max\qty{ c^{(1)}, c'_2, c^{(3)}, c^{(4)} }$.
\qed



\end{document}